\documentclass[10pt]{amsart}
\usepackage{amssymb,amsmath,latexsym,amscd,amsfonts,yfonts}
\usepackage{tikz}
    \usetikzlibrary{calc}

\usepackage{graphics}
\usepackage{graphicx}
\usepackage{wasysym}
\usepackage{booktabs}
\usepackage{appendix}

\def\ignore #1 {}

\newtheorem{thm}{Theorem}
\newtheorem{lem}[thm]{Lemma}

\newtheorem{defn}[thm]{Definition}
\newtheorem{prop}[thm]{Proposition}
\newtheorem{cor}[thm]{Corollary}
\newtheorem{question}{Question}

\theoremstyle{definition}
\newtheorem*{remark}{Remark}

\newtheorem{example}{Example}
\newtheorem{nug}{Nugget}

\newtheorem*{conj*}{Conjecture}

\newtheorem{innermainthm}{Main Theorem}
\newenvironment{mainthm}[1]
    {\renewcommand\theinnermainthm{#1}\innermainthm}
    {\endinnermainthm}

\def\hpic #1 #2 {\mbox{$\begin{array}[c]{l} \epsfig{file=#1,height=#2} \end{array}$}}
\def\vpic #1 #2 {\mbox{$\begin{array}[c]{l} \epsfig{file=#1,width=#2} \end{array}$}}

\def\R{\mbox{{\bf R}}}
\def\C{\mbox{{\bf C}}}
\newcommand{\CC}{{\mathcal{C}}}
\def\P{\mbox{{\bf P}}}

\def\encyc{\mbox{$\mathcal{E}$}}
\def\redencyc{\mbox{$\tilde{\mathcal{E}}$}}
\def\qstar{q^{\star}}

\newcommand{\bdry}{\partial}

\newcommand{\tot}{\sigma}
\DeclareMathOperator{\im}{Im}
\DeclareMathOperator{\In}{Inside}
\DeclareMathOperator{\Aff}{Aff}
\DeclareMathOperator{\area}{Area}
\DeclareMathOperator{\Vxs}{Vxs}

\newcommand{\D}{\Delta}

\newcommand{\s}{\sigma}

\newcommand{\T}{{\mathcal{T}}}

\def \sb #1 {\overline{\s}_{#1}}

\def\PP{\mathbf{p}}
\def\QQ{\mathbf{q}}
\def\RR{\mathbf{r}}
\def\SS{\mathbf{s}}

\newif\ifcorners
\newif\ifareas
\newif\ifinterior
\newif\ifcolors
\newif\ifnewareas

\def\poscolor{blue}
\def\negcolor{orange}

\def\colorone{black}
\def\colortwo{white}

\newcommand{\ace}{\begin{tikzpicture}
        \draw (0,0) -- (.2,0);
        \draw[fill=\colorone] (0,0) circle (1pt);
        \draw[fill=\colorone](.2,0) circle (1pt);
        \draw[fill=\colorone](.1,.2) circle (1pt);
        \end{tikzpicture}
        }
\newcommand{\threeforbidden}{\begin{tikzpicture}
        \draw (0,0) -- (.2,0);
        \draw (.1,.2) -- (.2,0);
        \draw (0,0) -- (.1,.2);
        \draw[fill=\colorone] (0,0) circle (1pt);
        \draw[fill=\colorone](.2,0) circle (1pt);
        \draw[fill=\colorone](.1,.2) circle (1pt);
        \end{tikzpicture}
        }
\newcommand{\acesub}{\begin{tikzpicture}
        \draw (0,.2) -- (0,0);
        \draw[fill=\colorone] (0,0) circle (1pt);
        \draw[fill=\colorone](.2,0) circle (1pt);
        \draw[fill=\colorone] (0,.2) circle (1pt);
        \draw[fill=\colorone](.2,.2) circle (1pt);
        \end{tikzpicture}
        }
\newcommand{\aceothersub}{\begin{tikzpicture}
        \draw (0,.2) -- (0,0) -- (.2,0);
        \draw[fill=\colorone] (0,0) circle (1pt);
        \draw[fill=\colorone](.2,0) circle (1pt);
        \draw[fill=\colorone] (0,.2) circle (1pt);
        \draw[fill=\colorone](.2,.2) circle (1pt);
        \end{tikzpicture}
        }
\newcommand{\fourone}{\begin{tikzpicture}
        \draw (0,.2) -- (0,0) -- (.2,0) -- (.2,.2);
        \draw[fill=\colorone] (0,0) circle (1pt);
        \draw[fill=\colorone](.2,0) circle (1pt);
        \draw[fill=\colorone] (0,.2) circle (1pt);
        \draw[fill=\colorone](.2,.2) circle (1pt);
        \end{tikzpicture}
        }
\newcommand{\fourtwo}{\begin{tikzpicture}
        \draw (0,0) -- (.2,0);
        \draw[fill=\colorone] (0,0) circle (1pt);
        \draw[fill=\colorone](.2,0) circle (1pt);
        \draw[fill=\colorone] (0,.2) circle (1pt);
        \draw[fill=\colortwo](.2,.2) circle (1pt);
        \end{tikzpicture}
        }
\newcommand{\six}{\begin{tikzpicture}
        \draw (0,.2) -- (.2,.2);
        \draw (0,0) -- (.2,0);
        \draw[fill=\colorone] (0,0) circle (1pt);
        \draw[fill=\colorone] (.2,0) circle (1pt);
        \draw[fill=\colorone] (.4,0) circle (1pt);
        \draw[fill=\colortwo] (0,.2) circle (1pt);
        \draw[fill=\colortwo](.2,.2) circle (1pt);
        \draw[fill=\colortwo](.4,.2) circle (1pt);
        \end{tikzpicture}
        }

\usepackage{microtype} % Improves spacing

\usepackage{scrextend}

\newcommand{\volume}[1]{
    \noindent\hrule
    \bigskip
    \begin{center}
        {\usefont{OT1}{ppl}{b}{n}
        {\fontsize{16pt}{0}\selectfont Volume #1 \\}
        }
    \end{center}
    \bigskip
    }
\newcommand{\degree}[1]{
    \begin{center}
        {\usefont{OT1}{ppl}{b}{n}
        {\fontsize{14pt}{0}\selectfont #1}\\
        }
    \end{center}
    \bigskip
    }
\newcommand{\entry}[4]{ \noindent $\mathbf{#1}$
    \hfill #2
    \begin{addmargin}[.25in]{.25in}
        {#3} \par
        {#4}
    \end{addmargin}
    \par
    \bigskip
    }

\begin{document}
%%%%%%%%%%%%%% TITLE STUFF %%%%%%%%%%%%%%%%%%%%
\title{An illustrated encyclopedia of area relations}
\author{Aaron Abrams and James Pommersheim}
%\date{\today}

\begin{abstract}
    To any combinatorial triangulation $T$ of a square, there is an associated polynomial relation $p_T$ among the areas of the triangles of $T$.  With the goal of understanding this polynomial, we consider polynomials obtained from $p_T$ by choosing $l$ of its variables and specializing $p_T$ to these variables by zeroing out the remaining variables.  We show that for fixed $l$, the set $\encyc_l$ of integer polynomials that appear as irreducible factors of such specializations is finite.  We compute this {\it area encyclopedia} $\encyc_l$ for $l\leq 4$.  We also show that in any dissection of a square into $l$ triangles, the areas of the triangles must satisfy a polynomial in $\encyc_l$.  Our results are obtained by studying the rational map that associates to each drawing of $T$ the tuple of areas of the triangles in that drawing.  By analyzing the ways of approaching the base locus, we derive restrictions on points of the closure of the image of this map. 
\end{abstract}
\maketitle

\setcounter{tocdepth}{1}
\tableofcontents

\section{Introduction and main results}
    \subsection{A polynomial relation among areas}
        Consider a triangulation of the unit square into $n$ triangles with areas $a_1,\dots, a_n$. Inspired by a theorem of Monsky \cite{monsky}, we proved in \cite{triangles1} the existence of a homogeneous integer polynomial in $n$ variables  that vanishes at the point $(a_1,\dots,a_n)$. By construction, the polynomial depends only on the combinatorics of the triangulation, so the same polynomial also vanishes at any tuple of areas arising by deforming the triangulation.  
        
        \begin{example}\label{ex:intro}  
            Let $T_1$ be the first triangulation shown in Figure \ref{fig:example1}.  If the boundary is drawn as a parallelogram, then regardless of where the interior vertex is drawn, the areas of the triangles satisfy the polynomial $$p_{T_1}=A-B+C-D,$$ where the variables correspond to the areas of the triangles in order around the square.  For the second triangulation $T_2$ in Figure \ref{fig:example1}, as long as the boundary is a parallelogram, the areas of the triangles satisfy the polynomial $$p_{T_2}=A^2-2AC+2AE+C^2+2CE+E^2-B^2-2BD-2BF-D^2+2DF-F^2.$$ 
            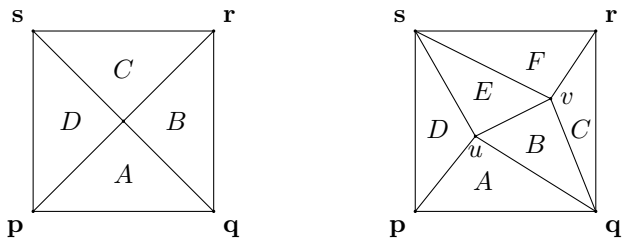
\begin{figure}[h]\label{fig:example1}
                \centering
                \cornerstrue
                \areastrue  
    \begin{tikzpicture}[scale=.4]
    \coordinate (P) at (0,0);
\coordinate (Q) at (6,0);
\coordinate (R) at (6,6);
\coordinate (S) at (0,6);

\draw[fill=black] (P) circle (1pt);
\draw[fill=black] (Q) circle (1pt);
\draw[fill=black] (R) circle (1pt);
\draw[fill=black] (S) circle (1pt);
\draw (P) -- (Q) -- (R) -- (S) -- cycle;

\coordinate (x) at (3,3);
\draw[fill=black] (x) circle (1pt);

\draw (P) -- (x);
\draw (Q) -- (x);
\draw (R) -- (x);
\draw (S) -- (x);

\ifcorners
    \draw (P) node[anchor=north east]{$\PP$};
    \draw (Q) node[anchor=north west]{$\QQ$};
    \draw (R) node[anchor=south west]{$\RR$};
    \draw (S) node[anchor=south east]{$\SS$};
\fi

\ifcolors
    \draw[line width=0, fill=\poscolor] (P) -- (Q) -- (x) -- cycle;
    \draw[line width=0, fill=\poscolor] (R) -- (S) -- (x) -- cycle;
    \draw[line width=0, fill=\negcolor] (Q) -- (R) -- (x) -- cycle;
    \draw[line width=0, fill=\negcolor] (S) -- (P) -- (x) -- cycle;
\fi

\ifareas
    \draw (3,1.25) node{$A$};
    \draw (4.75,3) node{$B$};
    \draw (3,4.75) node{$C$};
    \draw (1.25,3) node{$D$};
\fi

\ifnewareas
    \draw (3,1.25) node{\scriptsize $D$};
    \draw (4.75,3) node{\scriptsize $B$};
    \draw (3,4.75) node{\scriptsize $C$};
    \draw (1.25,3) node{\scriptsize $A$};
\fi
    \end{tikzpicture}
     \areasfalse \qquad\qquad \interiortrue\areastrue  
    \begin{tikzpicture}[scale=.4]
    \coordinate (P) at (0,0);
\coordinate (Q) at (6,0);
\coordinate (R) at (6,6);
\coordinate (S) at (0,6);

\draw[fill=black] (P) circle (1pt);
\draw[fill=black] (Q) circle (1pt);
\draw[fill=black] (R) circle (1pt);
\draw[fill=black] (S) circle (1pt);
\draw (P) -- (Q) -- (R) -- (S) -- cycle;

\coordinate (u) at (2,2.5);
\coordinate (v) at (4.5,3.75);
\draw[fill=black] (u) circle (1pt);
\draw[fill=black] (v) circle (1pt);

\draw (P) -- (u) -- (v) -- (R);
\draw (Q) -- (u) -- (S);
\draw (S) -- (v) -- (Q);

\ifcorners
    \draw (P) node[anchor=north east]{$\PP$};
    \draw (Q) node[anchor=north west]{$\QQ$};
    \draw (R) node[anchor=south west]{$\RR$};
    \draw (S) node[anchor=south east]{$\SS$};
\fi

\ifareas
    \draw (2.25,1) node{$A$};
    \draw (4,2.25) node{$B$};
    \draw (5.5,2.75) node{$C$};
    \draw (.75,2.75) node{$D$};
    \draw (2.25,4) node{$E$};
    \draw (4,5) node{$F$};
\fi

\ifinterior
    \draw (u) node[anchor=north]{$u$};
    \draw (v) node[anchor=west]{$v$};
\fi

\ifcolors
    \draw[line width=0, fill=\poscolor] (P) -- (Q) -- (u) -- cycle;
    \draw[line width=0, fill=\poscolor] (v) -- (Q) -- (R) -- cycle;
    \draw[line width=0, fill=\poscolor] (u) -- (v) -- (S) -- cycle;
    \draw[line width=0, fill=\negcolor] (u) -- (Q) -- (v) -- cycle;
    \draw[line width=0, fill=\negcolor] (v) -- (R) -- (S) -- cycle;
    \draw[line width=0, fill=\negcolor] (u) -- (S) -- (P) -- cycle;
\fi
    \end{tikzpicture}
     \interiorfalse\areasfalse\\
                \cornersfalse
                \caption{The areas of the four triangles on the left satisfy a homogeneous linear relation.  On the right, the areas of the six triangles satisfy a homogeneous quadratic polynomial.}
            \end{figure}
        \end{example}
        
        In this manner, to each combinatorial triangulation $T$ of a $4$-gon, there is associated an irreducible integer polynomial $p_T$, well defined up to sign, with variables corresponding to the triangles of $T$.  The polynomial $p_T$ has the property that for any piecewise linear map of $T$ to the plane such that the boundary forms a square (or more generally any parallelogram), $p_T$ vanishes when one replaces each variable by the area of the corresponding triangle.  This polynomial $p_T$ is known as the {\it area polynomial} of the combinatorial triangulation $T$.

        In previous work \cite{triangles1,chapter3} we have established various facts about this polynomial.  For instance, regarding the degree we have characterized when $p_T$ is linear and proved that its degree can be quite large.  However we still do not have a satisfactory way of understanding the degree of $p_T$ in general.  Regarding the coefficients, we used Monsky's theorem to show that the mod $2$ reduction of $p_T$ is equal to a power of the sum of the variables of $p_T$, and we have conjectured a kind of positivity which relates to the size of the coefficients of $p_T$ (see \cite{chapter3}).  Many mysteries about $p_T$ remain. The purpose of the current work is to further the study of this polynomial.
        
        Of particular interest here are polynomials obtained from $p_T$ by substituting 0 for some subset of the variables of $p_T$. We call such polynomials {\it specializations} of $p_T$.  Specializing the polynomial $p_T$ is natural from a geometric point of view; indeed, setting a variable of $p_T$ equal to $0$ reflects the geometric requirement that the corresponding triangle be degenerate.  Our notation for these specializations is $(p_T)|_L$, where $L$ is the set of variables that are not set to zero.
         
        \begin{example}\label{ex:specialize}
            Let $T$ be the triangulation into 6 triangles shown in Example \ref{ex:intro}.  Then $p_T|_{\{A,B\}}=(A+B)(A-B)$ and $p_T|_{\{A,C\}}=(A-C)^2$, etc.  The $3$-variable specializations include $(A+B-C)(A-B-C)$ and $A^2-2AC+2AE+C^2+2CE+E^2$, among others.
        \end{example}
 
    \subsection{The encyclopedia of area relations}   
        The polynomial $p_T$ for arbitrary $T$ seems quite complicated, but it turns out that the specializations of these polynomials to a fixed number of variables are severely restricted.  In fact, we show that for fixed $l$, there are only finitely many polynomials that can show up as an irreducible factor of an $l$-variable specialization of any $p_T$.  
     
        \begin{mainthm}{1}\label{mainthm:encyc} (See Theorem \ref{thm:encyc}.)
            \emph{
            For each positive integer $l$, there is a computable finite list $\encyc_l$ of integer polynomials in $l$ variables such that for any triangulation $T$, if $q$ is an irreducible factor of any $l$-variable specialization of $p_T$, then then $q$ is in $\encyc_l$, up to scaling and renaming of variables.
            }
        \end{mainthm}
        
        The union $\encyc=\bigcup\limits_{l=1}^{\infty} \encyc_l$ of these lists we call the {\it area encyclopedia}.  The set $\encyc_l$ is the \emph{$l^{\rm{th}}$ volume} of $\encyc$.
        \begin{example}\label{ex:introencyc}
            We will show that the first few volumes of the area encyclopedia are 
            \begin{align*}
                \encyc_1 = \ \{\ &A\ \}\\
                \encyc_2 = \ \{\ &A +B,\  A - B\ \}\\
                \encyc_3 = \ \{\ &A+B+C,\ A+B-C,\  A^2+2AB+2AC+B^2-2BC+C^2\ \}
            \end{align*} 
            The fourth volume $\encyc_4$ consists of $3$ linear polynomials, $4$ quadratics, and a quartic.  See Section \ref{sec:theencyc}.
        \end{example}
        
        The encyclopedia has many interesting and previously unproved consequences for the polynomial $p_T$.  For instance, from the $l=2$ case, it follows immediately that all the leading coefficients of $p_T$ are equal up to sign; that is, for any two variables $A$ and $B$, the coefficient of $A^d$ is equal to $\pm 1$ times the coefficient of $B^d$ (where $d$ is the degree of $p_T$).  Previously, we had established only that these coefficients must be odd.  This observation suggests a canonical $2$-coloring of the triangles in any triangulation $T$, based on the signs of these leading coefficients.  
        
        It turns out that for many of the polynomials $q$ of $\encyc$, there are two variables $X$ and $Y$ so that $q$ is a polynomial in $X+Y$ and the other variables. We call such polynomials {\it algebraic subdivisions}.  They are closely tied to geometric subdivisions of the triangulation.  By eliminating such algebraic subdivisions, we introduce an abridged version of the encyclopedia $\redencyc$, studied in Section \ref{sec:redencyc}.  
     
    \subsection{Consequences for dissections}  
        The area encyclopedia $\encyc$ helps us answer certain questions about the areas of triangles arising in a dissection of a square.  By a \emph{dissection} we mean a collection of triangles in $\R^2$ whose interiors are disjoint and whose union is the entire square.  Unlike in a triangulation, two triangles in a dissection are not required to intersect in a common face.  
        
        A celebrated theorem of Paul Monsky asserts that it is impossible to dissect a square into an odd number of triangles of equal area \cite{monsky}. A natural and motivating question to ask is the following:
        
        \begin{quote}
            For which tuples $(a_1, \dots, a_l)$ of positive real numbers is there a dissection of the unit square into $l$ triangles with areas $a_1, \dots, a_l$?
        \end{quote}
    
        By interpreting a dissection as the image of a triangulation under a PL map, we are able to show that this question is closely related to our area encyclopedia.
        \begin{mainthm}{2}\label{mainthm:dissections}(See Corollary \ref{cor:dissection}.)
            \emph{
            Given a dissection of a square into $l$ triangles with areas $a_1,\dots,a_l$, the $l$-tuple $(a_1,\dots,a_l)$ satisfies at least one polynomial in $\encyc_l$.
            }
        \end{mainthm}
        
        We also prove a partial converse.
        \begin{mainthm}{2.5}\label{mainthm:nonconverse} (See Corollary \ref{cor:EimpliesT}.)
            \emph{
            Let $(a_1,\dots,a_l)$ be an $l$-tuple of positive real numbers satisfying some polynomial $q\in\encyc_l$.  Then there is a triangulation $T$ with at most $3l-2$ triangles and a map from the vertices of $T$ into $\C^2$ such that the boundary of $T$ maps to a square and such that there are exactly $l$ nondegenerate triangles, the areas of which are $a_1,\dots,a_l$.
            }
        \end{mainthm}
        
        What distinguishes this from a complete converse to Main Theorem \ref{mainthm:dissections} is that the codomain of the promised mapping is $\C^2$, not $\R^2$.
        
        \begin{question}\label{question:realizerootasdiss}
            Given an $l$-tuple $(a_1,\dots,a_l)$ of positive real numbers satisfying some polynomial in the encyclopedia $\encyc_l$, must there exist a dissection of a square into triangles with areas $a_1,\dots,a_l$?   
        \end{question} 
        
        In trying to prove an affirmative answer to this question, we discovered an interesting counterexample to a slightly stronger statement. To wit, the vanishing of a particular $(p_T)|_L$ at the $l$-tuple $(a_1,\dots,a_l)$ of positive real numbers does not imply that the triangulation $T$ can be drawn in $\R^2$ with nonzero areas $a_1,\dots,a_l$.  This is demonstrated in Example \ref{ex:complexequi} by the triangulation $T=T_{2222}$ (see Table \ref{table:T.up.to.4.int}), which has $10$ triangles of which we have chosen a certain subset $L$ of size $4$. An irreducible quartic factor $q\in \encyc_4$ of the specialization $(p_T)|_L$ vanishes at the point $[1:1:1:1]$ of the area variety, and indeed there is a mapping of the vertices of $T$ to $\C^2$ whose only nondegenerate triangles are those of $L$ and have area $1$.  However, there is no such map to $\R^2$; in fact there is no map to $\R^2$ in which the four triangles in $L$ have any positive areas whatsoever (and the other six vanish).  In particular no drawing of $T$ gives a dissection realizing these areas.  However, this does not provide a negative answer to Question \ref{question:realizerootasdiss}, because the areas $1,1,1,1$ can easily be realized as the areas of the triangles in a dissection of a square.
        
    \subsection{Analysis of the area map}
        All the results in this paper stem from a detailed consideration of a certain rational map called the \emph{area map}. Given a combinatorial triangulation $T$ of a $4$-gon, the area map for $T$ associates to each drawing of $T$ in the plane the tuple of the areas of the triangles formed.  More precisely, we define the \emph{drawing space} $X(T)$ to be the set of all maps of the vertex set of $T$ to $\C^2$ such that the four boundary vertices of $T$ map to the vertices of a parallelogram, and we define the \emph{area space} $Y(T)$ to be the projective space with one coordinate for each triangle of $T$.  For any drawing $\rho\in X(T)$ we can measure the (complex-valued) area of each triangle of $T$; these areas are quadratic polynomials in the coordinates of $\rho$.  In this way the $i$th coordinate function of the rational map 
        $$\area_T:X(T)\dashrightarrow Y(T)$$
        records the area of the $i$th triangle in the drawing $\rho$.  In \cite{triangles1}, we showed that for any $T$, the closure of image of the map $\area_T$ is an irreducible hypersurface in $Y(T)$.  This is the \emph{area variety}, denoted $V(T)$.  The area polynomial $p_T$ is the defining polynomial of the hypersurface $V(T)$.

        Points in the image of the area map are represented by actual drawings $\rho\in X(T)$, which makes them relatively easy to study.  The challenge in this work is to understand points of $V$ that are not in the image of $\area$.  Such points are limits of areas of drawings obtained by approaching the base locus of the area map.  After a careful analysis of the ways one can approach the base locus, we arrive at the following theorem, which puts a significant restriction on points of $V(T)$ that are not in the image of $\area$.
        
        \begin{mainthm}{3}\label{mainthm:points}(See Theorem \ref{thm:characterization}.)
            \emph{
            Let $T$ be a triangulation, and let $w\in V(T)$.  Then either $w\in\im\area$ or there is a subset of the coordinates of $w$, not all of which are $0$, that sums to $0$.
            }
        \end{mainthm}
        
        This is the main theorem that we use to derive our conclusions about the polynomial $p_T$ and the encyclopedia. A richer version of this theorem appears in the body of this paper as Theorem \ref{thm:characterization}.
       
        Central to the proof of Main Theorem \ref{mainthm:points} is the idea of a  {\it bubble}.  For some drawings $\rho\in X(T)$ there is a cycle $C$ in the $1$-skeleton of $T$ with the property that all vertices of $C$ map to the same point $\star\in\C^2$ under $\rho$, but not all vertices inside $C$ map to $\star$.  In this case, we say that $\rho$ has a bubble.  It turns out that whenever we have a point of $w\in V(T)$ that is not in the image of $\area$, there must lurk a bubble.  More precisely, any point $w\in V(T)\setminus\im\area$ can be approached along a curve in $V(T)$ that is the image (under $\area$) of a curve in $X(T)$ that approaches a drawing $\rho\in X(T)$ that contains a bubble. We then show that the sum of the coordinates of $w$ corresponding to the triangles inside the bubble must vanish.  This puts us in position to deduce Main Theorem \ref{mainthm:points}.
   
    \subsection{Illustrating the encyclopedia}

        Main Theorem \ref{mainthm:points} has an important consequence that allows us to picture the polynomials in our encyclopedia. 
 
        Suppose $q\in\encyc_l$. Take a triangulation $T$ and a subset $L$ of the triangles of $T$ such that $q=(p_T)|_L$.  Denote by $Z_L\subset Y(T)$  the intersection of all coordinate hyperplanes corresponding to triangles not in $L$.  The irreducible polynomial $q$ defines a component $W$ of the intersection $V_L=V(T)\cap Z_L$.  A consequence of Main Theorem \ref{mainthm:points} is that a generic point of $W$ is in the image of $\area$.  Hence there is a family of drawings of $T$ whose areas realize the generic points of $W$.  In this way, we can draw pictures to illustrate each polynomial of the encyclopedia.
 
        \begin{example}
            The polynomial $q=A-B+C$ is in $\encyc_3$ by dint of being a specialization of the polynomial $p_T$ where $T$ is the triangulation from Figure \ref{fig:example1} (left).  Letting $D=0$, we imagine pictures in which the boundary vertices $\PP\QQ\RR\SS$ form a parallelogram and the interior vertex is constrained to lie on $\PP\SS$, as in Figure \ref{fig:2} (left).  This picture then illustrates the polynomial $q=A-B+C$ in the sense that for a generic zero $(a,b,c)$ of $q$, there is a drawing of this type with areas $(a, b, c)$.

            \begin{figure}\label{fig:introillustrate}
                \centering
                \cornerstrue
                \interiortrue
                 
    \begin{tikzpicture}[scale=.4]
    \coordinate (P) at (0,0);
\coordinate (Q) at (6,0);
\coordinate (R) at (6,6);
\coordinate (S) at (0,6);

\draw[fill=black] (P) circle (1pt);
\draw[fill=black] (Q) circle (1pt);
\draw[fill=black] (R) circle (1pt);
\draw[fill=black] (S) circle (1pt);
\draw (P) -- (Q) -- (R) -- (S) -- cycle;

\coordinate (x) at (0,3);
\draw[fill=black] (x) circle (1pt);

\draw (Q) -- (x) -- (R);

\ifcorners
    \draw (P) [anchor=north east] node{$\PP$} ;
    \draw (Q) [anchor=north west] node{$\QQ$} ;
    \draw (R) [anchor=south west] node{$\RR$} ;
    \draw (S) [anchor=south east] node{$\SS$} ;
\fi

\ifinterior
    \draw (x) [anchor=east] node{$v$} ;
\fi
    \end{tikzpicture}
     \qquad \qquad  
    \begin{tikzpicture}[scale=.4]
    \coordinate (P) at (0,0);
\coordinate (Q) at (6,0);
\coordinate (R) at (6,6);
\coordinate (S) at (0,6);

\draw[fill=black] (P) circle (1pt);
\draw[fill=black] (Q) circle (1pt);
\draw[fill=black] (R) circle (1pt);
\draw[fill=black] (S) circle (1pt);
\draw (P) -- (Q) -- (R) -- (S) -- cycle;

\coordinate (u) at (10,0);
\coordinate (v) at (6,2.4);
% \coordinate (v) at (-4,6);
% \coordinate (u) at (0,3.6);
\draw[fill=black] (u) circle (1pt);
\draw[fill=black] (v) circle (1pt);

\draw (P) -- (u) -- (v) -- (R);
\draw (Q) -- (u) -- (S) -- (v) -- cycle;

\ifcorners
    \draw (P) node[anchor=north east]{$\PP$};
    \draw (Q) node[anchor=north]{$\QQ$};
    \draw (R) node[anchor=south west]{$\RR$};
    \draw (S) node[anchor=south east]{$\SS$};
\fi

\ifinterior
    \draw (u) node[anchor=north west]{$u$};
    \draw (v) node[anchor=south west]{$v$};
\fi
    \end{tikzpicture}
     %\qquad \qquad \fig{.15}{D2ACE-alt}
                \interiorfalse
                \cornersfalse
                \caption{Illustrations of two of the polynomials from the encyclopedia.}
                \label{fig:2}
            \end{figure}
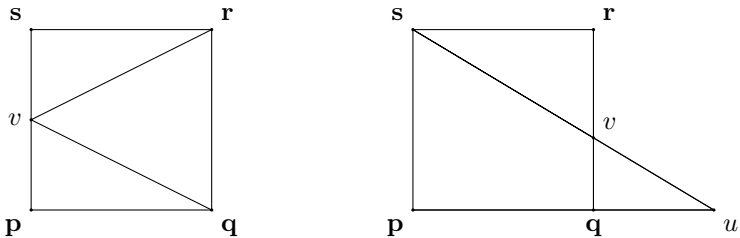
        
            In a similar manner the (unique) quadratic polynomial $q\in\encyc_3$ is illustrated by the picture in Figure \ref{fig:2} (right). After renaming the variables, the polynomial $q=A^2-2AC+2AE+C^2+2CE+E^2$ is the specialization of $p_{T_2}$ to the variables $A, C, E$, where $T_2$ is shown in Figure \ref{fig:example1} (right).  These variables represent the areas of the triangles $\PP\QQ u$, $\QQ\RR v$, and $\SS uv$.  Imagining that the other triangles are degenerate, one is led to the picture shown in Figure \ref{fig:2} (right). Here the vertex $u$ is constrained to lie on the line $\PP\QQ$ and $v$ is constrained to lie on both lines $\QQ \RR$ and $\SS u$.  One can verify that the areas $A,C,E$ of the three triangles in this drawing (one of which is negative) satisfy the polynomial $q$.  
        \end{example}
 
        In the preceding example, we considered triangulations $T$ with additional collinearity constraints imposed on certain subsets of the vertices.  These {\it constrained triangulations}, introduced in \cite{triangles1}, provide illustrations for every polynomial in the encyclopedia. 
        
        \begin{mainthm}{4}\label{mainthm:illustrate} (See Theorem \ref{thm:illustration}.)
            \emph{
            Each $q\in\encyc_l$ is illustrated by a constrained triangulation with at most $3l-2$ triangles.
            }
        \end{mainthm}
        
        We conclude this paper about the encyclopedia with (the beginning of) the encyclopedia itself.  For each polynomial, we give drawings of constrained triangulations that illustrate the polynomial.  We also give an abridged version of the encyclopedia from which the full encyclopedia can be reconstructed.
 
        This paper is organized in three parts. In part 1 we conduct our analysis of the base locus of the area map.  In part 2 we will apply this work to the area polynomials, proving our finiteness results.  In part 3 we give the first four volumes of the full illustrated encyclopedia as well as the abridged version.

\part{The area variety}
    
    A combinatorial triangulation $T$ can be drawn in the plane in a variety of ways.  In any given drawing, the areas of the triangles can be measured.  That is, there is an \emph{area map} from the space of drawings to the space of possible areas.  The closure of the image of this map is called the \emph{area variety} of $T$.  It is a hypersurface and therefore the zero set of a single polynomial, well defined up to scalar multiplication, called the \emph{area polynomial} of $T$.
    
    In the first part of this paper we study the area map and the area variety.  We are able to draw conclusions about area polynomials only after understanding the points of this variety. Of particular importance are the points that are not in the image of the area map.  After a careful technical analysis of the base locus of the area map, our efforts to understand such points culminates in Theorem \ref{thm:characterization}, which puts a significant restriction on which points can be in $V\setminus \im \area$.  
    
    This allows us to draw conclusions about the polynomial $p$. For example, one basic consequence of this analysis is that the point $[1:0:\cdots:0]$ is not in $V$, giving another proof (there being one in \cite{chapter3}, see below) that the each leading term $A^d$ appears in the polynomial $p$ with nonzero coefficient (where $d$ denotes the degree of $p$).  Additional consequences, including our main theorems, are fleshed out in the second part of the paper.

    \section{Background}
        First, here are some definitions.
    
        \begin{defn}[Triangulation, corner, interior vertex] \label{def:honest} 
            A \emph{triangulation} $T$ (of a square) is an oriented simplicial complex homeomorphic to a disk with four vertices on the boundary.  
            Vertices on the boundary of $T$ are called \emph{corners} and are labeled $\PP,\QQ,\RR,\SS$ in the cyclic order determined by the orientation of $T$.  Other vertices of $T$ are called \emph{interior vertices.}  
        \end{defn}
        
        We will be drawing triangulations in the plane.  It will be helpful to work over $\C$, so we use the ``complex plane'' $\C^2$.
        
        \begin{defn}[Drawing, nondegenerate triangle]\label{def:honestdrawings} 
            Let $T$ be a triangulation.  A \emph{drawing} of $T$ is a map $\rho:\Vxs(T)\to\C^2$ such that $\rho(\PP)+\rho(\RR)=\rho(\QQ)+\rho(\SS)$.  That is, the image of the boundary is a (possibly degenerate) parallelogram.
            The space of all drawings, topologized as a subspace of $(\C^2)^{\Vxs(T)}$, is denoted $X(T)$.
            
            The \emph{nondegenerate} triangles of $\rho$ are the triangles $\Delta$ of $T$ with the property that $\rho$ maps the vertices of $\Delta$ to three non-collinear points of $\C^2$.  Other triangles are \emph{degenerate} triangles of $\rho$.
        \end{defn}
        
        Note that $X(T)$ is isomorphic to the affine space $(\C^2)^{v-1}$ where $v$ is the number of vertices of $T$.  This is because we may specify $\rho$ arbitrarily on all vertices of $T$ except for one corner, whose image will then be determined by the other corners.
        
        Next we define the area map and the area variety.
            
        Let $p_i=(x_i,y_i)$ for $i=1,2,3$ be three points in the affine plane $\C^2$.  We define the area of the (ordered) triangle $\Delta=(p_1,p_2,p_3)$ to be 
        $$\area(\Delta)=\area(p_1p_2p_3)=\frac 12
        \left|
        \begin{matrix}
        1 & 1 & 1 \\
        x_1 & x_2 & x_3 \\
        y_1 & y_2 & y_3 
        \end{matrix}
        \right|.
        $$
        Note that if $\{p_i\}\subset\R^2$ then this is the usual signed area function.
        Note also that $\area(p_1p_2p_3)=0$ if and only if $p_1,p_2,p_3$ lie on a (complex) line in $\C^2$.
        
        Let $T$ be a fixed triangulation, and let $\Delta_1,\dots,\Delta_n$ be the triangles of $T$. 
        Note each $\Delta_i$ inherits an orientation from $T$.  
        Let $Y(T)$ be the projective space $\P^{n-1}$ with homogeneous coordinates $[\cdots:A_i:\cdots]$, $1\le i\le n$.
        
        \begin{defn}[Area map, area variety]\label{def:areamaphonest} 
            The rational map 
            $$\area=\area_T:X(T)\dashrightarrow Y(T)$$
            given by
            $\area(\rho)=[\cdots:\area(\rho(\Delta_i)):\cdots]$ is called the \emph{area map} associated to $T$.
            The \emph{area variety} $V=V(T)$ is the closure in $Y(T)$ of $\area(X(T))$.
        \end{defn}

        One should pause here to think through a dimension count.  The affine group $\Aff_2(\C)$ acts on $X$ and the area map is equivariant with respect to this action.
        Thus the generic fibers of the area map are at least 6-dimensional (that being the size of $\Aff_2$).
        Independently one easily counts that $\dim(Y)=n-1=\dim X-5$.
        Therefore, for a given $T$ the area variety $V(T)$ has codimension at least 1 in $Y(T)$, with equality if and only if generic fibers are exactly 6-dimensional.  
        
        \begin{thm}\cite{triangles1}\label{thm:honest}
            For any triangulation $T$, the area variety $V(T)$ is a hypersurface in $Y(T)$.
        \end{thm}
        
        In other words, almost every drawing is \emph{area-rigid}:  there is no 1-parameter family of area-preserving deformations of the drawing, other than those contained in an $\Aff_2$ orbit.  It is straightforward to prove Theorem \ref{thm:honest} inductively; in \cite{triangles1} a version of Theorem \ref{thm:honest} is proved in the more general and delicate context of \emph{constrained triangulations}, which are discussed in the present paper in Section \ref{sec:illustrating}.  
        
        A consequence of Theorem \ref{thm:honest} is that there is a unique (up to scaling) nonzero homogeneous polynomial $p=p_T$ that vanishes on $V$.  
        The polynomial $p$ is irreducible because $X$, and therefore $V$, is an irreducible variety. 
        Also $p$ has rational coefficients (because the coordinate functions of $\area$ do) so $p$ can be normalized to have integer coordinates with no common factor.  
        We assume this has been done; the polynomial $p$ is now well-defined up to sign for any $T$.
        
        \begin{defn}[Area polynomials]
            We call $\{\pm p\}$ the \emph{area polynomials} for $T$.
        \end{defn}
        
        We remark that the area polynomials are computable, e.g.~using Gr\"obner basis techniques, but that these computations quickly become intractable as the triangulation grows.
        
        In \cite{chapter3}, using a generalization of Monsky's theorem, we described the mod 2 reduction of $\pm p$.  Let $\sigma$ denote the linear form which is the sum of the variables.
        
        \begin{thm}[\cite{chapter3}, Mod 2 Theorem]\label{thm:mod2honest}
            The area polynomial $p$ satisfies $$p \equiv \sigma^d \mod 2,$$ where $d=\deg p$.
        \end{thm}
        
        \begin{cor}[\cite{chapter3}]
            Suppose the area polynomial $p$ of the triangulation $T$ has degree $d$.  Then all leading terms $A_i^d$ occur with odd (hence nonzero) coefficient. 
        \end{cor}
        
        In the present paper we show that the leading coefficients are in fact all equal, up to sign.  We suspect these coefficients are all $\pm 1$, as we conjectured already in \cite{chapter3}.

    \section{Base locus:  bubbles and bursting}\label{sec:baselocus}
        
        With the goal of understanding the polynomial $p=p_T$, we now embark on our study of the variety $V=V(T)$, which is the closure of the image of the map $\area=\area_T$.  Points of $\im\area$ are easy to understand, since they represent the actual areas of drawings of $T$. Points of $V$ that are not contained in $\im\area$ are not so readily pictured, as they are obtained only as limits of points in the image.  Our efforts to understand such points culminates in Theorem \ref{thm:characterization}, which puts a significant restriction on which points can be in $V\setminus \im \area$.  
        
        \subsection{Basic picture}\label{subsec:basic}
        
            \def\rhom{\rho_m}
            \def\rhoseq{ \{\rhom\}    }
            
            When contemplating a limit of drawings, one imagines a path $\gamma(s)$ in $X$ for $s\in(-\epsilon, \epsilon)$ that approaches a limiting drawing $\rho=\gamma(0)$ as $s\to 0$.  If there are triangles drawn by $\rho$ with nonzero area, then $\rho$ is in the domain of $\area$ and the limit of $\area(\gamma(s))$ is $\area(\rho)$, which is manifestly in the image of the area map.  Thus to see points of $V$ that are not in the image, we imagine that our limiting drawing $\rho$ consists entirely of triangles of area $0$.   That is, $\rho$ is in the base locus\footnote{The base locus is clearly contained in $L$; except in very trivial cases, $L$ is also contained in the base locus.}
            \begin{equation*}
                L = \{ \ \rho\in X\  | \  \area_i(\rho) = 0\ {\rm for\ all} \ i  \}.
            \end{equation*}
            where $\area_i$ denotes the area of the $i$th triangle. 
            Thus $\area(\rho)$ may not be defined, but if $w=\lim_{s\to 0} \area(\gamma(s))$ exists, then $w\in V$ even though $w$ may not be in the image of $\area$.  Indeed, all points of $V$ are obtained via such paths $\gamma(s)$, which we call {\it generating paths}, defined below.  
            
            We are particularly interested in approaching the highly degenerate drawing $\rho$ with drawings $\gamma(s)$ in which the boundary is drawn as a nondegenerate parallelogram.  For this purpose, we introduce the notation  $\tot(\rho)=\sum \area_i(\rho)$ to denote the total area of a drawing $\rho\in X$.  This depends
            only on where $\rho$ maps the corners of the square. We also introduce the hyperplane 
            $$H=\{y\in Y : \sum y_i=0\} \quad \subset \quad Y.$$
            Note that if $\rho$ is any drawing in the domain of $\area$ such that $\tot(\rho)=0$, then $\area(\rho)\in H$.
            
            \begin{defn}[Generating path]
                Let $w\in V$.
                A \emph{generating path} for $w$ is an analytic path $$\gamma:(-\epsilon,\epsilon)\to X$$ for some $\epsilon>0$ such that
                \begin{enumerate}
                    \item $\rho=\gamma(0)$ is not a constant map, i.e., not all vertices of the drawing $\rho$ are drawn at the same point of $\C^2$;
                    \item for all $s\neq 0$, we have $\tot(\gamma(s))\neq 0$, and hence in particular $\gamma(s)$ is in the domain of $\area$ and $\area(\gamma(s))\not\in H$;
                    \item $\area(\gamma(s))\to w$ as $s\to 0$.
                \end{enumerate}
            \end{defn}

            \begin{lem}\label{lem:existgen}
                Every $w\in V$ has a generating path.
            \end{lem}
            
            \begin{proof}
                Since $w$ is in the closure of $\im\area$, an irreducible variety not contained in $H$, $w$ is also in the closure of $\im\area\setminus H$.  Thus $w$ is the limit  of a sequence  $w_m$ in  $\im\area\setminus H$.  Let $\rho_m$ be any preimage of $w_m$.  For each $m$, translate $\rho_m$ so that $\rhom(\PP) = (0,0)$ and scale so that the maximum length of any edge in $\rhom$ is $1$.  Since translating and scaling do not affect the image under $\area$, we still have $\area(\rho_m)=w_m$ for all $m$.   Now $\rho_m$ is a bounded sequence in $X$, so by passing to a subsequence, we may assume that $\rhoseq\to\rho$  with $\area(\rho_m)\to w$.  Property (1) follows from the fact $\rhom$ has an edge of length 1.  Also, $w_m\notin H$ implies that  $\tot(\rho_m)\neq 0$ for all $m$. Once we have this $\rho\in X$ we may choose an analytic path $\gamma(s)$ satisfying (1), (2), and (3).  
            \end{proof}
            
            \begin{lem}\label{key.lemma}
                Suppose that $w\in V\setminus H$.  Let $\gamma$ be a generating path for $w$.  Then
                $\frac{\area_i(\gamma(s))}{\tot(\gamma(s))}$ is bounded in a neighborhood of $s=0$. 
            \end{lem}
            
            \begin{proof}
                Since $\area(\gamma(s))$ and $w$ are not in $H$, we can  restrict our attention to the affine
                piece $\sum A_i=1$ of $Y$.  Here the coordinates of $\area(\gamma(s))$ have the
                form $\frac{\area_i(\gamma(s))}{\tot(\gamma(s))}$, which approaches $\frac{w_i}{\sum w_i}$.  The lemma follows.
            \end{proof}

            Lemma \ref{key.lemma} helps us toward our goal of identifying points of $V$ that are not in $\im\area$, because it often allows us to conclude that such a point must be in $H$.  One may have the na\"ive (and, it turns out, incorrect) intuition, as we once did, that $V\setminus\im\area \subset H$; our thinking was that a point $w\in V\setminus\im\area$ is the limit of a generating path $\gamma$ along which all of the triangles of $\gamma(s)$ are shrinking in area as $s\to 0$, and so the coordinates of the limit $w$ should sum to $0$.  The truth, of course, is that coordinates of $w$ record the relative rates at which the various areas are shrinking.  The following example displays a phenomenon called a bubble, which we define precisely in the next section.
            
            \begin{example}[A bubble]\label{ex:bubble} 
                Let $T$ be the triangulation shown in Figure \ref{fig:bubble} (left).  We coordinatize the space $Y(T)\cong\P^9$ so that the first two coordinates correspond to the triangles labeled $A_1$ and $A_2$ and the last three correspond to the triangles surrounding the vertex $\star$.
                We claim that for any three distinct complex numbers $a,b,c$, the point $$w=[1:1:0:0:0:0:0:a-b:b-c:c-a]$$
                of $Y(T)$ is in $V$ but is not in $\im\area\cup H$.
                
                    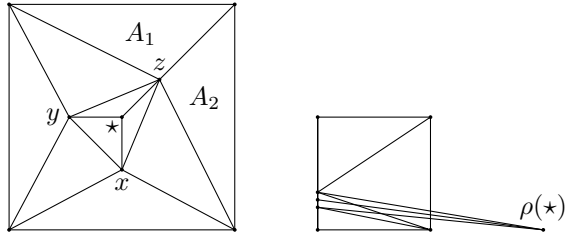
\begin{figure}
                        \centering
                         
    \begin{tikzpicture}[scale=.5]
    \coordinate (P) at (0,0);
\coordinate (Q) at (6,0);
\coordinate (R) at (6,6);
\coordinate (S) at (0,6);

\draw[fill=black] (P) circle (1pt);
\draw[fill=black] (Q) circle (1pt);
\draw[fill=black] (R) circle (1pt);
\draw[fill=black] (S) circle (1pt);
\draw (P) -- (Q) -- (R) -- (S) -- cycle;

\coordinate (x) at (3,1.6);
\coordinate (y) at (1.6,3);
\coordinate (z) at (4,4);
\draw[fill=black] (x) circle (1pt);
\draw[fill=black] (y) circle (1pt);
\draw[fill=black] (z) circle (1pt);

\draw (x) node[anchor=north]{$x$} -- (y) node[anchor=east]{$y$} -- (z) node[anchor=south]{$z$} -- cycle;
\draw (z) -- (Q) -- (x) -- (P) -- (y) -- (S) -- (z) -- (R) ;

\coordinate (star) at (3,3);
\draw (2.75,2.75) node{$\star$};
\draw[fill=black] (star) circle (1pt);
\draw (star) -- (x);
\draw (star) -- (y);
\draw (star) -- (z);

\draw (3.5,5.2) node{$A_1$};
\draw (5.2,3.5) node{$A_2$};
    \end{tikzpicture}
    
                        \qquad
                         
    \begin{tikzpicture}[scale=.5]
    \coordinate (P) at (0,0);
\coordinate (Q) at (3,0);
\coordinate (R) at (3,3);
\coordinate (S) at (0,3);

\draw[fill=black] (P) circle (1pt);
\draw[fill=black] (Q) circle (1pt);
\draw[fill=black] (R) circle (1pt);
\draw[fill=black] (S) circle (1pt);
\draw (P) -- (Q) -- (R) -- (S) -- cycle;

\coordinate (x) at (0,0.6);
\coordinate (y) at (0,0.8);
\coordinate (z) at (0,1);
\draw[fill=black] (x) circle (1pt);
\draw[fill=black] (y) circle (1pt);
\draw[fill=black] (z) circle (1pt);

\draw (x) -- (y) -- (z) -- cycle;
\draw (z) -- (Q) -- (x) -- (P) -- (y) -- (S) -- (z) -- (R) ;

\coordinate (star) at (6,0);
\draw[fill=black] (star) node[anchor=south]{$\rho(\star)$} circle (1pt);
\draw (star) -- (x);
\draw (star) -- (y);
\draw (star) -- (z);
    \end{tikzpicture}
    
                        \caption{A triangulation and a point on a path to the base locus.}
                        \label{fig:bubble}
                    \end{figure}
                
                Consider the path $\gamma(s)$ of drawings defined by
                    \begin{align*}
                    \PP&\longmapsto(0,0) \\
                    \QQ&\longmapsto(s,0) \\ %\ \clubsuit (1/m, 0)\  {\rm \ or\ some\ such, etc.\ } \\
                    \RR&\longmapsto(s,s) \\
                    \SS&\longmapsto(0,s) \\
                    x&\longmapsto(0,as^2) \\
                    y&\longmapsto(0,bs^2) \\
                    z&\longmapsto(0,cs^2) \\
                    \star&\longmapsto(1,0). \\
                    \end{align*}
                
                Figure \ref{fig:bubble} (right) shows the image of $\gamma(1/2)$, with $a,b,c$ small positive numbers.  One easily calculates that as $s\to 0$, $\area(\gamma(s))\to w$.  This establishes that $w\in V$ as claimed, and that $\gamma$ is a generating path for $w$.
                The point here is that the areas of the triangles around $\star$ shrink at the same rate (i.e., quadratically in $s$) as the total area.
                
                We now show $w\not\in \im\area\cup H$.
                Clearly $w\not\in H$.  Suppose for contradiction that $w=\area(\rho)$ for some $\rho\in X$.  Then since $w\not\in H$, $\rho$ is a drawing with nondegenerate boundary and by applying an element of $\Aff_2$ we may assume the boundary of $\rho$ is the unit square.  The coordinates of $w=\area(\rho)$ now impose many conditions on $\rho$.  For instance $\rho(x)$ must be on the $x$ axis, and $\rho(y)$ must be on the $y$ axis.  At least one of them must be the origin.  Also $\rho(z)$ is on both the $x$ axis and the $y$ axis, hence it is at the origin.  This makes (at least) one of the triangles around $\star$ degenerate.  But as $a,b,c$ are distinct, this is a contradiction.
                
                We therefore have a point $w\in V$ which is neither in $\im\area$ nor in $H$.
                
                Notice that all areas are shrinking
                to zero, but at different rates.  The slowest to shrink are ($A_1,A_2,A_8,A_9,A_{10}$), which shrink quadratically in $s$, as does the total area, so these five triangles are the ones which survive in the limit.  The area of each of these triangles is a nonzero finite fraction of the total area, in the limit.
                By contrast, the other areas shrink more rapidly, so that as $s$ tends to zero they each represent a vanishingly small fraction of the total area.  Lemma \ref{key.lemma} shows that if any individual triangle shrinks more slowly than the total, the limiting point of the area variety will be in $H$.
            \end{example}
            
            \begin{example}[Double bubble]\label{ex:doublebubble}
                The above phenomenon can happen at various speeds, producing nested bubbles.
                In this example, at time $s$ the boundary is a square of side length $s^2$.  The subscript $i$ is in $\{1,2,3\}$. Fix constants $a_i,b_i$.  The point $\star$ is in the center, the surrounding (orange) triangle has vertices $x_i$, and the next (blue) layer has vertices $y_i$.  (In reality the blue edges are straight.)  Watching a movie as $s\to 0$, the point $\star$ is drawn at $(1,1)$ throughout, while one observes the following sequence of events:  first the points $x_i$  cluster together on the $x$-axis; then the points $y_i$ cluster together on the $y$-axis and head toward the origin; then the whole square shrinks to the origin, not including the cluster of $x$'s; finally the cluster of $x$'s drifts toward the origin.  We encourage the interested reader to determine which points of $V$ have generating paths of this form.\\
                
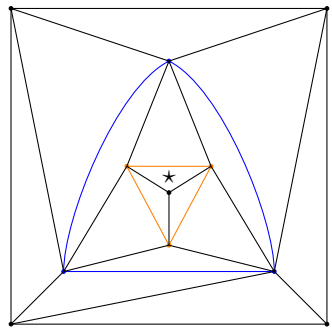
                \begin{equation*}
                    \vcenter{ 
    \begin{tikzpicture}[scale=.7]
    \coordinate (P) at (0,0);
\coordinate (Q) at (6,0);
\coordinate (R) at (6,6);
\coordinate (S) at (0,6);

\draw[fill=black] (P) circle (1pt);
\draw[fill=black] (Q) circle (1pt);
\draw[fill=black] (R) circle (1pt);
\draw[fill=black] (S) circle (1pt);
\draw (P) -- (Q) -- (R) -- (S) -- cycle;

\coordinate (x') at (5,1);
\coordinate (y') at (1,1);
\coordinate (z') at (3,5);
\coordinate (x) at (2.2,3);
\coordinate (y) at (3.8,3);
\coordinate (z) at (3,1.5);
\draw[orange, fill=orange] (x) circle (1pt);
\draw[orange, fill=orange] (y) circle (1pt);
\draw[orange, fill=orange] (z) circle (1pt);
\draw[blue, fill=blue] (x') circle (1pt);
\draw[blue, fill=blue] (y') circle (1pt);
\draw[blue, fill=blue] (z') circle (1pt);

\draw[blue] (y') .. controls +(0,1) and +(-1,-.5) .. (z') .. controls +(1,-.5) and +(0,1) .. (x') -- cycle;
\draw[orange] (x) -- (y) -- (z) -- cycle;
%\draw[blue] (y')-- (z')-- (x') -- cycle;

\coordinate (star) at (3,2.5);
\draw (3,2.8) node{$\star$};
\draw[fill=black] (star) circle (1pt);
\draw (star) -- (x);
\draw (star) -- (y);
\draw (star) -- (z);

\draw (P) -- (y') -- (S) -- (z') -- (R) -- (x') -- (Q);
\draw (P) -- (x');

\draw (y') -- (x) -- (z') -- (y) -- (x') -- (z) -- cycle;
    \end{tikzpicture}
    }
                    \hskip -2.5in
                    \begin{aligned}
                    \PP\QQ\RR\SS & \mapsto [0,s^2]\times[0,s^2]
                    \\
                    \textcolor{orange}{x_i} & \mapsto (s+a_is^4,0)
                    \\
                    \textcolor{blue}{y_i} & \mapsto (0,b_is^3)
                    \\
                    \star & \mapsto (1,1)
                    \end{aligned}
                \end{equation*}
                \begin{figure}[h]
                    \caption{Double bubble:  initially $\{\textcolor{orange}{x_i}\}$ forms a (the) bubble, but after it bursts a new bubble $\{\textcolor{blue}{y_i}\}$ forms.  (See Section \ref{subsec:bursting}.)}
                    \label{fig:double-bubble}
                \end{figure}
            \end{example}

        \subsection{Copointalities and bubbles}\label{subsec:purple}
                        
            The generating paths we consider have a limiting drawing $\rho$ at which all areas vanish.  What does this look like?  One way to make a drawing in which all the areas vanish is to draw all the vertices collinearly.  Indeed, if $\rho$ maps the vertices injectively to $\C^2$, then the vanishing of all areas implies that there is a single line containing the image of $\rho$ (see Lemma \ref{purple.lemma} below).  However, if $\rho$ has {\it copointalities}, i.e., distinct vertices of $T$ with the same image under $\rho$, then it is possible that the areas of the triangles all vanish even while there is no single line containing all the points $\rho(v)$.  For example, one can draw the triangulation described in Example \ref{ex:bubble} in such a way that all vertices except $\star$ are mapped onto a line, the vertices $x,y,z$ are copointal (on the line), and $\star$ is mapped somewhere off the line.
            
            To deal with this possibility, we introduce some language to discuss the combinatorics of the coincidences of points $\rho(v)$.  Recall that the domain of $\rho$ is $\Vxs(T)$.  Because we care about the triangles of an underlying triangulation, and coincidences involving vertices that are not connected by edges do not result in any triangles becoming degenerate, such coincidences do not concern us.  This motivates the following definitions.
        
            \begin{defn}[Cell-fiber, elastic complex]\label{def:cellfiber}
                Let $T$ be a triangulation, and let $\rho\in X$ be a drawing of $T$.  
                Let $\star\in\C^2$ be a point in the image of $\rho$. The subcomplex of $T$ spanned by $\rho^{-1}(\star)\subset\Vxs(T)$ is called the \emph{cell-fiber} of $\star$ and is denoted $\psi(\star)$; this is the union of all cells of $T$ (vertices, edges, or triangles) that are mapped entirely to $\star$ by a piecewise linear extension of $\rho$ to the 2-complex $T$.  (This is also a subset of the usual fiber over $\star$ of the same extension.)  
        
                We say that an edge $vv'$ of $T$ is \emph{elastic} in $\rho$ if $\rho(v)= \rho(v')$, or equivalently if the edge is contained in a cell-fiber.
        
                The \emph{elastic complex} of $\rho$ is the union of all its cell-fibers, minus any isolated vertices.  The elastic complex is denoted $\Psi(\rho)$.
            \end{defn}
        
            Note that it is possible for half or all of the boundary edges of $T$ to be elastic; if this happens then $\rho(T)$ is a degenerate parallelogram, i.e., the boundary of the square is mapped to points in the plane that are collinear.
            
            \begin{lem}[The elastic lemma]\label{purple.lemma}
                Fix $T$ and $\rho\in X$ 
                and suppose that $\area_i(\rho)=0$ for all $i$.
                Let $B$ be the closure (in $|T|$) of a component of $T\setminus \Psi(\rho)$.  
                Then there is a unique line $\ell$ in $\C^2$ such that $\rho(v)\in\ell$ for every vertex $v$ of $B$.
            \end{lem}
            
            \begin{proof}
                Let $\Delta$ be a triangle of $B$.  Note that at most one edge of $\Delta$ is elastic, since $B$ is not contained in any cell-fiber.  Hence, the three vertices of $\Delta$ map to either two or three distinct points.  Since $\area_{\Delta}(\rho)=0$, these points lie on a unique line $\ell$.  Now consider another triangle $\Delta'$ of $B$ that shares a non-elastic edge with $\Delta$.  For the same reason the images of the vertices of $\Delta'$ are also collinear, and as $\Delta'$ shares a non-elastic edge with $\Delta$ the third vertex of $\Delta'$ must also lie on $\ell$.  As all triangles of $B$ can be reached by crossing non-elastic edges, continuing this process eventually exhausts all triangles of $B$, which leads to the desired conclusion.
            \end{proof}
            
            So, if $\rho\in L$ then all triangles have area $0$, and there is a uniquely defined collection of lines $\ell_B$, one for each $B$ (the closure of a component of $T\setminus \Psi(\rho)$), such that $\rho$ maps all vertices lying in $B$ to the line $\ell_B$.  If $B$ and $B'$ intersect at the point $v$ then the lines $\ell_B$ and $\ell_{B'}$ intersect at the point $\rho(v)$.
            
            Let $C$ be a cycle in the 1-skeleton of $T$; this means that $C$ is a subgraph of the 1-skeleton that is homeomorphic to a circle.  We wish to refer to the ``inside'' of $C$.  Precisely, we introduce a subcomplex $\In(C)$ of $T$ as follows. 
            
            Let $\hat T$ be the cell complex $T$ together with a square 2-cell $Q$ 
            attached to the boundary of $T$, so that $|\hat T|$ is homeomorphic to a 2-sphere.  
            Then $\hat T \setminus C$ consists of two connected components, one of which does
            not intersect the interior of $Q$.  Denote the closure of this component by
            $\In(C)$.
            
            Fix a drawing $\rho$ in the base locus.  Let $C$ be a cycle in the 1-skeleton of $T$.  We adopt the following notation for the remainder of this section.
            Label the vertices of $C$ by $S_1, \dots, S_n$ in order and let $R_i$ be the vertex of $\In(C)$ which is adjacent to both $S_i$ and $S_{i+1}$ (indices mod $n$). Thus $S_iS_{i+1}R_i$ is a triangle of $\In(C)$. Note that the $R_i$ need not be distinct and that $R_i=S_j$ is possible.
            
            \begin{defn}[Bubble]
                We call $C$ a \emph{bubble} of $\rho$ if there is a point $\star\in \C^2$ such that
                \begin{itemize}
                    \item[(i)] for all $i$, $\rho(S_i) = \star$, 
                    \item[(ii)] for all $i$, $\rho(R_i)\neq\star$, and
                    \item[(iii)] all $R_i$ are in the same connected component of $T\setminus \psi(\star)$.  
                \end{itemize}
                If the cycle $C$ is a bubble of $\rho$ then we also refer to the subcomplex $\In(C)$ as a bubble of $\rho$.
            \end{defn}
            
            Notes:  The definition implies that the edges $R_iS_i$ (and therefore $R_iS_j$) are not elastic, that $R_i\ne S_j$, and that a bubble must contain at least one vertex (e.g., $R_1$) in its interior.  Also, (the inside of) a bubble is a contiguous set of triangles.  
        
            The triangles of the form $S_iS_{i+i}R_i$, with indices taken cyclically, play an important role.  We will refer to such triangles as $SSR$ triangles.
            
            The most important fact about bubbles is that they must be lurking whenever we have a point of $V$ that is not in the image of the area map. 
            
            \begin{thm}[Existence of Bubbles]\label{thm:frippfix}
                Let $T$ be a triangulation, let $V=V(T)$, and suppose $w\in V\setminus\im\area$. Then there exists a generating path $\gamma(s)$ for $w$ such that the drawing $\rho=\gamma(0)$ has a bubble.  
            \end{thm}
            
            \begin{proof} 
                Take any generating path $\gamma$ for $w$.  If $\rho=\gamma(0)$ has a bubble, we are done, so let us assume that $\rho$ has no bubble. 
                
                Suppose that the boundary $\PP\QQ\RR\SS$ is entirely elastic.
                Let $e$ be a non-elastic edge of $T$, which is guaranteed to exist by the definition of a generating path. 
                Among all cycles $B\subset\Psi$ with $e\subset\In(C)$ choose one, so that $\In(B)$ is minimal with respect to set inclusion.
                Such a $B$ is a bubble, contradicting our assumption that $\rho$ has no bubble.
                We conclude that the boundary cycle $\PP\QQ\RR\SS$ is not entirely elastic.
                
                Therefore, without loss of generality, we may assume $\rho(\PP)\ne \rho(\QQ)$, say $\rho(\PP)=(0,0)$ and $\rho(\QQ)=(1,0)$.  For each $s$ sufficiently close to $0$, there is the unique affine transformation $h(s)$ taking $\gamma(s)(\PP)$ to $(0,0)$, $\gamma(s)(\QQ)$ to $(1,0)$, and fixing $(0,1)$.  The transformation $h(s)$ converges to the identity as $s\to 0$. 
                By applying $h(s)$ to each $\gamma(s)$, we may assume that $\gamma(s)(\PP)=(0,0)$ and $\gamma(s)(\QQ)=(1,0)$ for all $s$.   Note that because $h(s)$ tends to the identity, we maintain the property that the drawings $\gamma(s)$ converge to the drawing $\rho$. 
                
                Since there is no bubble in $\rho$, it follows from the Elastic Lemma that  all vertices
                of $\rho$ are collinear.  Therefore by our choice of coordinates all vertices of $\rho$ lie 
                on the $x$ axis.  Since all of the $y$-coordinates of the vertices are approaching $0$, we may pick a vertex $v_0$ such that the order of vanishing of the $y$-coordinate of $v_0$ is minimal among all vertices. We then rescale in the $y$ direction so that $v_0$ approaches a point not on the $x$-axis and all vertices still have a limit.  To do this, for each $s\neq 0$, define a new drawing $\hat{\gamma}(s)$ by mapping vertex 
                $v$ of $T$ to the point $\left(x(\gamma(s)(v)), \frac{ y(\gamma(s)(v))}{y(\gamma(s)(v_0))} \right)\in\C^2$.
                We then see that by the choice of $v_0$, we still have  $\hat{\gamma}(s)$ approaching a limit, which we denote by $\hat{\rho}\in X$.  Since $\area(\hat{\gamma}(s))=\area(\gamma(s))$ for all $s\neq 0$, we see that $\hat{\gamma}$ is still a generating path for $w$.
                However, now the vertices of the limiting drawing $\hat{\rho}$ are not all collinear, since $\PP$ and $\QQ$ map to 
                $(0,0)$ and $(1,0)$, and the vertex $v_0$ maps to a point in $\C^2$ whose $y$-coordinate equals 1.  From the Elastic Lemma, it follows that $\hat{\rho}$ must contain a bubble, and $\hat{\gamma}$ is the desired generating path.
            \end{proof}
            
            We will need a slightly refined version of this theorem. To any generating path $\gamma$, we associate a nonnegative integer $r(\gamma)$, the order of vanishing of the total area $\sigma(\gamma(s))$ as $s\to 0$.  If $r=0$, then $\tot(\rho)\neq 0$, so $\rho$ is in the domain of $\area$, and so $w\in \im\area$. 
            
            We record the fact that in our proof of the Existence of Bubbles, the value of $r$ did not increase.
            
            \begin{thm}[Existence of Bubbles, Refined] 
                Let $T$ be a triangulation. Suppose $w\in V\setminus\im\area$,  and let $\gamma$ be a generating path for $w$.
                If $w\not\in\im\area$ then there exists a generating path $\tilde{\gamma}$ for $w$ such that $\tilde{\rho}=\tilde{\gamma}(0)$ has a bubble and $r(\tilde{\gamma})\leq r(\gamma)$.  
            \end{thm}
            \begin{proof}
                The proof of the Existence of Bubbles involves three modifications of the path $\gamma$.  The first is an affine transformation used to get $\rho(P)=(0,0)$ and $\rho(Q)=(1,0)$.  This is a fixed transformation, not depending on $s$, so the the value of $r$ is unchanged. The second is the affine transformation $h(s)$. Since $h(s)$ converges to the identity, it also does not change the value of $r$.  The third modification, the rescaling in the $y$-direction, decreases $r(\gamma)$ by the order of vanishing of $y(v_0)$.  The theorem follows. 
            \end{proof}
        
        \subsection{Bursting bubbles}\label{subsec:bursting}
        
            To deal with the presence of bubbles, we introduce the following procedure called bubble bursting.  The bubble exists because the vertices $S_i$ on $C$ are mapped close together, while some of the vertices in the interior of $C$ are mapped far away.  We burst the bubble by bringing these distant vertices close to the $S_i$.  In particular, we will map all the vertices in the interior of $C$ to the center of mass of the cluster of vertices $S_i$.  
            
            \begin{defn} [Bursting]
                Let $T$ be a triangulation. 
                Let $w\in V=V(T)$, and let $\gamma(s)$ be a generating path for $w$.  Suppose that $C$ is a bubble of $\rho=\gamma(0)$. We define a new path of drawings  $\bar{\gamma}$ which we refer to as {\em bursting the bubble} $C$. For any  $s$, define $\bar{\gamma}(s)(v)$ to agree with $\gamma(s)(v)$ on all vertices $v\notin\In(C)$ and also on all vertices of $C$, and to map all vertices  $v$ in the interior of $\In(C)$  to the center of mass of the $S_i$, that is we set  $\bar{\gamma}(s)(v) := \frac 1 n \sum_i \gamma(s)(S_i)$.
            \end{defn}
            
            Note that the new limiting drawing $\bar{\rho}=\bar{\gamma}(0)$  maps all vertices $v\in\In(C)$ to the same point and agrees with $\rho$ for all vertices $v\notin\In(C)$.  This has the effect of making $\In(C)$ entirely elastic;  indeed, we have $\Psi(\bar{\rho})=\Psi(\rho)\cup \In(C)$.
            
            Meanwhile, over in area space, bubble bursting has the effect of zeroing out all coordinates inside the bubble, as shown in the following theorem. 
            
            \begin{thm}[Bubble bursting] \label{thm:bubblebursting}
                Let $T$ be a triangulation, and let $V=V(T)$. 
                Let $w\in V\setminus H$.  Suppose that $\gamma$ is a generating path for $w$ that contains a bubble $C$,
                and let  $\bar{\gamma}$ be the result of bursting the bubble $C$. Then as $s\to 0$, $\area(\bar{\gamma}(s))$ converges to the point $\bar{w}$ whose coordinates outside $C$ agree with those of $w$ and whose coordinates inside $C$ are all zero.  In particular, this point $\bar{w}$ is in $V$.
            \end{thm}

            It will follow from the proof that indeed $\bar w\in Y(T)$, i.e., we have not set all coordinates to zero.  We also point out that $\bar{\gamma}$ may not be a generating path for $\bar w$, as the drawing $\bar{\gamma}(0)$ may be a constant map.  When we need a generating path, as we will in the proof of Theorem \ref{thm:characterization} in the next section, we will rescale $\bar{\gamma}$ as necessary.

            Before diving into the proof, we first prove a few lemmas. 
            
            \begin{lem}[Quad Lemma]
                Let $A,B,C,D$ be continuous paths $(-\epsilon, \epsilon) \to \C^2$ and let 
                $\sigma: (-\epsilon, \epsilon) \to \C$  be nonvanishing on a deleted neighborhood of $s=0$. Assume that
                \begin{itemize}
                    \item     $ A(0) \ne C(0)$;
                    \item Both $\frac{\area ACB}{\sigma}$ and $\frac{\area ACD}{\sigma}$ are bounded in a deleted neighborhood of $s=0$;
                \end{itemize}
                Then both $\frac{\area ABD}{\sigma}$ and $\frac{\area CBD}{\sigma}$ are also bounded in a deleted neighborhood of $s=0$.
            \end{lem}
            
            Under the hypotheses of the Quad Lemma we refer to $AC$ as the \emph{good diagonal} of the quadrilateral $ABCD$.
            
            \begin{proof}
                We may assume $A(0)=(0,0)$ and $C(0)=(1,0)$.  By applying an affine transformation $h(s)$ that converges to the identity we may further assume that $A(s)$ is constantly equal to $(0,0)$ and $C(s)=(1,0)$.  Note that the boundedness hypotheses and the desired boundedness conclusions are unchanged, because $h(s)$ approaches the identity.
                
                Now we may calculate.  We have $\frac{\area ACB}{\sigma}=\frac{y_B}{\sigma}$ and $\frac{\area ACD}{\sigma}=\frac{y_D}{\sigma}$ bounded by hypothesis.
                (Each variable is a function of $s$, though we have suppressed $s$ from the notation.)  Thus $\frac{\area(ABD)}{\sigma}=\frac{x_B y_D - x_D y_B}{\sigma}=x_B\frac{y_D}{\sigma} - x_D \frac{y_B}{\sigma}$ is bounded.
                It follows similarly that $\frac{\area(CBD)}{\sigma}$ is bounded, or alternatively one can observe that the areas of the four triangles formally sum to zero, so once three are bounded the fourth must be as well.
            \end{proof}
            
            \begin{lem}[Quad workout]
                Let $w\in V\setminus (\im\area\cup H)$ and let $\gamma$ be a generating path for $w$ with $\gamma(0)=\rho$.  Let $C$ be a bubble of $\rho$ with vertices $S_1,\dots,S_n$, and define $R_1,\dots,R_n$ as usual. Consider a triangle $\D$  of the form $S_iS_{i+1}R_j$, with $1\le i,j \le n$.
                Then along the path $\gamma$,
                \begin{equation}\label{badquad}
                    \frac{\area(\D)}{\sigma}\quad\mbox{ is bounded.}
                \end{equation}
            \end{lem}

            \begin{proof}
                Since $w\not\in H$, \eqref{badquad} holds for any triangle $\D$ of the triangulation by \ref{key.lemma}.  In particular \ref{badquad} holds already for $j=i$.  However when $j\ne i$, the triangle $\D$ may not be in the triangulation, and there is work to do.  For this purpose, we employ the quad lemma repeatedly, as follows.
                
                For simplicity we first suppose that there are no elastic edges of the form $S_i v$ where $v$ is in the interior of the bubble.
                
                Step 1.  
                    Let $R_1=x_0,x_1,\dots,x_p=R_2$ be a path in the link of $S_2$. See Figure \ref{fig:step12} (left).  
                    The first step of the proof is to show that \eqref{badquad} holds for all triangles of the form $S_1S_2x_k$ and $S_1 x_k x_{k+1}$.  
                    Apply the quad lemma successively to the quads $Q_k=S_2S_1x_kx_{k+1}$.  The good diagonal is $S_2x_k$; when $k=0$ this is because the triangles are in the triangulation, and for larger $k$ one triangle is in the triangulation and \eqref{badquad} holds for the other by induction.  Note that this is where we use the simplifying assumption that $S_2x_k$ is not elastic.  The quad lemma applied to each $Q_k$ ($0\le k<p$) now gives the desired result.  In particular we have shown that \eqref{badquad} holds for the triangle $S_1S_2R_2$.

                    The argument in the previous paragraph applies equally well anywhere along the cycle $C$. For example, letting $R_2=y_0,y_1,\dots,y_p=R_3$ be a path in the link of $S_3$ we have that \eqref{badquad} holds for the triangle $S_2y_ky_{k+1}$ for each $0\le k< p$.
                    
                    \begin{figure}
                        \centering
                        \begin{tikzpicture}[scale=.5]
                        \coordinate (S1) at (0,3);
\coordinate (S2) at (0,0);
\coordinate (S3) at (0,-3);
\coordinate (Sn) at (1,5);
\coordinate (R1) at (3,2);
\coordinate (R2) at (3,-2);
\coordinate (x1) at (4,1);
\coordinate (x2) at (4.2,0);
\coordinate (x3) at (4,-1);

\draw[fill=black] (S1) circle (1pt);
\draw[fill=black] (S2) circle (1pt);
\draw[fill=black] (S3) circle (1pt);
\draw[fill=black] (Sn) circle (1pt);
\draw[fill=black] (R1) circle (1pt);
\draw[fill=black] (R2) circle (1pt);
\draw[fill=black] (x1) circle (1pt);
\draw[fill=black] (x2) circle (1pt);
\draw[fill=black] (x3) circle (1pt);

\draw (S1) node[anchor=east]{$S_1$};
\draw (S2) node[anchor=east]{$S_2$};
\draw (S3) node[anchor=east]{$S_3$};
\draw (Sn) node[anchor=east]{$S_n$};
\draw (R1) node[anchor=south west]{$R_1=x_0$};
\draw (R2) node[anchor=north west]{$R_2=x_p$};
\draw (x1) node[anchor=west]{$x_1$};
\draw (x2) node[anchor=west]{$x_2$};

\draw (Sn) -- (S1) -- (S2) -- (S3);
\draw (S1) -- (R1) -- (S2) -- (R2) -- (S3);
\draw (R1)--(x1)--(x2)--(x3) --(R2);
\draw (S2) -- (x1);
\draw (S2) -- (x2);
\draw (S2) -- (x3);

\draw[red] (S2) -- (S1) -- (x1) -- (x2) -- cycle;
                        \end{tikzpicture}
                        \qquad\qquad
                        \begin{tikzpicture}[scale=.5]
                        \coordinate (S1) at (0,4);
\coordinate (S2) at (-1,0);
\coordinate (S3) at (-1,-3);
\coordinate (S4) at (0,-6);
\coordinate (R1) at (3,2);
\coordinate (R2) at (4,0);
\coordinate (R3) at (4,-5);
\coordinate (x1) at (4.4,-1);
\coordinate (x2) at (4.2,-2);
\coordinate (x3) at (4.4,-3);

\draw[fill=black] (S1) circle (1pt);
\draw[fill=black] (S2) circle (1pt);
\draw[fill=black] (S3) circle (1pt);
\draw[fill=black] (S4) circle (1pt);
%\draw[fill=black] (Sn) circle (1pt);
\draw[fill=black] (R1) circle (1pt);
\draw[fill=black] (R2) circle (1pt);
\draw[fill=black] (R3) circle (1pt);
\draw[fill=black] (x1) circle (1pt);
\draw[fill=black] (x2) circle (1pt);
\draw[fill=black] (x3) circle (1pt);

\draw (S1) node[anchor=east]{$S_1$};
\draw (S2) node[anchor=east]{$S_2$};
\draw (S3) node[anchor=east]{$S_3$};
\draw (S4) node[anchor=east]{$S_4$};
\draw (R1) node[anchor=south west]{$R_1$};
\draw (R2) node[anchor=south west]{$R_2$};
\draw (R3) node[anchor=north west]{$R_3$};
\draw (x1) node[anchor=west]{$y_1$};
\draw (x2) node[anchor=west]{$y_2$};
\draw (x3) node[anchor=west]{$y_3$};

\draw (S1) -- (S2) -- (S3) -- (S4);
\draw (S1) -- (R1) -- (S2) -- (R2) -- (S3) -- (R3) -- (S4);
\draw (R2)--(x1)--(x2)--(x3) --(R3);
\draw (S3) -- (x1);
\draw (S3) -- (x2);
\draw (S3) -- (x3);

\draw[red] (S2) -- (S1) -- (x1) -- (x2) -- cycle;
\draw[red,dotted] (S2)--(x1);
                        \end{tikzpicture}
                        \caption{Step 1, Step 2}
                        \label{fig:step12}
                    \end{figure}
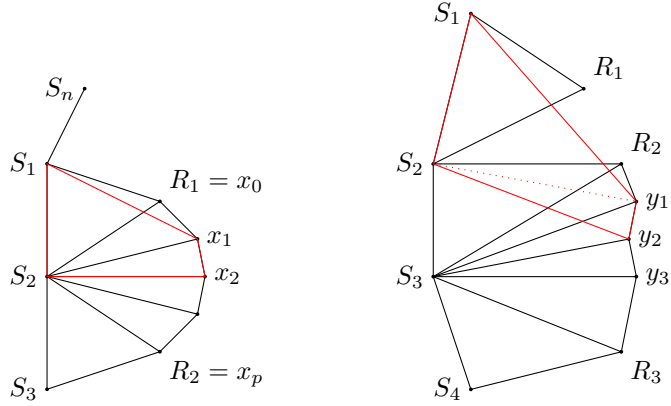
                
                Step 2.  
                    We now consider the quad (and reuse notation) $Q_k=S_2 S_1 y_k y_{k+1}$.  Starting with $k=0$ and proceeding step by step, we see inductively (and by referring to Step 1) that $S_2 y_k$ is a good diagonal of $Q_k$. (The quad $Q_1$ and its good diagonal are shown in red in Figure \ref{fig:step12}.)  To get started, for instance, the diagonal $S_2R_2=S_2y_0$ of $Q_0$ is good because Step 1 established that \eqref{badquad} holds for both $S_2y_0y_1$ and $S_2S_1y_0$. Applying the Quad Lemma to $Q_0$ shows that \eqref{badquad} holds for triangles $\Delta=S_1S_2y_1$ as well as $\Delta=S_1y_0y_1$.  The first of these, combined with our work in Step 1 showing that \eqref{badquad} holds for $S_2y_1y_2$, implies that $S_2 y_1$ is a good diagonal of $Q_1$.  Apply the Quad Lemma and repeating for each $k$  shows that \eqref{badquad} holds for each $\Delta=S_1S_2y_k$ as well as $\Delta=S_1y_ky_{k+1}$ (for $0\le k < p$).  In particular \eqref{badquad} holds for $\D=S_1S_2R_3$, and as in Step 1 this also applies generally to $\D=S_iS_{i+1}R_{i+2}$.
                    
                Step 3.  
                    Continuing in this fashion establishes the desired conclusion, namely that \eqref{badquad} holds for all $\D=S_iS_{i+1}R_j$.  As in Step 2, verifying that the hypotheses of the Quad Lemma hold in each step requires the result of all previous steps.  This completes the argument in the simplified case that none of the edges $S_i v$ with $v$ in the interior of the bubble are elastic.
                    
                    Suppose now we are in the general case.  We will create a setup in which the previous arguments can be carried out.  We will have $V$'s in place of $S$'s and $U$'s in place of $R$'s.
                    
                Step 0.  
                    Let $V_1,\dots,V_m$ be the closed walk that begins $S_1,S_2,\dots$ (counterclockwise around the outside of the bubble) and traverses the boundary of the component of the complement of the elastic complex whose closure forms the bubble.  That is, $V_1=S_1$ and $V_2=S_2$ and each edge $V_iV_{i+1}$ is elastic and, crucially, pushing the entire path $V_1,\dots,V_m$ slightly into the bubble (to the left relative to the direction of the path) gives a path disjoint from the elastic complex.  
                    
                    Note that the awkwardness of this definition arises because in general the $V_i$ need not be distinct, and the path of $V$'s may traverse some edges once in each direction.  If we happen to be in the simplified case we first considered, we will have $m=n$ and $V_i=S_i$ for each $i$.  
                    
                    This is all illustrated in Figure \ref{fig:step0}.  The $S$'s and $R$'s are as before, but the elastic complex (drawn more darkly than the rest) includes an edge from $S_3$ into the interior of the bubble as well as a triangle that is entirely contained in the interior of the bubble.  The first several $V$'s are labeled. The red path is the aforementioned embedded topological circle that is disjoint from the elastic complex. 

                    Next let $U_i$ be the vertex connected to both $V_i$ and $V_{i+1}$ on the left in the direction of travel. This is analogous to $R_i$, and is equal to $R_i$ in the simplified case we first considered.  The first several $U$'s are also labeled in Figure \ref{fig:step0}.
                    
                    \begin{figure}
                        \centering
                        \begin{tikzpicture}[scale=.5]
                        \coordinate (S1) at (0,4);
\coordinate (S2) at (-1,0);
\coordinate (S3) at (-1,-3);
\coordinate (S4) at (0,-6);
\coordinate (R1) at (3,2);
\coordinate (R2) at (4,0);
\coordinate (R3) at (4,-7);
\coordinate (V4) at (8,-3);
\coordinate (V5) at (10,0);
\coordinate (V6) at (13,-4);

\coordinate (U3) at (5,-2);
\coordinate (U4) at (8,-1);
\coordinate (U5) at (13,-1);
\coordinate (U6) at (10,-5);
\coordinate (U7) at (5,-5);

\draw[fill=black] (S1) circle (1pt);
\draw[fill=black] (S2) circle (1pt);
\draw[fill=black] (S3) circle (1pt);
\draw[fill=black] (S4) circle (1pt);
\draw[fill=black] (R1) circle (1pt);
\draw[fill=black] (R2) circle (1pt);
\draw[fill=black] (R3) circle (1pt);
\draw[fill=black] (V4) circle (1pt);
\draw[fill=black] (V5) circle (1pt);
\draw[fill=black] (V6) circle (1pt);
\draw[fill=black] (U3) circle (1pt);
\draw[fill=black] (U4) circle (1pt);
\draw[fill=black] (U5) circle (1pt);
\draw[fill=black] (U6) circle (1pt);
\draw[fill=black] (U7) circle (1pt);

\draw (S1) node[anchor=east]{$S_1=V_1$};
\draw (S2) node[anchor=east]{$S_2=V_2$};
\draw (S3) node[anchor=east]{$S_3=V_3=V_8$};
\draw (S4) node[anchor=east]{$S_4=V_9$};
\draw (R1) node[anchor=south west]{$R_1=U_1$};
\draw (R2) node[anchor=south west]{$R_2=U_2$};
\draw (R3) node[anchor=north west]{$R_3=U_8$};
\draw ($(V4)+(3,-5)$) node[anchor=north] {$V_4=V_7$};
    \draw [->] ($(V4)+(3,-5)$) .. controls ($(V4)+(0,-3)$) .. ($(V4)+(0,-1)$);

\draw (V5) node[anchor=south west]{$V_5$};
\draw (V6) node[anchor=south west]{$V_6$};
\draw (U3) node[anchor=south]{$U_3$};
\draw (U4) node[anchor=east]{$U_4$};
\draw (U5) node[anchor=west]{$U_5$};
\draw (U6) node[anchor=north]{$U_6$};
\draw (U7) node[anchor=north]{$U_7$};

\draw[very thick] (S1) -- (S2) -- (S3) -- (S4);
\draw (S1) -- (R1) -- (S2) -- (R2) -- (S3) -- (R3) -- (S4);
\draw[very thick] (S3)--(V4);
\draw[fill=gray] (V4)--(V5)--(V6) --(V4);

\draw (S3)--(U3)--(V4)--(U4)--(V5)--(U5)--(V6)--(U6)--(V4)--(U7)--(S3);
%  $(Q)+(0.5, 0)$

\draw[red, rounded corners] 
    ($(S1) +(.2,0)$) --  
    ($(S2) +(.2,0)$) -- 
    ($(S3) +(.2,.2)$) -- 
    ($(V4) +(-.1,.2)$) --
    ($(V5) +(0,.4)$) --
    ($(V6) +(.5,-.3)$) --
    ($(V4) +(0,-.2)$) -- 
    ($(S3) +(.3,-.2)$) -- 
    ($(S4) +(.2,0)$);
                        \end{tikzpicture}
                        \caption{Step 0}
                        \label{fig:step0}
                    \end{figure}
                                        
                    The previous arguments now translate directly to this situation, with $S$'s replaced by $V$'s and $R$'s replaced by $U$'s.  Recall that $V_1=S_1$ and $V_2=S_2$.  We first fill in a path $x_0,\dots,x_p$ from $U_1$ to $U_2$ in the link of $V_2$, which is possible since all $U$'s are in the same connected component of the complement of the elastic complex. Next we show that \eqref{badquad} holds for $V_1V_2U_2$ by using the quad lemma a bunch of times to show that \eqref{badquad} holds for each $V_1V_2x_k$ and $V_1x_kx_{k+1}$.  The quad lemma applies here because none of the edges $V_2x_k$ is elastic by construction.  By continuing to mimic the previous argument we obtain the corresponding conclusion that \eqref{badquad} holds for each $V_1V_2U_j$, hence for each $V_iV_{i+1}U_j$.
                
                The proof is completed by noting that the walk of $V$'s includes every edge $S_iS_{i+1}$, hence every $R_j$ is found among the $U_j$, and so the collection of triangles now known to satisfy \eqref{badquad} includes all those asserted to do so in the conclusion of the lemma.
            \end{proof}
            
            With the Quad Workout established, we now proceed to the proof of the Bubble Bursting Theorem.  
            
            \begin{proof}[Proof of bubble bursting theorem]
                First, normalize the coordinates of $w$ so that $\sum w_i=1$, which is possible since $w\notin H$.  Thus for all $i$, $\frac{\area_i(\gamma(s))}{\tot(\gamma(s))} \to w_i $.
                Now, if we replace $\gamma(s)$ with $\bar{\gamma}(s)$,  the values of $\tot$ have not changed, nor have the values $\area_i$ for triangles outside $C$.  Hence for these $i$, we have $\frac{\area_i(\bar{\gamma}(s))}{\tot(\gamma(s))} \to w_i = \bar{w}_i $.
                
                We now claim that for $i$ indexing triangles inside $C$,  $\frac{\area_i(\bar{\gamma}(s))}{\tot(\gamma(s))} \to 0$. 
                
                In the case that $i$ indexes a triangle with zero or one vertex on $C$, this triangle in $\bar{\gamma}(s)$ has at least two coincident vertices for all $s$, and hence  $\area_i(\bar{\gamma})(s)$ vanishes for all $s$, and our claim is true. 
                As noted previously, there are no triangles in the triangulation with three vertices on $C$.  This leaves only the triangles with two vertices on $C$. Since $C$ is a bubble, these triangles must be $SSR$ triangles.  Hence we wish to show that $\frac{\area(S_iS_{i+i}R')}{\tot}\to 0$, where $R'$ denotes the center of mass of all the $S$'s (for any value of the parameter $s$, which we now suppress from the notation). 
                
                Pick any of the $R$'s, say $R_1$, and apply the Corollary \ref{cor:nugget} of Nugget 1, found in the Appendix, to the wheel consisting of the $n$-gon $S_1, \dots, S_n$ with hub $R_1$.   This corollary tells us that there is an $i$ such that the triangle $U=S_iS_{i+i}R_1$ satisfies $\frac{\area(S_iS_jS_k)}{\area(U)}\to 0$ along $\gamma$ for all $i,j,k$. 
                  
                Note that the triangle $U$ may not be a triangle of $T$. However the Quad Workout implies that  $\frac{\area(U)}{\tot}$ is bounded, hence the product $\frac{\area(S_iS_jS_k)}{\tot}\to 0$ for all $i,j,k$.  We now average over $k\in \{1, \dots n\}$,  noting that the area of any $S_iS_jR'$ is the average of the areas of the $S_iS_jS_k$. We thus obtain $\frac{\area(S_iS_jR')}{\tot}\to 0$. Taking $j=i+1$ now yields our desired conclusion.
                
                We have now shown that for all $i$, $\frac{\area_i(\bar{\gamma}(s))}{\tot(\gamma(s))} \to \bar{w}_i$. Since $\sum \frac{\area_i(\bar{\gamma}(s))}{\tot(\gamma(s))} $ is identically $1$ for all  $s\ne 0$, we see that $\sum \bar{w}_i=1$.  Thus not all  $\bar{w}_i$ vanish.  Hence $\bar{w}$ defines a point in the projective space $Y$, and $\area(\bar{\gamma}(s))\to \bar{w}$.   This completes the proof.  
            \end{proof}
        
            \begin{cor}[Bubble Corollary] \label{tight.cycle}
                Let $T$ be a triangulation, and let $V=V(T)$. 
                Let $w\in V\setminus H$  and let $\gamma$ be a generating sequence for $w$.
                Suppose that $C$ is a bubble of $\rho=\gamma(0)$.
                Then $$\sum_{i\in \In(C)} w_i = 0.$$
            \end{cor}
            
            \begin{proof}
                Noting that the area inside $C$ is the same along $\gamma$ and $\bar{\gamma}$, we see that
                $$\sum_{i\in \In(C)} \frac{\area_i(\gamma(s))}{\tot(\gamma(s))} =\sum_{i\in \In(C)} \frac{\area_i(\bar{\gamma}(s))}{\tot(\gamma(s))}.$$ 
                As we have seen, as $s\to 0$, the terms on the left approach $w_i$, while each term on the right approaches $0$.
            \end{proof}
            
            One should revisit Example \ref{ex:bubble} in light of what we have now proved.
            There we had a point $w\not\in H$, and the elastic complex $\Psi(\rho)$ consisted of all of $T$ minus a neighborhood of $p_4$.  The only bubble is (the boundary of) this neighborhood, which consists of triangles 8, 9, and 10,
            and sure enough, in $w$ the coordinates $A_8, A_9$, and $A_{10}$ sum to zero.
            
            \begin{cor}
                Let $T$ be a triangulation, and let $w\in V(T)$. Let $\gamma$ be a generating path for $w$ and let $\rho=\gamma(0)$.
                If $\bdry T$ is a bubble in $\rho$ then $w\in H$.
                In particular if $\Psi(\rho)=\bdry T$ then $w\in H$.
            \end{cor}
            
            \begin{proof}
                If $\bdry T$ is a bubble then $\In(\bdry T)=T$. It then follows from the previous corollary that $w\in H$.
            \end{proof}

    \section{Points of the area variety}
    
    Repeatedly bursting bubbles now allows us to achieve our goal of analyzing points of $V\setminus H$ that are not in the image of $\area$.  Indeed, consider any $w\in V\setminus H$ that is not in $\im\area$. Take a generating path with a bubble, which exists by Theorem \ref{thm:frippfix}.  Apply bubble bursting to all the bubbles.  It is easy to see that the order in which we burst the bubbles does not matter, even if there are nested bubbles; the resulting $\bar{\gamma}(s)$ is the same. The Bubble Bursting theorem tells us that the resulting $\bar{w}$ is still in $V$, and the Bubble Corollary implies that $\bar{w}$ is still not in $H$.  Note that the elastic complex tends to get bigger as all the bubbles fill in.  If the elastic complex is not all of $T$, then because there are no longer any bubbles, the resulting $\bar{w}$ is in the image of $\area$ by the Existence of Bubbles. 
    
    However, it is possible that after bursting all the bubbles, the elastic complex is now the entire $T$.  In this case, we can scale up, \`a la Lemma \ref{lem:existgen}, to obtain a new generating path for $\bar{w}$.  In doing so we decrease the order of vanishing $r(\gamma)$ of the area $\sigma$ of the square.
    Thus this process ends after finitely many iterations, and we eventually end up with a point in $\im\area$.  
    
    \begin{thm}[Points of $V$, see Main Theorem \ref{mainthm:points}]\label{thm:characterization}
        Let $T$ be a triangulation.   Let $w\in V=V(T)$.  Then either:
        \begin{enumerate}
            \item $w\in\im\area$; or
            \item $w\in H$; or
            \item the process of repeated bubble bursting and rescaling terminates after a finite number of steps in a point $\bar{w}$ in $\im\area$.
        \end{enumerate}
        In particular if the first two possibilities do not hold then there exists a subset of the coordinates of $w$ that sums to $0$, and such that the point $\bar{w}$ obtained by replacing all these coordinates with zeroes is in $\im\area$.
    \end{thm}
    
    \begin{proof}
      Start with $w\in V$, which we may assume is not in $\im\area$ and not in $H$. Take a generating path $\gamma$ for $w$ and apply Existence of Bubbles, Refined.
      Then burst all the resulting bubbles. By the Bubble Corollary, the resulting point is still not in $H$.  If this point is in the image, then we may stop; we are in case (3).  Otherwise, since the drawing that results from bubble bursting has no bubble, it follows from the Existence of Bubbles that the entire $1$-skeleton must be elastic.  In this case, after translating to the origin, every coordinate of every vertex vanishes to order at least $1$.  Thus we may divide by a positive power of $s$ and get a new generating path with order of vanishing strictly smaller than $r(\gamma)$.  Again apply the Existence of Bubbles, Refined, and continue bursting. Since $r$ does not increase with the Existence of Bubbles, and $r$ strictly decreases whenever we hit an all-elastic situation, this process must eventually terminate. The result is the desired $\bar{w}\in\im\area$.  The final assertion follows from the Bubble Corollary.
    \end{proof}
    
    \begin{cor}
        Suppose $w=[\cdots:A_i:\cdots]$ where the $A_i$ are nonnegative real numbers.
        (At least one of the $A_i$ must be positive of course.)  If $w\in V$ then $w\in\im\area$.
    \end{cor}
    
    \begin{proof}
        Note $w\not\in H$.  Suppose $w\not\in\im\area$.  The theorem asserts that repeated bubble bursting results in a point of the image.  But repeated bubble bursting can only result in $w$, so $w$ is in the image, a contradiction.
    \end{proof}

    See Example \ref{ex:complexequi} below for a caution about applying this corollary to dissections.

\part{Existence of the encyclopedia}
    
    Armed with Theorem \ref{thm:characterization}, we are now able to deduce various properties of the polynomial $p$.  Ultimately we prove a finiteness result for area polynomials (Theorem \ref{thm:encyc}) which we stated in the introduction as Main Theorem \ref{mainthm:encyc}.  This relies on a more elementary combinatorial finiteness result that we prove next.

    Here are some definitions that we will need.
    \begin{defn}[Subdivision]\label{def:subdivision}
        Let $T$ be a triangulation.  A \emph{subdivision} is a set of three vertices $v,w,x$ such that each pair is joined by an edge, but there is no triangle of $T$ with vertices $v,w,x$.
    \end{defn}
        
    For example the vertices $x,y,z$ in Figure \ref{fig:bubble} form a subdivision in the triangulation shown there.  We will be primarily interested in triangulations without subdivisions, as it is easy to describe the relationship between the area polynomials $p_T$ and $p_{T'}$ if $T'$ is obtained from $T$ by subdividing.
    
    \begin{defn}[Edge contraction]\label{def:contraction}
        Let $T$ be a triangulation and $e=xy$ an interior edge of $T$.  Then $e$ is on the boundary of two triangles $\Delta,\Delta'$ of $T$; denote the third vertices of these triangles by $z$ and $z'$.  Consider the complex $T|_e$ obtained from $T$ by deleting the triangles $\Delta$ and $\Delta'$ and identifying edges $zx$ and $zy$ and identifying edges $z'x$ and $z'y$ (and also identifying vertices $x$ and $y$).  The complex $T|_e$ may or may not be simplicial; if it is, then $e$ can be \emph{contracted} and $T|_e$ is the result of \emph{contracting} the edge $e$.  If $T|_e$ is not a simplicial complex, then $e$ cannot be contracted.
    \end{defn}
    
    It is an exercise to see that the interior edge $e=xy$ can be contracted if and only if $\{x,y\}$ is not part of a subdivision, and in this case the resulting $T|_e$ is again a triangulation (Definition \ref{def:honest}).
    
    \section{Combinatorial simplification}
    
        It is easy to see that there are finitely many combinatorial types of triangulation with a fixed number $k$ of interior vertices.  For $k\le 4$ we have classified such triangulations up to isomorphism by hand and verified our calculation using the computer program {\tt plantri} \cite{plantri}.  After filtering out triangulations with subdivisions, the results are shown in Table \ref{table:T.up.to.4.int} and summarized here:
            \begin{itemize}
                \item When $k=0$ there is one triangulation, $T_0$.
                \item When $k=1$ there is one triangulation, $T_1$.
                \item When $k=2$ there is one triangulation, $T_2$.
                \item When $k=3$ there are two triangulations, $T_3$ and the exploded diagonal $T_{2,1}$ (see \cite{triangles1} or Table \ref{table:T.up.to.4.int}).
                \item When $k=4$ there are exactly six triangulations:  $T_4$, $T_{3,1}$, and the four others shown and named in Table \ref{table:T.up.to.4.int}.
            \end{itemize}
        
            \begin{table}
                \begin{tabular}{c c c c c c}
                    \toprule
                    \addlinespace[5pt]
                     
    \begin{tikzpicture}[scale=.2]
    \coordinate (P) at (0,0);
\coordinate (Q) at (6,0);
\coordinate (R) at (6,6);
\coordinate (S) at (0,6);

\draw[fill=black] (P) circle (1pt);
\draw[fill=black] (Q) circle (1pt);
\draw[fill=black] (R) circle (1pt);
\draw[fill=black] (S) circle (1pt);
\draw (P) -- (Q) -- (R) -- (S) -- cycle;

\draw (P) -- (R);

\ifcorners
    \draw (P) node[anchor=north east]{$\PP$};
    \draw (Q) node[anchor=north west]{$\QQ$};
    \draw (R) node[anchor=south west]{$\RR$};
    \draw (S) node[anchor=south east]{$\SS$};
\fi

\ifareas
    \draw (4,2) node{$A$};
    \draw (1.8,4) node{$B$};
\fi

\ifnewareas
    \draw (4,2) node{\scriptsize $A$};
    \draw (1.8,4) node{\scriptsize $B$};
\fi

\ifcolors
    \draw[line width=0, fill=\poscolor] (P) -- (Q) -- (R) -- cycle;
    \draw[line width=0, fill=\negcolor] (R) -- (S) -- (P) -- cycle;
\fi
    \end{tikzpicture}
     &  &  &  &  & \\
                    $T_0$ & & & & & \\
                    \addlinespace[5pt]
                     
    \begin{tikzpicture}[scale=.2]
    \coordinate (P) at (0,0);
\coordinate (Q) at (6,0);
\coordinate (R) at (6,6);
\coordinate (S) at (0,6);

\draw[fill=black] (P) circle (1pt);
\draw[fill=black] (Q) circle (1pt);
\draw[fill=black] (R) circle (1pt);
\draw[fill=black] (S) circle (1pt);
\draw (P) -- (Q) -- (R) -- (S) -- cycle;

\coordinate (x) at (3,3);
\draw[fill=black] (x) circle (1pt);

\draw (P) -- (x);
\draw (Q) -- (x);
\draw (R) -- (x);
\draw (S) -- (x);

\ifcorners
    \draw (P) node[anchor=north east]{$\PP$};
    \draw (Q) node[anchor=north west]{$\QQ$};
    \draw (R) node[anchor=south west]{$\RR$};
    \draw (S) node[anchor=south east]{$\SS$};
\fi

\ifcolors
    \draw[line width=0, fill=\poscolor] (P) -- (Q) -- (x) -- cycle;
    \draw[line width=0, fill=\poscolor] (R) -- (S) -- (x) -- cycle;
    \draw[line width=0, fill=\negcolor] (Q) -- (R) -- (x) -- cycle;
    \draw[line width=0, fill=\negcolor] (S) -- (P) -- (x) -- cycle;
\fi

\ifareas
    \draw (3,1.25) node{$A$};
    \draw (4.75,3) node{$B$};
    \draw (3,4.75) node{$C$};
    \draw (1.25,3) node{$D$};
\fi

\ifnewareas
    \draw (3,1.25) node{\scriptsize $D$};
    \draw (4.75,3) node{\scriptsize $B$};
    \draw (3,4.75) node{\scriptsize $C$};
    \draw (1.25,3) node{\scriptsize $A$};
\fi
    \end{tikzpicture}
     &  &  &  &  & \\
                    $T_1$ & & & & & \\
                    \addlinespace[5pt]
                     
    \begin{tikzpicture}[scale=.2]
    \coordinate (P) at (0,0);
\coordinate (Q) at (6,0);
\coordinate (R) at (6,6);
\coordinate (S) at (0,6);

\draw[fill=black] (P) circle (1pt);
\draw[fill=black] (Q) circle (1pt);
\draw[fill=black] (R) circle (1pt);
\draw[fill=black] (S) circle (1pt);
\draw (P) -- (Q) -- (R) -- (S) -- cycle;

\coordinate (u) at (2,2.5);
\coordinate (v) at (4.5,3.75);
\draw[fill=black] (u) circle (1pt);
\draw[fill=black] (v) circle (1pt);

\draw (P) -- (u) -- (v) -- (R);
\draw (Q) -- (u) -- (S);
\draw (S) -- (v) -- (Q);

\ifcorners
    \draw (P) node[anchor=north east]{$\PP$};
    \draw (Q) node[anchor=north west]{$\QQ$};
    \draw (R) node[anchor=south west]{$\RR$};
    \draw (S) node[anchor=south east]{$\SS$};
\fi

\ifareas
    \draw (2.25,1) node{$A$};
    \draw (4,2.25) node{$B$};
    \draw (5.5,2.75) node{$C$};
    \draw (.75,2.75) node{$D$};
    \draw (2.25,4) node{$E$};
    \draw (4,5) node{$F$};
\fi

\ifinterior
    \draw (u) node[anchor=north]{$u$};
    \draw (v) node[anchor=west]{$v$};
\fi

\ifcolors
    \draw[line width=0, fill=\poscolor] (P) -- (Q) -- (u) -- cycle;
    \draw[line width=0, fill=\poscolor] (v) -- (Q) -- (R) -- cycle;
    \draw[line width=0, fill=\poscolor] (u) -- (v) -- (S) -- cycle;
    \draw[line width=0, fill=\negcolor] (u) -- (Q) -- (v) -- cycle;
    \draw[line width=0, fill=\negcolor] (v) -- (R) -- (S) -- cycle;
    \draw[line width=0, fill=\negcolor] (u) -- (S) -- (P) -- cycle;
\fi
    \end{tikzpicture}
     &  &  &  &  & \\
                    $T_2$ & & & & & \\
                    \addlinespace[5pt]
                     
    \begin{tikzpicture}[scale=.2]
    \coordinate (P) at (0,0);
\coordinate (Q) at (6,0);
\coordinate (R) at (6,6);
\coordinate (S) at (0,6);

\draw[fill=black] (P) circle (1pt);
\draw[fill=black] (Q) circle (1pt);
\draw[fill=black] (R) circle (1pt);
\draw[fill=black] (S) circle (1pt);
\draw (P) -- (Q) -- (R) -- (S) -- cycle;

\coordinate (w) at (1.5,1.5);
\coordinate (x) at (3,3);
\coordinate (y) at (4.5,4.5);
\draw[fill=black] (w) circle (1pt);
\draw[fill=black] (x) circle (1pt);
\draw[fill=black] (y) circle (1pt);

\draw (P) -- (w) -- (x) -- (y) -- (R);
\draw (Q) -- (w) ;
\draw (Q) -- (x) ;
\draw (Q) -- (y) ;
\draw (S) -- (w) ;
\draw (S) -- (x) ;
\draw (S) -- (y) ;

\ifcorners
    \draw (P) node[anchor=north east]{$\PP$};
    \draw (Q) node[anchor=north west]{$\QQ$};
    \draw (R) node[anchor=south west]{$\RR$};
    \draw (S) node[anchor=south east]{$\SS$};
\fi

\ifcolors
    \draw[line width=0, fill=\poscolor] (P) -- (Q) -- (w) -- cycle;
    \draw[line width=0, fill=\poscolor] (x) -- (Q) -- (y) -- cycle;
    \draw[line width=0, fill=\poscolor] (R) -- (S) -- (y) -- cycle;
    \draw[line width=0, fill=\poscolor] (x) -- (S) -- (w) -- cycle;
    \draw[line width=0, fill=\negcolor] (w) -- (Q) -- (x) -- cycle;
    \draw[line width=0, fill=\negcolor] (y) -- (Q) -- (R) -- cycle;
    \draw[line width=0, fill=\negcolor] (y) -- (S) -- (x) -- cycle;
    \draw[line width=0, fill=\negcolor] (w) -- (S) -- (P) -- cycle;
\fi
    \end{tikzpicture}
     &  
    \begin{tikzpicture}[scale=.2]
    \coordinate (P) at (0,0);
\coordinate (Q) at (6,0);
\coordinate (R) at (6,6);
\coordinate (S) at (0,6);

\draw[fill=black] (P) circle (1pt);
\draw[fill=black] (Q) circle (1pt);
\draw[fill=black] (R) circle (1pt);
\draw[fill=black] (S) circle (1pt);
\draw (P) -- (Q) -- (R) -- (S) -- cycle;

\coordinate (w) at (3,1.6);
\coordinate (x) at (1.6,3);
\coordinate (y) at (4,4);
\draw[fill=black] (w) circle (1pt);
\draw[fill=black] (x) circle (1pt);
\draw[fill=black] (y) circle (1pt);

\draw (w) -- (y) -- (x) -- cycle;
\draw (Q) -- (w) -- (P);
\draw (Q) -- (y);
\draw (P) -- (x) -- (S);
\draw (S) -- (y) -- (R) ;

\ifcorners
    \draw (P) node[anchor=north east]{$\PP$};
    \draw (Q) node[anchor=north west]{$\QQ$};
    \draw (R) node[anchor=south west]{$\RR$};
    \draw (S) node[anchor=south east]{$\SS$};
\fi

\ifcolors
    \draw[line width=0, fill=\poscolor] (P) -- (Q) -- (w) -- cycle;
    \draw[line width=0, fill=\poscolor] (w) -- (x) -- (y) -- cycle;
    \draw[line width=0, fill=\poscolor] (x) -- (S) -- (P) -- cycle;
    \draw[line width=0, fill=\negcolor] (P) -- (w) -- (x) -- cycle;
    \draw[line width=0, fill=\negcolor] (w) -- (Q) -- (y) -- cycle;
    \draw[line width=0, fill=\negcolor] (y) -- (Q) -- (R) -- cycle;
    \draw[line width=0, fill=\negcolor] (R) -- (S) -- (y) -- cycle;
    \draw[line width=0, fill=\negcolor] (x) -- (S) -- (y) -- cycle;
\fi
    \end{tikzpicture}
     &  &  & \\
                    $T_3$ & $T_{2,1}$ & & & & \\
                    %\midrule
                    \addlinespace[5pt]
                     
    \begin{tikzpicture}[scale=.2]
    \coordinate (P) at (0,0);
\coordinate (Q) at (6,0);
\coordinate (R) at (6,6);
\coordinate (S) at (0,6);

\draw[fill=black] (P) circle (1pt);
\draw[fill=black] (Q) circle (1pt);
\draw[fill=black] (R) circle (1pt);
\draw[fill=black] (S) circle (1pt);
\draw (P) -- (Q) -- (R) -- (S) -- cycle;

\coordinate (w) at (1.2,1.2);
\coordinate (x) at (2.4,2.4);
\coordinate (y) at (3.6,3.6);
\coordinate (z) at (4.8,4.8);
\draw[fill=black] (w) circle (1pt);
\draw[fill=black] (x) circle (1pt);
\draw[fill=black] (y) circle (1pt);
\draw[fill=black] (z) circle (1pt);

\draw (P) -- (w) -- (x) -- (y) -- (z) -- (R);
\draw (Q) -- (w) ;
\draw (Q) -- (x) ;
\draw (Q) -- (y) ;
\draw (Q) -- (z) ;
\draw (S) -- (w) ;
\draw (S) -- (x) ;
\draw (S) -- (y) ;
\draw (S) -- (z) ;

\ifcorners
    \draw (P) node[anchor=north east]{$\PP$};
    \draw (Q) node[anchor=north west]{$\QQ$};
    \draw (R) node[anchor=south west]{$\RR$};
    \draw (S) node[anchor=south east]{$\SS$};
\fi

\ifcolors
    \draw[line width=0,fill=\poscolor] (P) -- (Q) -- (w) -- cycle;
    \draw[line width=0,fill=\poscolor] (x) -- (Q) -- (y) -- cycle;
    \draw[line width=0,fill=\poscolor] (z) -- (Q) -- (R) -- cycle;
    \draw[line width=0,fill=\poscolor] (z) -- (S) -- (y) -- cycle;
    \draw[line width=0,fill=\poscolor] (x) -- (S) -- (w) -- cycle;
    \draw[line width=0,fill=\negcolor] (w) -- (Q) -- (x) -- cycle;
    \draw[line width=0,fill=\negcolor] (y) -- (Q) -- (z) -- cycle;
    \draw[line width=0,fill=\negcolor] (z) -- (R) -- (S) -- cycle;
    \draw[line width=0,fill=\negcolor] (y) -- (S) -- (x) -- cycle;
    \draw[line width=0,fill=\negcolor] (w) -- (S) -- (P) -- cycle;
\fi
    \end{tikzpicture}
     &  
    \begin{tikzpicture}[scale=.2]
    \coordinate (P) at (0,0);
\coordinate (Q) at (6,0);
\coordinate (R) at (6,6);
\coordinate (S) at (0,6);

\draw[fill=black] (P) circle (1pt);
\draw[fill=black] (Q) circle (1pt);
\draw[fill=black] (R) circle (1pt);
\draw[fill=black] (S) circle (1pt);
\draw (P) -- (Q) -- (R) -- (S) -- cycle;

\coordinate (w) at (1.5,2.5);
\coordinate (x) at (2.5,1.5);
\coordinate (y) at (3.5,3.5);
\coordinate (z) at (4.75,4.75);
\draw[fill=black] (w) circle (1pt);
\draw[fill=black] (x) circle (1pt);
\draw[fill=black] (y) circle (1pt);
\draw[fill=black] (z) circle (1pt);

\draw (y) -- (x) -- (P);
\draw (P) -- (w) -- (y) -- (z) -- (R) ;
\draw (w) -- (x);
\draw (S) -- (w);
\draw (S) -- (y);
\draw (S) -- (z);
\draw (Q) -- (x);
\draw (Q) -- (y);
\draw (Q) -- (z);

\ifcorners
    \draw (P) node[anchor=north east]{$\PP$};
    \draw (Q) node[anchor=north west]{$\QQ$};
    \draw (R) node[anchor=south west]{$\RR$};
    \draw (S) node[anchor=south east]{$\SS$};
\fi

\ifcolors
    \draw[line width=0,fill=\poscolor] (P) -- (Q) -- (x) -- cycle;
    \draw[line width=0,fill=\poscolor] (P) -- (w) -- (S) -- cycle;
    \draw[line width=0,fill=\poscolor] (w) -- (x) -- (y) -- cycle;
    \draw[line width=0,fill=\negcolor] (P) -- (x) -- (w) -- cycle;
    \draw[line width=0,fill=\negcolor] (x) -- (Q) -- (y) -- cycle;
    \draw[line width=0,fill=\negcolor] (y) -- (Q) -- (z) -- cycle;
    \draw[line width=0,fill=\negcolor] (z) -- (Q) -- (R) -- cycle;
    \draw[line width=0,fill=\negcolor] (z) -- (R) -- (S) -- cycle;
    \draw[line width=0,fill=\negcolor] (S) -- (y) -- (z) -- cycle;
    \draw[line width=0,fill=\negcolor] (S) -- (w) -- (y) -- cycle;
\fi
    \end{tikzpicture}
     &  
    \begin{tikzpicture}[scale=.2]
    \coordinate (P) at (0,0);
\coordinate (Q) at (6,0);
\coordinate (R) at (6,6);
\coordinate (S) at (0,6);

\draw[fill=black] (P) circle (1pt);
\draw[fill=black] (Q) circle (1pt);
\draw[fill=black] (R) circle (1pt);
\draw[fill=black] (S) circle (1pt);
\draw (P) -- (Q) -- (R) -- (S) -- cycle;

\coordinate (w) at (1.75,2.25);
\coordinate (x) at (3.75,1.75);
\coordinate (y) at (4.25,3.75);
\coordinate (z) at (2.25,4.25);
\draw[fill=black] (w) circle (1pt);
\draw[fill=black] (x) circle (1pt);
\draw[fill=black] (y) circle (1pt);
\draw[fill=black] (z) circle (1pt);

\draw (w) -- (x) -- (y) -- (z) -- cycle;
\draw (w) -- (y);

\draw (P) -- (w);
\draw (P) -- (x);
\draw (Q) -- (x);
\draw (Q) -- (y);
\draw (R) -- (y);
\draw (R) -- (z);
\draw (S) -- (z);
\draw (S) -- (w);

\ifcorners
    \draw (P) node[anchor=north east]{$\PP$};
    \draw (Q) node[anchor=north west]{$\QQ$};
    \draw (R) node[anchor=south west]{$\RR$};
    \draw (S) node[anchor=south east]{$\SS$};
\fi

\ifcolors
    \draw[line width=0,fill=\poscolor] (P) -- (Q) -- (x) -- cycle;
    \draw[line width=0,fill=\poscolor] (x) -- (Q) -- (y) -- cycle;
    \draw[line width=0,fill=\poscolor] (y) -- (R) -- (z) -- cycle;
    \draw[line width=0,fill=\poscolor] (P) -- (x) -- (w) -- cycle;
    \draw[line width=0,fill=\poscolor] (w) -- (x) -- (y) -- cycle;
    \draw[line width=0,fill=\poscolor] (w) -- (y) -- (z) -- cycle;
    \draw[line width=0,fill=\poscolor] (w) -- (z) -- (S) -- cycle;
    \draw[line width=0,fill=\poscolor] (z) -- (R) -- (S) -- cycle;
    \draw[line width=0,fill=\negcolor] (Q) -- (R) -- (y) -- cycle;
    \draw[line width=0,fill=\negcolor] (S) -- (P) -- (w) -- cycle;
\fi
    \end{tikzpicture}
     &  
    \begin{tikzpicture}[scale=.2]
    \coordinate (P) at (0,0);
\coordinate (Q) at (6,0);
\coordinate (R) at (6,6);
\coordinate (S) at (0,6);

\draw[fill=black] (P) circle (1pt);
\draw[fill=black] (Q) circle (1pt);
\draw[fill=black] (R) circle (1pt);
\draw[fill=black] (S) circle (1pt);
\draw (P) -- (Q) -- (R) -- (S) -- cycle;

\coordinate (w) at (1.75,2);
\coordinate (x) at (4.25,1.5);
\coordinate (y) at (4.25,4);
\coordinate (z) at (1.75,4.5);
\draw[fill=black] (w) circle (1pt);
\draw[fill=black] (x) circle (1pt);
\draw[fill=black] (y) circle (1pt);
\draw[fill=black] (z) circle (1pt);

\draw (w) -- (x) -- (y) -- (z) -- cycle;
\draw (w) -- (y);

\draw (P) -- (z);
\draw (P) -- (w);
\draw (P) -- (x);
\draw (Q) -- (x);
\draw (R) -- (x);
\draw (R) -- (y);
\draw (R) -- (z);
\draw (S) -- (z);

\ifcorners
    \draw (P) node[anchor=north east]{$\PP$};
    \draw (Q) node[anchor=north west]{$\QQ$};
    \draw (R) node[anchor=south west]{$\RR$};
    \draw (S) node[anchor=south east]{$\SS$};
\fi

\ifcolors
    \draw[line width=0,fill=\poscolor] (P) -- (Q) -- (x) -- cycle;
    \draw[line width=0,fill=\poscolor] (P) -- (x) -- (w) -- cycle;
    \draw[line width=0,fill=\poscolor] (x) -- (Q) -- (R) -- cycle;
    \draw[line width=0,fill=\poscolor] (x) -- (R) -- (y) -- cycle;
    \draw[line width=0,fill=\poscolor] (w) -- (y) -- (z) -- cycle;
    \draw[line width=0,fill=\negcolor] (w) -- (x) -- (y) -- cycle;
    \draw[line width=0,fill=\negcolor] (y) -- (R) -- (z) -- cycle;
    \draw[line width=0,fill=\negcolor] (z) -- (R) -- (S) -- cycle;
    \draw[line width=0,fill=\negcolor] (S) -- (w) -- (z) -- cycle;
    \draw[line width=0,fill=\negcolor] (S) -- (P) -- (w) -- cycle;
\fi
    \end{tikzpicture}
     &  
    \begin{tikzpicture}[scale=.2]
    \coordinate (P) at (0,0);
\coordinate (Q) at (6,0);
\coordinate (R) at (6,6);
\coordinate (S) at (0,6);

\draw[fill=black] (P) circle (1pt);
\draw[fill=black] (Q) circle (1pt);
\draw[fill=black] (R) circle (1pt);
\draw[fill=black] (S) circle (1pt);
\draw (P) -- (Q) -- (R) -- (S) -- cycle;

\coordinate (w) at (2,2.25);
\coordinate (x) at (4.25,1.75);
\coordinate (y) at (4.5,4.25);
\coordinate (z) at (1.75,4.25);
\draw[fill=black] (w) circle (1pt);
\draw[fill=black] (x) circle (1pt);
\draw[fill=black] (y) circle (1pt);
\draw[fill=black] (z) circle (1pt);

\draw (w) -- (x) -- (y) -- (z) -- cycle;
\draw (w) -- (y);

\draw (P) -- (z);
\draw (P) -- (w);
\draw (P) -- (x);
\draw (Q) -- (x);
\draw (Q) -- (y);
\draw (R) -- (y);
\draw (S) -- (y);
\draw (S) -- (z);

\ifcorners
    \draw (P) node[anchor=north east]{$\PP$};
    \draw (Q) node[anchor=north west]{$\QQ$};
    \draw (R) node[anchor=south west]{$\RR$};
    \draw (S) node[anchor=south east]{$\SS$};
\fi

\ifcolors
    \draw[line width=0,fill=\poscolor] (P) -- (Q) -- (x) -- cycle;
    \draw[line width=0,fill=\poscolor] (Q) -- (R) -- (y) -- cycle;
    \draw[line width=0,fill=\poscolor] (w) -- (x) -- (y) -- cycle;
    \draw[line width=0,fill=\poscolor] (y) -- (S) -- (z) -- cycle;
    \draw[line width=0,fill=\poscolor] (z) -- (P) -- (w) -- cycle;
    \draw[line width=0,fill=\negcolor] (P) -- (x) -- (w) -- cycle;
    \draw[line width=0,fill=\negcolor] (x) -- (Q) -- (y) -- cycle;
    \draw[line width=0,fill=\negcolor] (y) -- (R) -- (S) -- cycle;
    \draw[line width=0,fill=\negcolor] (w) -- (y) -- (z) -- cycle;
    \draw[line width=0,fill=\negcolor] (S) -- (P) -- (z) -- cycle;
\fi
    \end{tikzpicture}
     &  
    \begin{tikzpicture}[scale=.2]
    \coordinate (P) at (0,0);
\coordinate (Q) at (6,0);
\coordinate (R) at (6,6);
\coordinate (S) at (0,6);

\draw[fill=black] (P) circle (1pt);
\draw[fill=black] (Q) circle (1pt);
\draw[fill=black] (R) circle (1pt);
\draw[fill=black] (S) circle (1pt);
\draw (P) -- (Q) -- (R) -- (S) -- cycle;

\coordinate (w) at (1.75,2);
\coordinate (x) at (4.25,1.5);
\coordinate (y) at (4.25,4);
\coordinate (z) at (1.75,4.5);
\draw[fill=black] (w) circle (1pt);
\draw[fill=black] (x) circle (1pt);
\draw[fill=black] (y) circle (1pt);
\draw[fill=black] (z) circle (1pt);

\draw (w) -- (x) -- (y) -- (z) -- cycle;
\draw (w) -- (y);

\draw (P) -- (w);
\draw (P) -- (x);
\draw (Q) -- (x);
\draw (Q) -- (y);
\draw (R) -- (y);
\draw (R) -- (z);
\draw (S) -- (z);
\draw (P) -- (z);

\ifcorners
    \draw (P) node[anchor=north east]{$\PP$};
    \draw (Q) node[anchor=north west]{$\QQ$};
    \draw (R) node[anchor=south west]{$\RR$};
    \draw (S) node[anchor=south east]{$\SS$};
\fi

\ifcolors
    \draw[line width=0,fill=\poscolor] (P) -- (Q) -- (x) -- cycle;
    \draw[line width=0,fill=\poscolor] (P) -- (w) -- (z) -- cycle;
    \draw[line width=0,fill=\poscolor] (w) -- (x) -- (y) -- cycle;
    \draw[line width=0,fill=\poscolor] (z) -- (R) -- (S) -- cycle;
    \draw[line width=0,fill=\negcolor] (P) -- (x) -- (w) -- cycle;
    \draw[line width=0,fill=\negcolor] (Q) -- (y) -- (x) -- cycle;
    \draw[line width=0,fill=\negcolor] (Q) -- (R) -- (y) -- cycle;
    \draw[line width=0,fill=\negcolor] (R) -- (z) -- (y) -- cycle;
    \draw[line width=0,fill=\negcolor] (w) -- (y) -- (z) -- cycle;
    \draw[line width=0,fill=\negcolor] (S) -- (P) -- (z) -- cycle;
\fi
    \end{tikzpicture}
    \\
                    $T_4$ & $T_{3,1}$ & $T_{2,2,2,2}$ & $T_{3,1,3,1}$ & $T_{3,2,1,2}$ & $T_{3,2,2,1}$\\
                    \addlinespace[5pt]
                    \bottomrule
                    \addlinespace[5pt]
                \end{tabular}
                \caption{Triangulations with up to 4 interior vertices and no subdivisions.}
                \label{table:T.up.to.4.int}
            \end{table}

        We will be interested in drawings $\rho$ that have a fixed number $l$ of nondegenerate triangles.  The triangulation being drawn may have many degenerate triangles, but we show that the nondegenerate triangles of such a drawing can be reinterpreted as the nondegenerate triangles in the image of a \emph{simple} drawing, which will automatically have a bounded number of triangles.

        For the following we use the language of elasticity (see Definition \ref{def:cellfiber}).
        \begin{defn}[Simple drawing]\label{def:simple}
            Let $T$ be a triangulation and $\rho$ a drawing of $T$.  We call $\rho$ \emph{simple} if all of the following are true:
            \begin{enumerate}
                \item[(i)\ ] every subdivision of $T$ contains at least one nondegenerate triangle of $\rho$;
                \item[(ii)\ ] there is no interior edge $e=vw$ that is elastic and that can be contracted;
                \item[(iii)\ ] every interior vertex $v$ is contained in at least one nondegenerate triangle of $\rho$.
            \end{enumerate}
        \end{defn}
        
        Note that the set of simple drawings is open in $X(T)$.
        
        \begin{lem}\label{lem:simplebound}
            Let $T$ be a triangulation and $\rho$ a simple drawing of $T$.  Suppose $\rho$ has $l>0$ nondegenerate triangles and that no (nonempty) subset of the nondegenerate triangles has areas summing to zero.  Then $T$ has at most $\frac{3}{2}l-2$ interior vertices and at most $3l-2$ triangles.
        \end{lem}
        
        \begin{proof}
            We show that our hypotheses imply the following:    
            \begin{enumerate}
                \item $\rho$ has no elastic edges, and
                \item every corner of $T$ is contained in at least one nondegenerate triangle of $\rho$, and
                \item every interior vertex of $T$ is contained in at least two nondegenerate triangles of $\rho$.
            \end{enumerate}
            The lemma follows by counting vertices with multiplicity:  there are $3l$ vertices of nondegenerate triangles in the image of $\rho$, so if $T$ has $k$ interior vertices then we have $3l\ge 4+2k$.  The number of triangles of $T'$ is exactly $2k+2$.

            To show (1), first observe that no boundary edge can be elastic, since this would force all areas to sum to zero.
            Suppose the interior edge $e=xy$ is elastic.  As $\rho$ is simple, the edge $e$ cannot be contracted, which means that $\{x,y\}$ is part of a subdivision.  Again using simplicity, there must be at least one nondegenerate triangle inside this subdivision.  However $e$ being elastic implies that all areas inside the subdivision sum to zero, a contradiction.
            
            To show (2) we suppose there is a corner, say $\PP$, that is not in any nondegenerate triangle of $\rho$.  By (1) there are no elastic edges, so $\rho$ maps $\PP$ and all its neighbors to a set of collinear points.  As $\QQ$ and $\SS$ are among the neighbors of $\PP$, this forces the sum of all the areas to be zero, a contradiction.
            
            To show (3) let $v$ be an interior vertex.  By simplicity $v$ is contained in at least one nondegenerate triangle of $\rho$.  Suppose $v$ is contained in exactly one nondegenerate triangle $\Delta$.  Then all triangles surrounding $v$ other than $\Delta$ are degenerate triangles of $\rho$, and by (1) there are no elastic edges, so all vertices of all these triangles are mapped collinearly by $\rho$.  But this causes $\Delta$ to be degenerate also, a contradiction.  This completes the proof.
        \end{proof}
        
        \begin{thm}[Finiteness of drawings]\label{thm:simplification}
            Let $T$ be a triangulation and let $\rho$ be a drawing of $T$ with $l>0$ nondegenerate triangles.  Suppose that no (nonempty) subset of nondegenerate triangles has areas summing to zero.  There is a triangulation $T'$ and a simple drawing $\rho'$ of $T'$ such that the nondegenerate triangles of $\rho'$ are the same as the nondegenerate triangles of $\rho$.
        \end{thm}
        
        \begin{proof}
            If $\rho$ is already simple then of course we are done.  We assume not and we describe a process for simplifying $\rho$.
            
            Suppose $T$ contains a subdivision with only degenerate triangles of $\rho$ inside the subdivision.  Then we may replace all triangles inside the subdivision with a single triangle, and we replace $\rho$ with its restriction $\rho'$ to the remaining vertices.  Note that $\rho'$ has the same nondegenerate triangles as $\rho$.  We assume this has been done for all such subdivisions.
            
            Next suppose that $\rho$ has an elastic edge $e=vw$.  Note that $e$ cannot connect two corners, because this would cause the total area of all nondegenerate triangles to be zero, contrary to hypothesis.  We claim that $\{v,w\}$ is not part of a subdivision:  if it were, then there could not be any nondegenerate triangles inside the subdivision, since their areas would sum to the area of the subdivided triangle, which is zero by elasticity of $vw$; yet, we have already eliminated all subdivisions containing only degenerate triangles.  So $\{v,w\}$ is not part of a subdivision, and we may contract the edge obtaining the triangulation $T|_e$. Note that the same map $\rho$ defines a drawing of $T|_e$, and the nondegenerate triangles have not changed.  
            
            It is possible that contracting an edge creates a subdivision.  If all triangles inside this subdivision are degenerate then we can return to the previous step and simplify the situation.
            
            Repeat this process until no elastic edges remain.
             
            Finally we argue that we can delete any vertices of $T$ that are part of only degenerate triangles. Let $v$ be such a vertex, with neighbors $w_i$ labeled in cyclic order around $v$.  Note that $v$ cannot be a corner:  if it were, then since we have eliminated elasticity, $v$ and all its neighbors (including two other corners) would be mapped by $\rho$ to a set of collinear points, forcing the sum of all areas of all triangles in the whole triangulation to be zero, contrary to hypothesis.  Thus the $w_i$ form a polygon $\pentagon$ enclosing $v$.  Since there are no elastic edges, it follows that there is a single line $\ell$ that contains (the images of) $v$ and all $w_i$. 
            
            We would like to delete $v$; the only thing to check is that we can re-triangulate the polygon $\pentagon$ in such a way that the result gives a new triangulation of the square.  Indeed the only way this could fail is if we create an edge inside $\pentagon$ joining vertices $w_i$ and $w_j$ that are already connected by an edge outside $\pentagon$.  If there is no edge $w_iw_j$ outside $\pentagon$ then we may triangulate $\pentagon$ arbitrarily.  Otherwise suppose $w_i$ is connected to $w_j$ outside $\pentagon$, with $|i-j|$ minimal.  Then we may choose $w$ between $w_i$ and $w_j$ and triangulate $\pentagon$ by adding edges from $w$ to every other $w_i$.  This results in a new triangulation $T'$ of the square.  
            
            Restricting $\rho$ to the vertices of $T'$ produces a drawing $\rho'$ of the new triangulation.  Observe that all of the newly defined triangles inside $\pentagon$ are degenerate, because all of the $w_i$ lie on the line $\ell$.  Thus $\rho'$ has the same nondegenerate triangles as $\rho$.  
            
            It is possible that this step produces subdivisions and/or elasticity, in which case we return to the previous cases as needed.
            Each time we need to modify $T$, the number of vertices goes down, so this cannot go on forever.  We repeat this process until every vertex is part of at least one nondegenerate triangle.
            
            After all this we have a new triangulation $T'$ and a new drawing $\rho'$ which satisfies the conclusion of the theorem.
        \end{proof}

        Drawings with one or two nondegenerate triangles are addressed by the following corollary.
        
        \begin{cor}\label{cor:oneortwo} Let $T$ be a triangulation. 
            \begin{enumerate}
                \item There is no drawing of $T$ with exactly one nondegenerate triangle.
                \item Suppose $T$ has a drawing  with exactly two nondegenerate triangles of areas $a$ and $b$.  Then $a=\pm b$.
            \end{enumerate}
        \end{cor}
        
        \begin{proof}
            If $l=1$, then $\frac32 l - 2 <0$, proving (1).
            
            If $l=2$, then there is a simple drawing of a triangulation $T'$ that has at most one interior vertex.   By Table \ref{table:T.up.to.4.int}, $T'$ must be either $T_0$ or $T_1$.  In the first case since $p_{T_0}=A-B$ we must have $a=b$; in the second case two of the variables of $p_{T_1}=A-B+C-D$ must be zero so the only possibilities are $a=\pm b$.
        \end{proof}
    
     The first part of this corollary may seem obvious from a geometric point of view.  However, it is not completely frivolous since we  allow our drawings to be complex (for instance see Example \ref{ex:complexequi}).  The statement implies that for any $T$, the point $[1:0:\dots:0]$ is not in the image of the area map.

    \section{Coefficients and a canonical 2-coloring of triangles}
    
        We illustrate the power of Theorem \ref{thm:characterization} with two observations about points of $V$ and the coefficients of area polynomials.  The first does not come as a surprise, and in fact it is weaker than Theorem \ref{thm:mod2honest}, which was first proved in \cite{chapter3} and which immediately implies that all leading coefficients of $p_T$ are odd.  We state and prove it here anyway, as it is a useful fact to keep in mind.  Our second observation, on the other hand, is new and improves significantly on Theorem \ref{thm:mod2honest}.  It also points to a canonical 2-coloring of the triangles in any triangulation.

        \begin{thm}[Triangulum non quadratum est] \label{thm:encyc1}
            Let $T$ be a triangulation.  Then \\
            (a) $[1:0:\cdots:0]\notin V(T)$, and \\
            (b) the leading coefficients of $p_{T}$ are nonzero.
        \end{thm}
        
        \begin{proof}
          Assume $w=[1:0:\cdots:0]\in V(T)$.  Any subset of the coordinates that sums to $0$ must consist entirely of $0$'s, so Theorem \ref{thm:characterization} implies that $w$ is in the image of the area map for $T$.  This contradicts Corollary \ref{cor:oneortwo}.
          
          Assertion (b) follows from (a). 
        \end{proof}
        
        \begin{thm}\label{thm:encyc2}
            Let $T$ be a triangulation.  Then \\
            (a) if $[a:b:0:\cdots:0]\in V(T)$ with $a,b\neq 0$, then $a=\pm b$, and \\
            (b) the polynomial $p_T(A_1, A_2, 0, \dots, 0)$ has the form $c(A_1+A_2)^e(A_1-A_2)^f$.  
            
            In particular, all leading coefficients of $p(T)$ are equal up to sign. 
        \end{thm}
        
        \begin{proof}
            Assume $w=[a:b:\cdots:0]\in V(T)$ with $a, b\neq 0$ and assume $a\neq -b$. Then by Theorem \ref{thm:characterization}, it follows that $w$ is in the image of the area map for $T$.  By Corollary \ref{cor:oneortwo}, it follows that $a=\pm b$. 
          
            The assertion (b) now follows straightaway.  If $q$ is an irreducible factor of $p_T(A_1, A_2, 0, \dots, 0)$ then the its zero set of $q$ is contained in the union of the zero sets of $A_1+A_2$ and $A_1-A_2$. Thus $q$ must equal one of $A_1+A_2$ or $A_1-A_2$ up to scalar multiples. 
        \end{proof}
        
        We do not know if the leading coefficients of $p_{T}$ are equal to $\pm 1$ for every $T$; this has held for every example we have calculated.  This is a special case of the positivity phenomenon discussed in \cite{chapter3}.
        
        Theorem \ref{thm:encyc2} points to a natural $2$-coloring defined on the triangles of $T$: each triangle is assigned a color based on the sign of the corresponding leading coefficient in $p_T$.  We call this the \emph{canonical 2-coloring} of $T$.  We note that if $T$ is the diagonal triangulation $T_n$, then by work carried out in \cite{triangles1} the canonical 2-coloring of $T$ is the checkerboard coloring.  Table \ref{table:2color} shows the triangulations from Table \ref{table:T.up.to.4.int} with their canonical 2-colorings. 
        
        \begin{question}
            Is there an efficient algorithm for determining the canonical 2-coloring of a triangulation $T$?  If $T$ is a triangulation that admits a checkerboard coloring, is this coloring necessarily the same as the canonical coloring?
        \end{question} 
        
        Of course, two triangles have the same or different colors according to whether the exponent $f$ in Theorem \ref{thm:encyc2}(b) is even or odd.
                
        Corollary \ref{cor:EimpliesT} in the next section (also Main Theorem \ref{mainthm:nonconverse}) allows one to carry out a similar analysis of points of $V$ that have more than two nonzero coordinates.

        \colorstrue
    
        \begin{table}%\label{table:2color}
            \begin{tabular}{p{.7in}p{.7in}p{.7in}p{.7in}p{.7in}p{.7in}}
            \toprule
            \addlinespace[5pt]
                  
    \begin{tikzpicture}[scale=.2]
    \coordinate (P) at (0,0);
\coordinate (Q) at (6,0);
\coordinate (R) at (6,6);
\coordinate (S) at (0,6);

\draw[fill=black] (P) circle (1pt);
\draw[fill=black] (Q) circle (1pt);
\draw[fill=black] (R) circle (1pt);
\draw[fill=black] (S) circle (1pt);
\draw (P) -- (Q) -- (R) -- (S) -- cycle;

\draw (P) -- (R);

\ifcorners
    \draw (P) node[anchor=north east]{$\PP$};
    \draw (Q) node[anchor=north west]{$\QQ$};
    \draw (R) node[anchor=south west]{$\RR$};
    \draw (S) node[anchor=south east]{$\SS$};
\fi

\ifareas
    \draw (4,2) node{$A$};
    \draw (1.8,4) node{$B$};
\fi

\ifnewareas
    \draw (4,2) node{\scriptsize $A$};
    \draw (1.8,4) node{\scriptsize $B$};
\fi

\ifcolors
    \draw[line width=0, fill=\poscolor] (P) -- (Q) -- (R) -- cycle;
    \draw[line width=0, fill=\negcolor] (R) -- (S) -- (P) -- cycle;
\fi
    \end{tikzpicture}
     &  
    \begin{tikzpicture}[scale=.2]
    \coordinate (P) at (0,0);
\coordinate (Q) at (6,0);
\coordinate (R) at (6,6);
\coordinate (S) at (0,6);

\draw[fill=black] (P) circle (1pt);
\draw[fill=black] (Q) circle (1pt);
\draw[fill=black] (R) circle (1pt);
\draw[fill=black] (S) circle (1pt);
\draw (P) -- (Q) -- (R) -- (S) -- cycle;

\coordinate (x) at (3,3);
\draw[fill=black] (x) circle (1pt);

\draw (P) -- (x);
\draw (Q) -- (x);
\draw (R) -- (x);
\draw (S) -- (x);

\ifcorners
    \draw (P) node[anchor=north east]{$\PP$};
    \draw (Q) node[anchor=north west]{$\QQ$};
    \draw (R) node[anchor=south west]{$\RR$};
    \draw (S) node[anchor=south east]{$\SS$};
\fi

\ifcolors
    \draw[line width=0, fill=\poscolor] (P) -- (Q) -- (x) -- cycle;
    \draw[line width=0, fill=\poscolor] (R) -- (S) -- (x) -- cycle;
    \draw[line width=0, fill=\negcolor] (Q) -- (R) -- (x) -- cycle;
    \draw[line width=0, fill=\negcolor] (S) -- (P) -- (x) -- cycle;
\fi

\ifareas
    \draw (3,1.25) node{$A$};
    \draw (4.75,3) node{$B$};
    \draw (3,4.75) node{$C$};
    \draw (1.25,3) node{$D$};
\fi

\ifnewareas
    \draw (3,1.25) node{\scriptsize $D$};
    \draw (4.75,3) node{\scriptsize $B$};
    \draw (3,4.75) node{\scriptsize $C$};
    \draw (1.25,3) node{\scriptsize $A$};
\fi
    \end{tikzpicture}
     &  
    \begin{tikzpicture}[scale=.2]
    \coordinate (P) at (0,0);
\coordinate (Q) at (6,0);
\coordinate (R) at (6,6);
\coordinate (S) at (0,6);

\draw[fill=black] (P) circle (1pt);
\draw[fill=black] (Q) circle (1pt);
\draw[fill=black] (R) circle (1pt);
\draw[fill=black] (S) circle (1pt);
\draw (P) -- (Q) -- (R) -- (S) -- cycle;

\coordinate (u) at (2,2.5);
\coordinate (v) at (4.5,3.75);
\draw[fill=black] (u) circle (1pt);
\draw[fill=black] (v) circle (1pt);

\draw (P) -- (u) -- (v) -- (R);
\draw (Q) -- (u) -- (S);
\draw (S) -- (v) -- (Q);

\ifcorners
    \draw (P) node[anchor=north east]{$\PP$};
    \draw (Q) node[anchor=north west]{$\QQ$};
    \draw (R) node[anchor=south west]{$\RR$};
    \draw (S) node[anchor=south east]{$\SS$};
\fi

\ifareas
    \draw (2.25,1) node{$A$};
    \draw (4,2.25) node{$B$};
    \draw (5.5,2.75) node{$C$};
    \draw (.75,2.75) node{$D$};
    \draw (2.25,4) node{$E$};
    \draw (4,5) node{$F$};
\fi

\ifinterior
    \draw (u) node[anchor=north]{$u$};
    \draw (v) node[anchor=west]{$v$};
\fi

\ifcolors
    \draw[line width=0, fill=\poscolor] (P) -- (Q) -- (u) -- cycle;
    \draw[line width=0, fill=\poscolor] (v) -- (Q) -- (R) -- cycle;
    \draw[line width=0, fill=\poscolor] (u) -- (v) -- (S) -- cycle;
    \draw[line width=0, fill=\negcolor] (u) -- (Q) -- (v) -- cycle;
    \draw[line width=0, fill=\negcolor] (v) -- (R) -- (S) -- cycle;
    \draw[line width=0, fill=\negcolor] (u) -- (S) -- (P) -- cycle;
\fi
    \end{tikzpicture}
     &  
    \begin{tikzpicture}[scale=.2]
    \coordinate (P) at (0,0);
\coordinate (Q) at (6,0);
\coordinate (R) at (6,6);
\coordinate (S) at (0,6);

\draw[fill=black] (P) circle (1pt);
\draw[fill=black] (Q) circle (1pt);
\draw[fill=black] (R) circle (1pt);
\draw[fill=black] (S) circle (1pt);
\draw (P) -- (Q) -- (R) -- (S) -- cycle;

\coordinate (w) at (1.5,1.5);
\coordinate (x) at (3,3);
\coordinate (y) at (4.5,4.5);
\draw[fill=black] (w) circle (1pt);
\draw[fill=black] (x) circle (1pt);
\draw[fill=black] (y) circle (1pt);

\draw (P) -- (w) -- (x) -- (y) -- (R);
\draw (Q) -- (w) ;
\draw (Q) -- (x) ;
\draw (Q) -- (y) ;
\draw (S) -- (w) ;
\draw (S) -- (x) ;
\draw (S) -- (y) ;

\ifcorners
    \draw (P) node[anchor=north east]{$\PP$};
    \draw (Q) node[anchor=north west]{$\QQ$};
    \draw (R) node[anchor=south west]{$\RR$};
    \draw (S) node[anchor=south east]{$\SS$};
\fi

\ifcolors
    \draw[line width=0, fill=\poscolor] (P) -- (Q) -- (w) -- cycle;
    \draw[line width=0, fill=\poscolor] (x) -- (Q) -- (y) -- cycle;
    \draw[line width=0, fill=\poscolor] (R) -- (S) -- (y) -- cycle;
    \draw[line width=0, fill=\poscolor] (x) -- (S) -- (w) -- cycle;
    \draw[line width=0, fill=\negcolor] (w) -- (Q) -- (x) -- cycle;
    \draw[line width=0, fill=\negcolor] (y) -- (Q) -- (R) -- cycle;
    \draw[line width=0, fill=\negcolor] (y) -- (S) -- (x) -- cycle;
    \draw[line width=0, fill=\negcolor] (w) -- (S) -- (P) -- cycle;
\fi
    \end{tikzpicture}
     &  
    \begin{tikzpicture}[scale=.2]
    \coordinate (P) at (0,0);
\coordinate (Q) at (6,0);
\coordinate (R) at (6,6);
\coordinate (S) at (0,6);

\draw[fill=black] (P) circle (1pt);
\draw[fill=black] (Q) circle (1pt);
\draw[fill=black] (R) circle (1pt);
\draw[fill=black] (S) circle (1pt);
\draw (P) -- (Q) -- (R) -- (S) -- cycle;

\coordinate (w) at (3,1.6);
\coordinate (x) at (1.6,3);
\coordinate (y) at (4,4);
\draw[fill=black] (w) circle (1pt);
\draw[fill=black] (x) circle (1pt);
\draw[fill=black] (y) circle (1pt);

\draw (w) -- (y) -- (x) -- cycle;
\draw (Q) -- (w) -- (P);
\draw (Q) -- (y);
\draw (P) -- (x) -- (S);
\draw (S) -- (y) -- (R) ;

\ifcorners
    \draw (P) node[anchor=north east]{$\PP$};
    \draw (Q) node[anchor=north west]{$\QQ$};
    \draw (R) node[anchor=south west]{$\RR$};
    \draw (S) node[anchor=south east]{$\SS$};
\fi

\ifcolors
    \draw[line width=0, fill=\poscolor] (P) -- (Q) -- (w) -- cycle;
    \draw[line width=0, fill=\poscolor] (w) -- (x) -- (y) -- cycle;
    \draw[line width=0, fill=\poscolor] (x) -- (S) -- (P) -- cycle;
    \draw[line width=0, fill=\negcolor] (P) -- (w) -- (x) -- cycle;
    \draw[line width=0, fill=\negcolor] (w) -- (Q) -- (y) -- cycle;
    \draw[line width=0, fill=\negcolor] (y) -- (Q) -- (R) -- cycle;
    \draw[line width=0, fill=\negcolor] (R) -- (S) -- (y) -- cycle;
    \draw[line width=0, fill=\negcolor] (x) -- (S) -- (y) -- cycle;
\fi
    \end{tikzpicture}
     & \\
            %\midrule
            \addlinespace[5pt]
                  
    \begin{tikzpicture}[scale=.2]
    \coordinate (P) at (0,0);
\coordinate (Q) at (6,0);
\coordinate (R) at (6,6);
\coordinate (S) at (0,6);

\draw[fill=black] (P) circle (1pt);
\draw[fill=black] (Q) circle (1pt);
\draw[fill=black] (R) circle (1pt);
\draw[fill=black] (S) circle (1pt);
\draw (P) -- (Q) -- (R) -- (S) -- cycle;

\coordinate (w) at (1.2,1.2);
\coordinate (x) at (2.4,2.4);
\coordinate (y) at (3.6,3.6);
\coordinate (z) at (4.8,4.8);
\draw[fill=black] (w) circle (1pt);
\draw[fill=black] (x) circle (1pt);
\draw[fill=black] (y) circle (1pt);
\draw[fill=black] (z) circle (1pt);

\draw (P) -- (w) -- (x) -- (y) -- (z) -- (R);
\draw (Q) -- (w) ;
\draw (Q) -- (x) ;
\draw (Q) -- (y) ;
\draw (Q) -- (z) ;
\draw (S) -- (w) ;
\draw (S) -- (x) ;
\draw (S) -- (y) ;
\draw (S) -- (z) ;

\ifcorners
    \draw (P) node[anchor=north east]{$\PP$};
    \draw (Q) node[anchor=north west]{$\QQ$};
    \draw (R) node[anchor=south west]{$\RR$};
    \draw (S) node[anchor=south east]{$\SS$};
\fi

\ifcolors
    \draw[line width=0,fill=\poscolor] (P) -- (Q) -- (w) -- cycle;
    \draw[line width=0,fill=\poscolor] (x) -- (Q) -- (y) -- cycle;
    \draw[line width=0,fill=\poscolor] (z) -- (Q) -- (R) -- cycle;
    \draw[line width=0,fill=\poscolor] (z) -- (S) -- (y) -- cycle;
    \draw[line width=0,fill=\poscolor] (x) -- (S) -- (w) -- cycle;
    \draw[line width=0,fill=\negcolor] (w) -- (Q) -- (x) -- cycle;
    \draw[line width=0,fill=\negcolor] (y) -- (Q) -- (z) -- cycle;
    \draw[line width=0,fill=\negcolor] (z) -- (R) -- (S) -- cycle;
    \draw[line width=0,fill=\negcolor] (y) -- (S) -- (x) -- cycle;
    \draw[line width=0,fill=\negcolor] (w) -- (S) -- (P) -- cycle;
\fi
    \end{tikzpicture}
     &  
    \begin{tikzpicture}[scale=.2]
    \coordinate (P) at (0,0);
\coordinate (Q) at (6,0);
\coordinate (R) at (6,6);
\coordinate (S) at (0,6);

\draw[fill=black] (P) circle (1pt);
\draw[fill=black] (Q) circle (1pt);
\draw[fill=black] (R) circle (1pt);
\draw[fill=black] (S) circle (1pt);
\draw (P) -- (Q) -- (R) -- (S) -- cycle;

\coordinate (w) at (1.5,2.5);
\coordinate (x) at (2.5,1.5);
\coordinate (y) at (3.5,3.5);
\coordinate (z) at (4.75,4.75);
\draw[fill=black] (w) circle (1pt);
\draw[fill=black] (x) circle (1pt);
\draw[fill=black] (y) circle (1pt);
\draw[fill=black] (z) circle (1pt);

\draw (y) -- (x) -- (P);
\draw (P) -- (w) -- (y) -- (z) -- (R) ;
\draw (w) -- (x);
\draw (S) -- (w);
\draw (S) -- (y);
\draw (S) -- (z);
\draw (Q) -- (x);
\draw (Q) -- (y);
\draw (Q) -- (z);

\ifcorners
    \draw (P) node[anchor=north east]{$\PP$};
    \draw (Q) node[anchor=north west]{$\QQ$};
    \draw (R) node[anchor=south west]{$\RR$};
    \draw (S) node[anchor=south east]{$\SS$};
\fi

\ifcolors
    \draw[line width=0,fill=\poscolor] (P) -- (Q) -- (x) -- cycle;
    \draw[line width=0,fill=\poscolor] (P) -- (w) -- (S) -- cycle;
    \draw[line width=0,fill=\poscolor] (w) -- (x) -- (y) -- cycle;
    \draw[line width=0,fill=\negcolor] (P) -- (x) -- (w) -- cycle;
    \draw[line width=0,fill=\negcolor] (x) -- (Q) -- (y) -- cycle;
    \draw[line width=0,fill=\negcolor] (y) -- (Q) -- (z) -- cycle;
    \draw[line width=0,fill=\negcolor] (z) -- (Q) -- (R) -- cycle;
    \draw[line width=0,fill=\negcolor] (z) -- (R) -- (S) -- cycle;
    \draw[line width=0,fill=\negcolor] (S) -- (y) -- (z) -- cycle;
    \draw[line width=0,fill=\negcolor] (S) -- (w) -- (y) -- cycle;
\fi
    \end{tikzpicture}
     &  
    \begin{tikzpicture}[scale=.2]
    \coordinate (P) at (0,0);
\coordinate (Q) at (6,0);
\coordinate (R) at (6,6);
\coordinate (S) at (0,6);

\draw[fill=black] (P) circle (1pt);
\draw[fill=black] (Q) circle (1pt);
\draw[fill=black] (R) circle (1pt);
\draw[fill=black] (S) circle (1pt);
\draw (P) -- (Q) -- (R) -- (S) -- cycle;

\coordinate (w) at (1.75,2.25);
\coordinate (x) at (3.75,1.75);
\coordinate (y) at (4.25,3.75);
\coordinate (z) at (2.25,4.25);
\draw[fill=black] (w) circle (1pt);
\draw[fill=black] (x) circle (1pt);
\draw[fill=black] (y) circle (1pt);
\draw[fill=black] (z) circle (1pt);

\draw (w) -- (x) -- (y) -- (z) -- cycle;
\draw (w) -- (y);

\draw (P) -- (w);
\draw (P) -- (x);
\draw (Q) -- (x);
\draw (Q) -- (y);
\draw (R) -- (y);
\draw (R) -- (z);
\draw (S) -- (z);
\draw (S) -- (w);

\ifcorners
    \draw (P) node[anchor=north east]{$\PP$};
    \draw (Q) node[anchor=north west]{$\QQ$};
    \draw (R) node[anchor=south west]{$\RR$};
    \draw (S) node[anchor=south east]{$\SS$};
\fi

\ifcolors
    \draw[line width=0,fill=\poscolor] (P) -- (Q) -- (x) -- cycle;
    \draw[line width=0,fill=\poscolor] (x) -- (Q) -- (y) -- cycle;
    \draw[line width=0,fill=\poscolor] (y) -- (R) -- (z) -- cycle;
    \draw[line width=0,fill=\poscolor] (P) -- (x) -- (w) -- cycle;
    \draw[line width=0,fill=\poscolor] (w) -- (x) -- (y) -- cycle;
    \draw[line width=0,fill=\poscolor] (w) -- (y) -- (z) -- cycle;
    \draw[line width=0,fill=\poscolor] (w) -- (z) -- (S) -- cycle;
    \draw[line width=0,fill=\poscolor] (z) -- (R) -- (S) -- cycle;
    \draw[line width=0,fill=\negcolor] (Q) -- (R) -- (y) -- cycle;
    \draw[line width=0,fill=\negcolor] (S) -- (P) -- (w) -- cycle;
\fi
    \end{tikzpicture}
     &  
    \begin{tikzpicture}[scale=.2]
    \coordinate (P) at (0,0);
\coordinate (Q) at (6,0);
\coordinate (R) at (6,6);
\coordinate (S) at (0,6);

\draw[fill=black] (P) circle (1pt);
\draw[fill=black] (Q) circle (1pt);
\draw[fill=black] (R) circle (1pt);
\draw[fill=black] (S) circle (1pt);
\draw (P) -- (Q) -- (R) -- (S) -- cycle;

\coordinate (w) at (1.75,2);
\coordinate (x) at (4.25,1.5);
\coordinate (y) at (4.25,4);
\coordinate (z) at (1.75,4.5);
\draw[fill=black] (w) circle (1pt);
\draw[fill=black] (x) circle (1pt);
\draw[fill=black] (y) circle (1pt);
\draw[fill=black] (z) circle (1pt);

\draw (w) -- (x) -- (y) -- (z) -- cycle;
\draw (w) -- (y);

\draw (P) -- (z);
\draw (P) -- (w);
\draw (P) -- (x);
\draw (Q) -- (x);
\draw (R) -- (x);
\draw (R) -- (y);
\draw (R) -- (z);
\draw (S) -- (z);

\ifcorners
    \draw (P) node[anchor=north east]{$\PP$};
    \draw (Q) node[anchor=north west]{$\QQ$};
    \draw (R) node[anchor=south west]{$\RR$};
    \draw (S) node[anchor=south east]{$\SS$};
\fi

\ifcolors
    \draw[line width=0,fill=\poscolor] (P) -- (Q) -- (x) -- cycle;
    \draw[line width=0,fill=\poscolor] (P) -- (x) -- (w) -- cycle;
    \draw[line width=0,fill=\poscolor] (x) -- (Q) -- (R) -- cycle;
    \draw[line width=0,fill=\poscolor] (x) -- (R) -- (y) -- cycle;
    \draw[line width=0,fill=\poscolor] (w) -- (y) -- (z) -- cycle;
    \draw[line width=0,fill=\negcolor] (w) -- (x) -- (y) -- cycle;
    \draw[line width=0,fill=\negcolor] (y) -- (R) -- (z) -- cycle;
    \draw[line width=0,fill=\negcolor] (z) -- (R) -- (S) -- cycle;
    \draw[line width=0,fill=\negcolor] (S) -- (w) -- (z) -- cycle;
    \draw[line width=0,fill=\negcolor] (S) -- (P) -- (w) -- cycle;
\fi
    \end{tikzpicture}
     &  
    \begin{tikzpicture}[scale=.2]
    \coordinate (P) at (0,0);
\coordinate (Q) at (6,0);
\coordinate (R) at (6,6);
\coordinate (S) at (0,6);

\draw[fill=black] (P) circle (1pt);
\draw[fill=black] (Q) circle (1pt);
\draw[fill=black] (R) circle (1pt);
\draw[fill=black] (S) circle (1pt);
\draw (P) -- (Q) -- (R) -- (S) -- cycle;

\coordinate (w) at (2,2.25);
\coordinate (x) at (4.25,1.75);
\coordinate (y) at (4.5,4.25);
\coordinate (z) at (1.75,4.25);
\draw[fill=black] (w) circle (1pt);
\draw[fill=black] (x) circle (1pt);
\draw[fill=black] (y) circle (1pt);
\draw[fill=black] (z) circle (1pt);

\draw (w) -- (x) -- (y) -- (z) -- cycle;
\draw (w) -- (y);

\draw (P) -- (z);
\draw (P) -- (w);
\draw (P) -- (x);
\draw (Q) -- (x);
\draw (Q) -- (y);
\draw (R) -- (y);
\draw (S) -- (y);
\draw (S) -- (z);

\ifcorners
    \draw (P) node[anchor=north east]{$\PP$};
    \draw (Q) node[anchor=north west]{$\QQ$};
    \draw (R) node[anchor=south west]{$\RR$};
    \draw (S) node[anchor=south east]{$\SS$};
\fi

\ifcolors
    \draw[line width=0,fill=\poscolor] (P) -- (Q) -- (x) -- cycle;
    \draw[line width=0,fill=\poscolor] (Q) -- (R) -- (y) -- cycle;
    \draw[line width=0,fill=\poscolor] (w) -- (x) -- (y) -- cycle;
    \draw[line width=0,fill=\poscolor] (y) -- (S) -- (z) -- cycle;
    \draw[line width=0,fill=\poscolor] (z) -- (P) -- (w) -- cycle;
    \draw[line width=0,fill=\negcolor] (P) -- (x) -- (w) -- cycle;
    \draw[line width=0,fill=\negcolor] (x) -- (Q) -- (y) -- cycle;
    \draw[line width=0,fill=\negcolor] (y) -- (R) -- (S) -- cycle;
    \draw[line width=0,fill=\negcolor] (w) -- (y) -- (z) -- cycle;
    \draw[line width=0,fill=\negcolor] (S) -- (P) -- (z) -- cycle;
\fi
    \end{tikzpicture}
     &  
    \begin{tikzpicture}[scale=.2]
    \coordinate (P) at (0,0);
\coordinate (Q) at (6,0);
\coordinate (R) at (6,6);
\coordinate (S) at (0,6);

\draw[fill=black] (P) circle (1pt);
\draw[fill=black] (Q) circle (1pt);
\draw[fill=black] (R) circle (1pt);
\draw[fill=black] (S) circle (1pt);
\draw (P) -- (Q) -- (R) -- (S) -- cycle;

\coordinate (w) at (1.75,2);
\coordinate (x) at (4.25,1.5);
\coordinate (y) at (4.25,4);
\coordinate (z) at (1.75,4.5);
\draw[fill=black] (w) circle (1pt);
\draw[fill=black] (x) circle (1pt);
\draw[fill=black] (y) circle (1pt);
\draw[fill=black] (z) circle (1pt);

\draw (w) -- (x) -- (y) -- (z) -- cycle;
\draw (w) -- (y);

\draw (P) -- (w);
\draw (P) -- (x);
\draw (Q) -- (x);
\draw (Q) -- (y);
\draw (R) -- (y);
\draw (R) -- (z);
\draw (S) -- (z);
\draw (P) -- (z);

\ifcorners
    \draw (P) node[anchor=north east]{$\PP$};
    \draw (Q) node[anchor=north west]{$\QQ$};
    \draw (R) node[anchor=south west]{$\RR$};
    \draw (S) node[anchor=south east]{$\SS$};
\fi

\ifcolors
    \draw[line width=0,fill=\poscolor] (P) -- (Q) -- (x) -- cycle;
    \draw[line width=0,fill=\poscolor] (P) -- (w) -- (z) -- cycle;
    \draw[line width=0,fill=\poscolor] (w) -- (x) -- (y) -- cycle;
    \draw[line width=0,fill=\poscolor] (z) -- (R) -- (S) -- cycle;
    \draw[line width=0,fill=\negcolor] (P) -- (x) -- (w) -- cycle;
    \draw[line width=0,fill=\negcolor] (Q) -- (y) -- (x) -- cycle;
    \draw[line width=0,fill=\negcolor] (Q) -- (R) -- (y) -- cycle;
    \draw[line width=0,fill=\negcolor] (R) -- (z) -- (y) -- cycle;
    \draw[line width=0,fill=\negcolor] (w) -- (y) -- (z) -- cycle;
    \draw[line width=0,fill=\negcolor] (S) -- (P) -- (z) -- cycle;
\fi
    \end{tikzpicture}
    \\
            \bottomrule
            \addlinespace[5pt]
            \end{tabular}
        \caption{Triangulations with up to 4 interior vertices, canonically 2-colored.}
        \label{table:2color}
        \end{table}
        
        \colorsfalse    
        
    \section{Finiteness and the encyclopedia}\label{sec:encyc}
    
        Continuing in the same vein as Theorems \ref{thm:encyc1} and \ref{thm:encyc2}, we can ask about those points on the variety $V(\T)$ with at most $l$ nonzero coordinates, for any fixed $l$.  This is equivalent to studying the polynomial obtained from $p_{T}$ by setting all but $l$ variables equal to $0$.   It turns out that for fixed $l$, there is a finite list $\encyc_l$ of polynomials in $l$ variables that contains every irreducible factor of every polynomial obtained in this way, up to scalar multiples and renaming of the variables.  This is the Finiteness Theorem of this section.
        
        \begin{defn}[Specialization]
            Let $p$ be a polynomial in $n$ variables.  Let $L$ be any subset of the variables, and let $l=|L|$.  Let $p|_L$ be the polynomial in $l$ variables obtained from $p$ by substituting $0$ for all variables not in $L$.  We call $p|_L$ an {\em $l$-variable specialization} of $p$. 
        \end{defn}
        
        \begin{defn}[Area encyclopedia]
            Let $l$ be a positive integer. Consider two polynomials in $l$ variables to be equivalent if they are equal up to multiplication by a nonzero scalar and permutation of the variables.  The {\em area encyclopedia} $\encyc_l$ is defined as the set of equivalence classes of irreducible factors of the $l$-variable specializations of the area polynomials $p_T$, where $T$ can be any triangulation of a square.
        \end{defn}

        \begin{thm}[Finiteness, see Main Theorem \ref{mainthm:encyc}]\label{thm:encyc}
            For each positive integer $l$, the area encyclopedia $\encyc_l$ is finite.
            
            In particular, for any $l>1$, every polynomial in $\encyc_l$ is equivalent to an irreducible factor of an $l$-variable specialization of $p_T$, where $T$ has at most $3l-2$ triangles.   
        \end{thm}
        
        \begin{example}
             It is easy to see that  
             $$\encyc_1=\{A\},$$
             since any homogeneous polynomial in one variable is equivalent to $A^d$ for some $d$.
             Theorem \ref{thm:encyc2} tells us that $$\encyc_2=\{A+B,\ A-B\}.$$
            For $l=3$, we have calculated that 
             \begin{equation*}\label{eq:encyc3}
                   \encyc_3=\{A+B+C,\ A+B-C,\ A^2+2AB+2AC+B^2-2BC+C^2\}.
             \end{equation*}
             To verify this, we must all examine triangulations $T$ with $3l-2=7$ or fewer triangles (and no subdivisions). We see from Table \ref{table:T.up.to.4.int} that there are only three of these. For each such $T$, we compute $p_T$ using Gr\"obner bases (or by hand). We then find all $3$-variable specializations of these $p_T$ and factor these specializations into irreducibles. The factors we obtain, up to equivalence, are the three polynomials listed in  $\encyc_3$.  It then follows from Theorem \ref{thm:encyc} that for any triangulation $T$, every $3$-variable specialization of $p_T$ is equivalent to one of these three polynomials. 
             
             Similarly, to compute $\encyc_4$, one must examine all triangulations with up to $10$ triangles, namely all $11$ of the triangulations appearing in Table \ref{table:T.up.to.4.int}.  The resulting $\encyc_4$ consists of the $8$ polynomials---3 linears, 4 quadratics, and a quartic---shown in Encyclopedia 1, Vol. 4,  at the end of this paper. 
        \end{example}
        
        \begin{remark}
             Since the number of triangles in a triangulation is always even, the bound $3l-2$ in the Finiteness Theorem can be reduced to the even integer $2\lfloor\frac {3l-2}{2}\rfloor$. With this modification, the bound is then sharp for all $l\leq 4$. For example, when $l=3$, the unique quadratic $q\in \encyc_3$ requires a triangulation $T$ with at least $6$ triangles.  Similarly, the quartic on $\encyc_4$ requires $T$ with at least $10$ triangles. 
        \end{remark}

        The proof of the Finiteness Theorem relies on the following proposition that is a consequence of the Bubble Theorem.
        
        \begin{prop}\label{prop:density}
            Let $T$ be a triangulation with $n$ triangles.  Suppose $L$ is an $l$-element subset of the triangles of $T$, and let $\P(L)\simeq \P^{l-1}$ be the corresponding coordinate subspace of $Y(T)\simeq \P^{n-1}$.  Let $W$ be a component of the intersection $V(T)\cap \P(L)$, and assume that $W$ is not equal to the hyperplane $H$ in $\P(L)$ defined by the vanishing of the sum $\sigma$ of all $l$ coordinates.  Then 
            \begin{itemize}
                \item A generic point of $W$ has the property that no nonempty subset of its coordinates sums to $0$.
                \item $W\cap \im\area$ is dense in $W$.
            \end{itemize}
        \end{prop}
        
        \begin{proof}
            Let $q$ be the polynomial in $l$ variables that defines $W$.  Note that $q$ is a factor of the specialization $p|_L$.  From  Theorem \ref{thm:encyc1}, we see that that all of the leading terms of $p$ are nonzero, and hence all of the leading terms of $q$ are also nonzero.  Hence $q$ cannot equal the linear polynomial $ \sum_{s\in S} A_s$ for any proper subset $S$ of $L$.  From this and our assumption that $W\neq H$, we see that a generic point $w\in W$ has the property that no subset of the coordinates of $w$ sums to $0$. Using the Bubble theorem, we conclude that $w$ is in $\im\area$, verifying the second assertion in the Proposition. 
        \end{proof}

        \begin{proof}[Proof of Finiteness]
            For $l=1$, we easily see that $\encyc_1=\{A\}$ as noted above.
        
            Let $l>1$.  We first show explicitly that $\sigma$ (the sum of the variables) is in $\encyc_l$.  If $l=2$, then we may take $T$ to be the triangulation with $4$ triangles with $L$ a pair of opposite triangles.  In this case $p|_L=\sigma$. Repeatedly subdividing the triangles of $L$, it is easy to find for each $l>2$, a triangulation $T$ with at most $l+3$  triangles such that $p_T$ is linear and has $\sigma$ among its specializations.  
        
            Now let $q\in\encyc_l$, and assume that $q\neq \sigma$. Since $q\in \encyc_l$, there exists a triangulation $T$ and an $l$-element subset $L$ of the triangles of $T$ such that $q$ is an irreducible factor of $p|_L$. Our goal is to show that there is a triangulation $T'$ with at most $3l-2$ triangles such that some irreducible factor of an $l$-variable specialization of $p_{T'}$ is equivalent to $q$.  Let $n$ be the number of triangles of $T$.
     
            As in Proposition \ref{prop:density}, we let $\P(L)\simeq \P^{l-1}$ be the corresponding coordinate subspace of $Y\simeq \P^{n-1}$.  Let $W$ be the component $V(T)\cap \P(L)$ of defined by $q$.    Since $q\neq \sigma$, Proposition \ref{prop:density} implies that $W\cap\im\area$ is dense in $W$. 
                        
            We now claim that  
            $$W\subset \bigcup_{T',L',h} \iota_h(V(T')\cap \P(L')),$$
            where the union is taken over all triangulations $T'$ with at most $3l-2$ triangles, $L'$ is an $l$-element subset of $T'$ and $h:L'\to L$ is a bijection, giving rise to an identification $\iota_h$ of $\P(L')$ with $\P(L)$. 
            
            To see this, take a generic point $w\in W$. Using Proposition \ref{prop:density}, we may assume that that no nonempty subset of the coordinates of $w$ sums to zero and $w\in\im\area$.  Thus $w=\area(\rho)$  for some drawing $\rho$ of $T$ with exactly $l$ nondegenerate triangles. By Theorem \ref{thm:simplification}, there is a triangulation $T'$ with at most $3l-2$ triangles and a point $w'$ of $V(T')$ whose $l$ nonzero coordinates agree with those of $w$.  Hence $w$ is in the union above.  Therefore a dense subset of $W$ is contained in the union above.  Thus all of $W$ is contained in the union, as desired.
            
            Now since  $W$ is irreducible, there must be a single set  $\iota_h(V(T')\cap \P(L'))$ that contains $W$. Therefore the irreducible $q$ defining $W$ divides $(p_{T'})|_{L'}$ with the substitution of variables determined by $h$.  This completes the proof.
        \end{proof}

        \begin{cor}\label{cor:areassatisfyencyc}
            Let $T$ be a triangulation and let $\rho$ be a drawing of $T$ with $l$ nondegenerate triangles. If the areas of these $l$ triangles are $a_1, \dots, a_l$, then $(a_1, \dots, a_l)$ satisfies some polynomial of $\encyc_l$.
        \end{cor}

        \begin{proof}
           Since $\rho$ is a drawing of $T$, some permutation of $(a_1, \dots, a_l, 0, \dots, 0)$ satisfies $p_T$.  Hence if $p_T|_L$ is the specialization of $p_T$ to the $l$ nondegenerate triangles, the point $(a_1, \dots, a_l)$ satisfies $p_T|_L$, and hence satisfies some polynomial in $\encyc_l$.
        \end{proof}

        \begin{cor} [See Main Theorem \ref{mainthm:nonconverse}] \label{cor:EimpliesT}
            If an $l$-tuple $(a_1,\dots,a_l)$ of complex numbers satisfies some polynomial of $\encyc_l$, and no subset sums to $0$, then there exists a triangulation $T$ with at most $3l-2$ triangles and a drawing $\rho$ of $T$ with exactly $l$ nondegenerate triangles such that $\area(\rho)=(a_1,\dots,a_l)$.
        \end{cor}
        
        \begin{proof}
            For some $T$, there is a point $w\in V(T)$ whose nonzero coordinates are exactly $a_1,\dots,a_l$. Applying Theorem \ref{thm:characterization}, we see that $w$ is in $\im \area_{T}$, so there's a drawing $\rho$ of $T$ with $\area(\rho)=w$.  By Theorem \ref{thm:simplification} we may assume $T$ satisfies the stated bound on the number of triangles.
        \end{proof} 
        
        Combining the previous corollaries, we obtain a nearly-complete characterization of those tuples that are realized as the areas of a drawing of a triangulation.
        
        \begin{cor}\label{cor:realizingareas}
            Suppose $(a_1,\dots,a_l)$ is an $l$-tuple of complex numbers such that no subset sums to $0$.  Then there is drawing of some triangulation $T$ realizing the areas $(a_1,\dots,a_l)$ if and only if $(a_1,\dots,a_l)$ satisfies some polynomial of $\encyc_l$.
        \end{cor}
        
        Regarding classical dissections, Corollary \ref{cor:areassatisfyencyc} also implies the following.
        
        \begin{cor} [See Main Theorem \ref{mainthm:dissections}] \label{cor:dissection}
            Given a dissection of a square into $l$ triangles with areas $(a_1,\dots,a_l)$, the tuple $(a_1,\dots,a_l)$ satisfies some polynomial of $\encyc_l$.
        \end{cor}

        \begin{proof}
            By results of \cite{chapter3}, any dissection (``simplicial'' or not) can be viewed as a drawing of a triangulation, and so the previous corollary applies.
        \end{proof}
                
        \begin{question}
            Is the converse of the previous corollary true?  That is if $(a_1,\dots,a_l)$ is an $l$-tuple of positive real numbers satisfying some polynomial of $\encyc_l$, does there necessarily exist a dissection of a square into $l$ triangles with those areas? 
        \end{question}
        
        Regarding the latter question, we know from Corollary \ref{cor:dissection} that there is a drawing of a triangulation realizing the positive areas $(a_1,\dots,a_l)$. However, such a drawing is only a dissection if the vertices are drawn in $\R^2$.

        \begin{example}[A complex equidissection]\label{ex:complexequi}
            For the triangulation shown in Figure \ref{fig:complexequi}, label the unmarked triangles $A_1$ through $A_4$.  There is no drawing of this $T$ in the (real) plane in which these four triangles have positive areas (let alone equal areas) and the six marked triangles are degenerate, so it surprised us to notice that the point $[1:1:1:1:0:0:0:0:0:0]$ satisfies the polynomial $p_T$ (as does more generally $[1:x:1:x:0:0:0:0:0:0]$).  By Corollary \ref{cor:realizingareas} there must be a drawing realizing these areas, and indeed there is:  one maps the interior vertices to the points $(\frac{1+i}{2},0), (1,\frac{1+i}{2}), (\frac{1-i}{2},1)$, and $(0,\frac{1-i}{2})$. 
            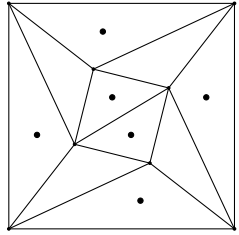
\begin{figure}
                \centering
                 
    \begin{tikzpicture}[scale=.5]
    \coordinate (P) at (0,0);
\coordinate (Q) at (6,0);
\coordinate (R) at (6,6);
\coordinate (S) at (0,6);

\draw[fill=black] (P) circle (1pt);
\draw[fill=black] (Q) circle (1pt);
\draw[fill=black] (R) circle (1pt);
\draw[fill=black] (S) circle (1pt);
\draw (P) -- (Q) -- (R) -- (S) -- cycle;

\coordinate (w) at (1.75,2.25);
\coordinate (x) at (3.75,1.75);
\coordinate (y) at (4.25,3.75);
\coordinate (z) at (2.25,4.25);
\draw[fill=black] (w) circle (1pt);
\draw[fill=black] (x) circle (1pt);
\draw[fill=black] (y) circle (1pt);
\draw[fill=black] (z) circle (1pt);

\draw (w) -- (x) -- (y) -- (z) -- cycle;
\draw (w) -- (y);

\draw (P) -- (w);
\draw (P) -- (x);
\draw (Q) -- (x);
\draw (Q) -- (y);
\draw (R) -- (y);
\draw (R) -- (z);
\draw (S) -- (z);
\draw (S) -- (w);

\ifcorners
    \draw (P) node[anchor=north east]{$\PP$};
    \draw (Q) node[anchor=north west]{$\QQ$};
    \draw (R) node[anchor=south west]{$\RR$};
    \draw (S) node[anchor=south east]{$\SS$};
\fi

\ifcolors
    \draw[line width=0,fill=\poscolor] (P) -- (Q) -- (x) -- cycle;
    \draw[line width=0,fill=\poscolor] (x) -- (Q) -- (y) -- cycle;
    \draw[line width=0,fill=\poscolor] (y) -- (R) -- (z) -- cycle;
    \draw[line width=0,fill=\poscolor] (P) -- (x) -- (w) -- cycle;
    \draw[line width=0,fill=\poscolor] (w) -- (x) -- (y) -- cycle;
    \draw[line width=0,fill=\poscolor] (w) -- (y) -- (z) -- cycle;
    \draw[line width=0,fill=\poscolor] (w) -- (z) -- (S) -- cycle;
    \draw[line width=0,fill=\poscolor] (z) -- (R) -- (S) -- cycle;
    \draw[line width=0,fill=\negcolor] (Q) -- (R) -- (y) -- cycle;
    \draw[line width=0,fill=\negcolor] (S) -- (P) -- (w) -- cycle;
\fi

\draw[fill=black] (2.5,5.25) circle (2pt);
\draw[fill=black] (3.5,.75) circle (2pt);
\draw[fill=black] (5.25,3.5) circle (2pt);
\draw[fill=black] (.75,2.5) circle (2pt);
\draw[fill=black] (3.25,2.5) circle (2pt);
\draw[fill=black] (2.75,3.5) circle (2pt);

\ifnewareas
    \draw (1.75,1.5) node{\scriptsize $A$};
    \draw (4.5,2) node{\scriptsize $B$};
    \draw (4,4.5) node{\scriptsize $C$};
    \draw (1.5,4) node{\scriptsize $D$};
\fi
    \end{tikzpicture}
    
                \caption{This triangulation can be drawn with the six marked triangles degenerate and the four marked triangles all having area $1$, but not while giving the interior vertices real coordinates.}
                \label{fig:complexequi}
            \end{figure}
        \end{example}

    \section{Illustrating the encyclopedia}\label{sec:illustrating}
    
        In this section, we discuss how to associate pictures to the polynomials in the encyclopedia. In particular, we show that every polynomial $p\in \encyc$ can be illustrated by a constrained triangulation, an idea introduced in \cite{chapter3}.  For $l=2, 3, 4$, we give an illustrated encyclopedia that shows the pictures corresponding to each polynomial.

        Theorem \ref{thm:encyc} shows that each polynomial on $\encyc_l$ is associated to a variety of drawings with $l$ nondegenerate triangles.  In the language of \cite{chapter3}, we will associate constrained triangulations $\T=(T,\CC)$ to each polynomial of $\encyc_l$.
    
        \subsection{Constrained triangulations}
            The ideas in this subsection are treated more thoroughly in \cite{chapter3}.  Here we briefly give the relevant definitions.
            
            \begin{defn} [Constrained triangulation, living triangle, drawing]\label{def:constrained} 
                A \emph{constrained triangulation} $\T$ is a pair $\T=(T,\CC)$, where $T$ is a triangulation and where $\CC=\{C_i\}$ is a  set of \emph{(collinearity) constraints}.  
                Each collinearity constraint $C_i$ is a set of vertices of $T$ of the form $\Vxs(S_i)$ where $S_i$ is a contiguous set of triangles of $T$.  
                (This means that there is a connected subgraph of the dual graph to $T$ whose vertices are the triangles of $S_i$.)  
                We require the sets $S_i$ of triangles to be disjoint, although the constraints $C_i$ need not be.  
                
                A 2-cell of $T$ is called \emph{alive} or \emph{living} if there is no constraint containing all of its vertices.  
                                
                A \emph{drawing} of the constrained triangulation $\T=(T,\CC)$ is a drawing $\rho$ of the (unconstrained) triangulation $T$ (cf.~Definition \ref{def:honestdrawings}) such that for each $C\in\CC$ there is a line $\ell_C\subset\C^2$ with $\rho(v)\in\ell_C$ for each $v\in C$.
            \end{defn}

            The living triangles of $\T$ are the only triangles whose image in a drawing of $\T$ can be nondegenerate.  
            
            We next introduce a combinatorial notion of ``factoring'' for constrained triangulations.  See also \cite{chapter3}.
            
            \begin{defn} [Combinatorially irreducible]\label{def:combirr}
                A constrained triangulation $\T=(T,\CC)$ is \emph{combinatorially irreducible} if $|C\cap C'|\leq 1$ for any two distinct constraints $C,C'\in\CC$.
            \end{defn}
        
            \begin{defn}[Amalgamation]
                Let $\T=(T,\CC)$ be a constrained triangulation, suppose $uv$ is an edge of $T$, and suppose $u,v\in C \cap C'$ (where $C,C'\in\CC$ are distinct).  The \emph{amalgamation of $\T$} (or of $C$ and $C'$) is the constrained triangulation $\T'=(T',\CC')$ with $T'=T$ and $\CC'=\CC$ except that $C, C'$ have been replaced by their union.
                
                The \emph{maximal amalgamation} of $\T$ is the constrained triangulation obtained from $\T$ obtained by applying all possible amalgamations to $\T$.
            \end{defn}

            Observe that if $\bar \T$ is an amalgamation of $\T$ (including the maximal amalgamation) then any drawing of $\bar \T$ is also a drawing of $\T$.  A drawing of $\T$ may or may not be a drawing of $\bar\T$; see Example \ref{ex:combred}.
            
            An amalgamation of $\T$ is a type of \emph{combinatorial factor} of $\T$.  There is another type, called a contraction, obtained by contracting the edge $uv$ in the definition above (if this results in a triangulation); see \cite{chapter3}.  

            \begin{example}[Combinatorial reduction and factors] \label{ex:combred}
                The basic example to keep in mind is the $\T$ shown in Figure \ref{fig:BE-reduction} (left) with two constraints, each indicated by a mark inside a doomed triangle.  The figure also shows two drawings of the underlying triangulation $T$, one simple (middle) and one not (right).  The first drawing (middle) is also a drawing of the maximal amalgamation $\T'$ of $\T$, which has only one constraint (consisting of four vertices).  The second drawing (right) is a simple drawing of the other factor $\T''$, which is an unconstrained triangulation.  Both $\T'$ and $\T''$ are combinatorially irreducible, and any drawing of $\T$ is either a drawing of $\T'$ or a drawing of $\T''$ (or both).
                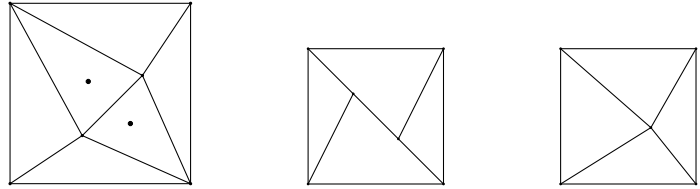
\begin{figure}
                    \centering
                    \begin{tikzpicture}[scale=.4]
                        \draw (0,0) -- (0,6) -- (6,6) -- (6,0) -- cycle;
\draw (0,0) -- (2.4,1.6) -- (4.4,3.6) -- (6,6);
\draw (6,0) -- (2.4,1.6) -- (0,6) -- (4.4,3.6) -- cycle;
\draw[fill=black] (0,0) circle (1pt);
\draw[fill=black] (0,6) circle (1pt);
\draw[fill=black] (6,0) circle (1pt);
\draw[fill=black] (6,6) circle (1pt);
\draw[fill=black] (2.4,1.6) circle (1pt);
\draw[fill=black] (4.4,3.6) circle (1pt);
% \draw (3,.5) node{$A$};
 \draw [fill=black] (4,2) circle (2pt);
% \draw (5.2,2.8) node{$C$};
% \draw (1,2.2) node{$D$};
 \draw [fill=black] (2.6,3.4) circle (2pt);
% \draw[dotted] (4,2) -- (2.6,3.4);
% \draw (2.6,3.4) node{$E$};
% \draw (3.8,5) node{$F$};
% \draw (0,0) node[anchor=north east]{$\PP$};
% \draw (6,0) node[anchor=north west]{$\QQ$};
% \draw (6,6) node[anchor=south west]{$\RR$};
% \draw (0,6) node[anchor=south east]{$\SS$};
% \draw (2.4,1.6) node[anchor=north]{$u$};
% \draw (4.4,3.6) node[anchor=west]{$v$};
                    \end{tikzpicture}
                    \qquad \qquad
                    \begin{tikzpicture}[scale=.3]
                        \draw (0,0) -- (0,6) -- (6,6) -- (6,0) -- cycle;
\draw (6,0) -- (0,6) ;
\draw (0,0) -- (2,4);
\draw (4,2) -- (6,6);
\draw [fill] (0,0) circle (1pt);
\draw [fill] (0,6) circle (1pt);
\draw [fill] (6,6) circle (1pt);
\draw [fill] (6,0) circle (1pt);
\draw [fill] (4,2) circle (1pt);
\draw [fill] (2,4) circle (1pt);
                    \end{tikzpicture}
                    \qquad \qquad
                    \begin{tikzpicture}[scale=.3]
                        \draw (0,0) -- (0,6) -- (6,6) -- (6,0) -- cycle;
\draw (4,2.5) -- (0,0);
\draw (4,2.5) -- (6,0);
\draw (4,2.5) -- (0,6);
\draw (4,2.5) -- (6,6);
\draw [fill] (0,0) circle (1pt);
\draw [fill] (0,6) circle (1pt);
\draw [fill] (6,6) circle (1pt);
\draw [fill] (6,0) circle (1pt);
\draw [fill] (4,2.5) circle (1pt);
                    \end{tikzpicture}
                    \caption{A combinatorially reducible $\T$ and its (two) irreducible factors, $\T'$ and $\T''$.}
                    \label{fig:BE-reduction}
                \end{figure}
            \end{example}

            Building on the previous example, the following lemma says that every simple drawing of a constrained triangulation $\T$ either has a subset of living triangles whose areas sum to zero or else is a drawing of a particular combinatorially irreducible factor of $\T$, namely the maximal amalgamation.  (A drawing $\rho$ of $\T=(T,\CC)$ is \emph{simple} if it is simple when viewed as a drawing of $T$; see Definition \ref{def:simple}.)

            \begin{lem}\label{lem:reducing}
                Let $\T=(T,\CC)$ be a constrained triangulation.  Suppose there is a simple drawing $\rho$ with the property that no nonempty subset of the living triangles of $\T$ is drawn with areas summing to zero.  Then (a) the maximal amalgamation $\bar \T$ of $\T$ is combinatorially irreducible, and (b) $\rho$ is a simple drawing of $\bar \T$.  
            \end{lem}
            
            \begin{proof}
                If $\T$ is combinatorially irreducible then $\bar \T=\T$ and we are done, so assume on the contrary that $\T$ is combinatorially reducible.  We first prove (a) by showing that there are constraints $C,C'\in\CC$ that share an edge, hence can be amalgamated.  It follows that (under the hypothesis of the lemma) the maximal amalgamation $\bar\T$ is combinatorially irreducible.
                
                As $\T$ is combinatorially reducible there must be distinct constraints $C,C'\in\CC$ and distinct vertices $x,y\in C\cap C'$.  If $xy$ is an edge of $T$ then we are done, so suppose not. There is a path $\gamma$ from $x$ to $y$ in the 1-skeleton of $T$ that uses only vertices of $C$; by rechoosing $y$ if necessary we may assume this path doesn't use any vertices of $C\cap C'$ other than $x$ and $y$. There is also a path $\gamma'$ from $x$ to $y$ using only vertices of $C'$.  The union of these paths is a cycle $\Gamma$ that encloses some of the triangles of $T$.  As $\rho$ maps each of $\gamma,\gamma'$ to a set of collinear points, the total area enclosed by the cycle is zero.  (To see this, note that if $\rho(x)\ne\rho(y)$ then the lines $\ell_C$ and $\ell_{C'}$ coincide, and the assertion follows.  Otherwise $\rho(x)=\rho(y)$, and in this case the lines $\ell_C$ and $\ell_{C'}$ may be distinct but each of $\gamma,\gamma'$ individually maps to a loop enclosing zero area.)  Our hypothesis therefore implies that none of the enclosed triangles is living in $\T$.  If there were no available amalgamation then all vertices of triangles enclosed by $\Gamma$ would be part of a single constraint, along with all vertices of $C$ and of $C'$.  This contradicts our assumption that $C$ and $C'$ are distinct.

                We now prove (b).
                Suppose $xy$ is an edge and $x,y\in C\cap C'$.  We claim that $\rho(x)\ne\rho(y)$.  Suppose $\rho(x)=\rho(y)$, i.e., the edge $xy$ is elastic.  Note that $xy$ is an interior edge, since it is contained in two constraints.  By simplicity, this edge cannot be contracted, so it must be part of a subdivision.  But then the sum of the areas inside the subdivision is zero, hence by hypothesis all of these triangles are degenerate, which again contradicts simplicity of $\rho$.  We conclude that $\rho(x)\ne\rho(y)$, as claimed.
                
                Let $\ell$ be the line determined by the points $\rho(x)$ and $\rho(y)$.  As $\rho$ is a drawing of $\T$, and $C$ and $C'$ are constraints of $\T$, we must have that all points of $C\cup C'$ are drawn on $\ell$.  Therefore $\rho$ is a drawing of the amalgamation of $C$ and $C'$.
                
                Repeating this process until there are no remaining amalgamations shows that $\rho$ is a (simple) drawing of $\bar\T$.
            \end{proof}

        \subsection{Maximal amalgamations illustrate the encyclopedia} 
            We now show that if $q\in \encyc_l$ with $l>1$, then there is maximally amalgamated combinatorially irreducible $\T$ that {\it illustrates} $q$ in the sense that almost every point in the zero set of $q$ is realized by a drawing of $\T$.  
            
            \begin{defn}[Illustration]
                Let $q\in \encyc_l$.  Given a constrained triangulation $\T$ with $l$ living triangles, we say that $\T$ {\em illustrates} $q$ if a generic point of the zero set of $q$ is equal to $\area(\rho)$ for a simple drawing $\rho$ of $\T$.
             \end{defn}
            
            Note that if $\T=(T,\CC)$ illustrates $q\in\encyc_l$ then by the simplicity of $\rho$ the triangulation $T$ has at most $3l-2$ triangles.

            \begin{thm}[Illustration, see Main Theorem \ref{mainthm:illustrate}]\label{thm:illustration}  
                Let $l\ge 2$ and let $q\in\encyc_l$.  Then there is a maximally amalgamated $\T=(T,\CC)$ that is combinatorially irreducible and that illustrates $q$.
            \end{thm}

            \begin{proof}
                We suppose that $q\ne \sigma$; the case $q=\sigma$ can easily be treated separately.
                
                Choose a triangulation $T$ and a subset $L$ of the triangles of $T$ such that $q$ is an irreducible factor of the specialization $p_T|_L$.  
                Let $W$ be the irreducible component of $V(T)\cap \P(L)$ corresponding to $q$, and let $w$ be a generic point of $W$. By Proposition \ref{prop:density}, we may assume that there is a drawing $\rho$ with $w=\area(\rho)$ and that no subset of the areas of nondegenerate triangles of $\rho$ sums to zero.  
                
                By Theorem \ref{thm:simplification} there is also a triangulation $T'$ and a simple drawing $\rho'$ with $\area(\rho')=w$.  Let $L'$ be the set of nondegenerate triangles of $\rho'$.  Let $\T'=(T',\CC)$ be the constrained triangulation with $\CC=\{\Vxs(\Delta)\ |\ \Delta\not\in L'\}$ consisting of one constraint for each triangle $\Delta$ of $T'$ that is degenerate in the drawing $\rho'$.
                Lemma \ref{lem:reducing} shows that $\rho'$ is a drawing of the maximal amalgamation $\bar\T'$ of $\T'$, and that $\bar\T'$ is combinatorially irreducible.  
                
                As $\area(\rho')=w$, we conclude that for any generic point $w$ of $W$, there exists a triangulation $T'$ with at most $3l-2$ triangles and a maximally amalgamated, combinatorially irreducible $\T'=(T',\CC)$ with exactly $l$ living triangles such that $w$ is the area of a simple drawing of $\T'$.  As $W$ is irreducible, it follows that there is a single such $\T'$ that illustrates $q$. 
            \end{proof}

        \subsection{Epilogue:  drawability}\label{sec:epilogue}
            Our interest in combinatorial irreducibility stems from the fact that, in our experience, the relevant space of drawings of a combinatorially irreducible constrained triangulation is itself typically an irreducible variety. However Theorem \ref{thm:illustration} leaves open the possibility that the constrained triangulation $\T$, despite being combinatorially irreducible, could illustrate multiple polynomials.
            
            We know of only one scenario in which this phenomenon occurs, but it is well-controlled. Namely, suppose that in the constrained triangulation $\T$, either there are two distinct constraints, one containing $\PP$ and $\QQ$ and the other containing $\RR$ and $\SS$, or else there are two distinct constraints, one containing $\PP$ and $\SS$ and the other containing $\QQ$ and $\RR$.  In this situation $\T$ is called \emph{toroidally reducible}, and the drawing space splits into one component $X_0$ whose drawings $\rho$ satisfy $\rho(\PP)=\rho(\QQ)$ and $\rho(\RR)=\rho(\SS)$ (or the other possibility), and the rest which will typically consist of just one further component, if $\T$ is combinatorially irreducible.  The area of any $\rho\in X_0$ is in $H$, and this component satisfies the polynomial $\sigma$, so while $\T$ does illustrate multiple polynomials, one of them is $\sigma$.  This  phenomenon arises in Volume 4 of the encyclopedia.  See \cite{chapter3} for further discussion of toroidal reducibility.
            
            Other than toroidal reducibility, however, it is possible that each $\T$ arising from Theorem \ref{thm:illustration} illustrates a unique polynomial.  In previous work \cite{chapter3}, we have shown that if $\T$ is combinatorially irreducible and has any drawing that is \emph{generic} in a certain sense, then the relevant space of drawings can be paramaterized as a rational variety and is therefore irreducible.  We call such a $\T$ \emph{drawable}.  See \cite{chapter3} for the definition of generic; for our purpose here it suffices to say that the drawings $\rho\in X_0$ of a toroidally reducible $\T$ are not generic. 
            
            So, starting with $\sigma \ne q\in \encyc$, if $\T$ is a drawable constrained triangulation illustrating $q$, then in fact every generic drawing of $\T$ has areas that satisfy $q$.  In this sense, an arbitrary generic drawing of $\T$ provides a perfect illustration of the polynomial $q$.  In the language of \cite{chapter3}, it also follows that such a $\T$ is \emph{hyper}. 
              
            Every $\T$ arising from Theorem \ref{thm:illustration} with $l\le 4$ is drawable, and the drawings we have shown are generic (except the toroidally reducible illustrations of $\sigma$).  In general the existence of generic drawings of constrained triangulations seems to be a matter of incidence geometry that we do not know how to resolve.  See Questions $1$, $1'$, $2$, $2'$ in \cite{chapter3} for more on this.

\part{The encyclopedia}

    \section{Volumes 1 through 4}\label{sec:theencyc}
    
            We present here the first four volumes of the illustrated area encylopedia $\encyc$.  For each polynomial in $\encyc_2\cup\encyc_3\cup\encyc_4$ we have found all constrained triangulations that illustrate $q$; by Theorem \ref{thm:illustration} each $q$ has such a $\T=(T,\CC)$, and the underlying triangulation $T$ appears on Table \ref{table:T.up.to.4.int}.  (Note however that not all $T$'s from Table \ref{table:T.up.to.4.int} appear.)
            
            Some explanation of the encyclopedia contents is appropriate.  We use a dot to indicate triangles whose vertices are part of a constraint.  Constrained triangulations are assumed to be maximally amalgamated.
            
            Each volume of $\encyc$ contains the polynomial $\sigma$, the sum of the variables.  In any drawing $\rho$ illustrating $\sigma$, the boundary $\PP\QQ\RR\SS$ maps to a degenerate parallelogram.  We have drawn this boundary as a segment, with $\rho(\PP)=\rho(\QQ)$ and $\rho(\RR)=\rho(\SS)$. 
            
            In some instances there are multiple $\T$'s illustrating a polynomial $q$.  We include multiple illustrations only when the drawings are essentially different.

            Whenever possible the constrained triangulation $\T$ is chosen without subdivisions.  We also omit drawings that contain subdivisions, unless this is the only way to illustrate a particular $q$.
            Defining ``subdivision'' in the context of drawings requires some care, because a drawing $\rho\in X(T)$ can contain a subdivision even when $T$ doesn't.  Roughly, $\rho\in X(T)$ \emph{contains a subdivision} if the nondegenerate triangles of $\rho$ can be interpreted as (the image of) a drawing of a possibly different triangulation $T'$ that does have a subdivision, in which some of the nondegenerate triangles of $\rho$ are contained in the subdivision of $T'$.  This situation is easy to understand algebraically; in particular any $q$ that is illustrated by a drawing with a subdivision is itself subdivided in an algebraic sense, discussed in the next section.  However when there is a more interesting illustration of $q$, we give it.
            
            One can adapt Theorem \ref{thm:encyc2} to constrained triangulations with the same conclusion; in particular a constrained triangulation (if it is hyper) also has a canonical 2-coloring.  The signs of leading terms of a polynomial are preserved by specialization but not by factoring, and the canonical 2-coloring of $\T=(T,\CC)$ is not in general inherited from the canonical 2-coloring of $T$ shown in Table \ref{table:2color}.

            Without further ado, here are Volumes 1--4 of the area encyclopedia.
            
            \bigskip
%           \newpage
            %%%%%%%%%%%%%%%%%%%%%%%%%%%%%%%%%%%%%%%%%
% Dictionary
% LaTeX Template
% Version 1.0 (20/12/14)
%
% This template has been downloaded from:
% http://www.LaTeXTemplates.com
%
% Original author:
% Vel (vel@latextemplates.com) inspired by a template by Marc Lavaud
%
% License:
% CC BY-NC-SA 3.0 (http://creativecommons.org/licenses/by-nc-sa/3.0/)
%
%%%%%%%%%%%%%%%%%%%%%%%%%%%%%%%%%%%%%%%%%

\hrule
\volume{1}
    \degree{Linear}
        \entry{A}{}
            {Also denoted $\sigma$.  This polynomial is in the abridged encyclopedia.}
            {The zero set of this polynomial in $\P^0$ is empty, so there are no points to illustrate. This is the only polynomial in the encyclopedia that is not illustrated by a constrained triangulation. 
            }

\volume{2}
    \degree{Linear}
        \entry{A+B}{}
            {Also denoted $\sigma$.  This polynomial is a subdivision of $A$.}
            {\begin{center}
                \begin{tabular}{ c p{1in} c }
                    Illustrated by & & with drawing \\
                    \cornerstrue  
    \newareastrue
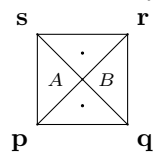
    \begin{tikzpicture}[scale=.2]
    \coordinate (P) at (0,0);
\coordinate (Q) at (6,0);
\coordinate (R) at (6,6);
\coordinate (S) at (0,6);

\draw[fill=black] (P) circle (1pt);
\draw[fill=black] (Q) circle (1pt);
\draw[fill=black] (R) circle (1pt);
\draw[fill=black] (S) circle (1pt);
\draw (P) -- (Q) -- (R) -- (S) -- cycle;

\coordinate (x) at (3,3);
\draw[fill=black] (x) circle (1pt);

\draw (P) -- (x);
\draw (Q) -- (x);
\draw (R) -- (x);
\draw (S) -- (x);

\draw[fill=black] (3,1.25) circle (2pt);
\draw[fill=black] (3,4.75) circle (2pt);

\ifnewareas
    \draw (1.2,3) node{\scriptsize $A$};
    \draw (4.6,3) node{\scriptsize $B$};
\fi

\ifcorners
    \draw (P) node[anchor=north east]{$\PP$};
    \draw (Q) node[anchor=north west]{$\QQ$};
    \draw (R) node[anchor=south west]{$\RR$};
    \draw (S) node[anchor=south east]{$\SS$};
\fi
    \end{tikzpicture}
    \newareasfalse
     &  & \cornerstrue  
    
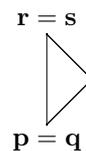
\begin{tikzpicture}[scale=.2]
    \coordinate (P) at (0,0);
\coordinate (Q) at (0,0);
\coordinate (R) at (0,6);
\coordinate (S) at (0,6);

\draw[fill=black] (P) circle (1pt);
\draw[fill=black] (Q) circle (1pt);
\draw[fill=black] (R) circle (1pt);
\draw[fill=black] (S) circle (1pt);
\draw (P) -- (Q) -- (R) -- (S) -- cycle;

\coordinate (x) at (3,3);
\draw[fill=black] (x) circle (1pt);

\draw (P) -- (x);
\draw (Q) -- (x);
\draw (R) -- (x);
\draw (S) -- (x);

\ifcorners
    \draw (P) node[anchor=north]{$\PP=\QQ$};
%    \draw (Q) node[anchor=north west]{$\QQ$};
%    \draw (R) node[anchor=south west]{$\RR$};
    \draw (S) node[anchor=south]{$\RR=\SS$};
\fi

\ifcolors
    \draw[line width=0, fill=\poscolor] (P) -- (Q) -- (x) -- cycle;
    \draw[line width=0, fill=\poscolor] (R) -- (S) -- (x) -- cycle;
    \draw[line width=0, fill=\negcolor] (Q) -- (R) -- (x) -- cycle;
    \draw[line width=0, fill=\negcolor] (S) -- (P) -- (x) -- cycle;
\fi

\ifareas
    \draw (3,1.5) node{$A$};
    \draw (4.5,3) node{$B$};
    \draw (3,4.5) node{$C$};
    \draw (1.3,3) node{$D$};
\fi

\ifnewareas
    \draw (3,1.25) node{\scriptsize $D$};
    \draw (4.75,3) node{\scriptsize $B$};
    \draw (3,4.75) node{\scriptsize $C$};
    \draw (1.25,3) node{\scriptsize $A$};
\fi
    \end{tikzpicture}
    
                    \cornersfalse
                \end{tabular}
            \end{center}
            This drawing should be interpreted as a ``triangulation of a segment.''  There are two non-degenerate triangles, arranged back-to-back, whose areas cancel.
            }

        \entry{A-B}{}
            {This polynomial is in the abridged encyclopedia.}
            {\begin{center}
                \begin{tabular}{ c p{1in} c }
                    Illustrated by & & with drawing \\
                     
    \newareastrue
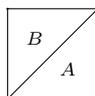
    \begin{tikzpicture}[scale=.2]
    \coordinate (P) at (0,0);
\coordinate (Q) at (6,0);
\coordinate (R) at (6,6);
\coordinate (S) at (0,6);

\draw[fill=black] (P) circle (1pt);
\draw[fill=black] (Q) circle (1pt);
\draw[fill=black] (R) circle (1pt);
\draw[fill=black] (S) circle (1pt);
\draw (P) -- (Q) -- (R) -- (S) -- cycle;

\draw (P) -- (R);

\ifcorners
    \draw (P) node[anchor=north east]{$\PP$};
    \draw (Q) node[anchor=north west]{$\QQ$};
    \draw (R) node[anchor=south west]{$\RR$};
    \draw (S) node[anchor=south east]{$\SS$};
\fi

\ifareas
    \draw (4,2) node{$A$};
    \draw (1.8,4) node{$B$};
\fi

\ifnewareas
    \draw (4,2) node{\scriptsize $A$};
    \draw (1.8,4) node{\scriptsize $B$};
\fi

\ifcolors
    \draw[line width=0, fill=\poscolor] (P) -- (Q) -- (R) -- cycle;
    \draw[line width=0, fill=\negcolor] (R) -- (S) -- (P) -- cycle;
\fi
    \end{tikzpicture}
    \newareasfalse
     & &  
    
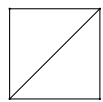
\begin{tikzpicture}[scale=.2]
    \coordinate (P) at (0,0);
\coordinate (Q) at (6,0);
\coordinate (R) at (6,6);
\coordinate (S) at (0,6);

\draw[fill=black] (P) circle (1pt);
\draw[fill=black] (Q) circle (1pt);
\draw[fill=black] (R) circle (1pt);
\draw[fill=black] (S) circle (1pt);
\draw (P) -- (Q) -- (R) -- (S) -- cycle;

\draw (P) -- (R);

\ifcorners
    \draw (P) node[anchor=north east]{$\PP$};
    \draw (Q) node[anchor=north west]{$\QQ$};
    \draw (R) node[anchor=south west]{$\RR$};
    \draw (S) node[anchor=south east]{$\SS$};
\fi

\ifareas
    \draw (4,2) node{$A$};
    \draw (1.8,4) node{$B$};
\fi

\ifnewareas
    \draw (4,2) node{\scriptsize $A$};
    \draw (1.8,4) node{\scriptsize $B$};
\fi

\ifcolors
    \draw[line width=0, fill=\poscolor] (P) -- (Q) -- (R) -- cycle;
    \draw[line width=0, fill=\negcolor] (R) -- (S) -- (P) -- cycle;
\fi
    \end{tikzpicture}
    
                \end{tabular}
            \end{center}
            }

\volume{3}
    \degree{Linear}
        \entry{A+B+C}{}
            {Also denoted $\sigma$.  This polynomial is a subdivision of $A$.}
            {\begin{center}
                \begin{tabular}{ c p{1in} c }
                    Illustrated by & & with drawing \\
                     
    \newareastrue
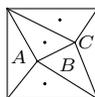
    \begin{tikzpicture}[scale=.2]
    \coordinate (P) at (0,0);
\coordinate (Q) at (6,0);
\coordinate (R) at (6,6);
\coordinate (S) at (0,6);

\draw[fill=black] (P) circle (1pt);
\draw[fill=black] (Q) circle (1pt);
\draw[fill=black] (R) circle (1pt);
\draw[fill=black] (S) circle (1pt);
\draw (P) -- (Q) -- (R) -- (S) -- cycle;

\coordinate (u) at (2,2.5);
\coordinate (v) at (4.5,3.75);
\draw[fill=black] (u) circle (1pt);
\draw[fill=black] (v) circle (1pt);

\draw (P) -- (u) -- (v) -- (R);
\draw (Q) -- (u) -- (S);
\draw (S) -- (v) -- (Q);

\ifcorners
    \draw (P) node[anchor=north east]{$\PP$};
    \draw (Q) node[anchor=north west]{$\QQ$};
    \draw (R) node[anchor=south west]{$\RR$};
    \draw (S) node[anchor=south east]{$\SS$};
\fi

\ifareas
    \draw (2.25,1) node{$A$};
    \draw (4,2.25) node{$B$};
    \draw (5.5,2.75) node{$C$};
    \draw (.75,2.75) node{$D$};
    \draw (2.25,4) node{$E$};
    \draw (4,5) node{$F$};
\fi

\ifinterior
    \draw (u) node[anchor=north]{$u$};
    \draw (v) node[anchor=west]{$v$};
\fi

\ifcolors
    \draw[line width=0, fill=\poscolor] (P) -- (Q) -- (u) -- cycle;
    \draw[line width=0, fill=\poscolor] (v) -- (Q) -- (R) -- cycle;
    \draw[line width=0, fill=\poscolor] (u) -- (v) -- (S) -- cycle;
    \draw[line width=0, fill=\negcolor] (u) -- (Q) -- (v) -- cycle;
    \draw[line width=0, fill=\negcolor] (v) -- (R) -- (S) -- cycle;
    \draw[line width=0, fill=\negcolor] (u) -- (S) -- (P) -- cycle;
\fi

\draw[fill=black] (2.5,1) circle (2pt);
\draw[fill=black] (3.5,5.25) circle (2pt);
\draw[fill=black] (2.5,3.75) circle (2pt);

\ifnewareas
    \draw (.75,2.75) node{\scriptsize $A$};
    \draw (4,2.25) node{\scriptsize $B$};
    \draw (5.25,3.75) node{\scriptsize $C$};
\fi

    \end{tikzpicture}
    \newareasfalse
     &  &  
    
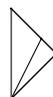
\begin{tikzpicture}[scale=.2]
    \coordinate (P) at (0,0);
\coordinate (Q) at (0,0);
\coordinate (R) at (0,6);
\coordinate (S) at (0,6);

\draw[fill=black] (P) circle (1pt);
\draw[fill=black] (Q) circle (1pt);
\draw[fill=black] (R) circle (1pt);
\draw[fill=black] (S) circle (1pt);
\draw (P) -- (Q) -- (R) -- (S) -- cycle;

\coordinate (x) at (3,3);
\draw[fill=black] (x) circle (1pt);

\coordinate (y) at (2,4);
\draw[fill=black] (y) circle (1pt);

\draw (P) -- (x);
\draw (Q) -- (x);
\draw (R) -- (x);
\draw (S) -- (x);

\draw (P) -- (y);

\ifcorners
    \draw (P) node[anchor=north]{$\PP=\QQ$};
%    \draw (Q) node[anchor=north west]{$\QQ$};
%    \draw (R) node[anchor=south west]{$\RR$};
    \draw (S) node[anchor=south]{$\RR=\SS$};
\fi

\ifcolors
    \draw[line width=0, fill=\poscolor] (P) -- (Q) -- (x) -- cycle;
    \draw[line width=0, fill=\poscolor] (R) -- (S) -- (x) -- cycle;
    \draw[line width=0, fill=\negcolor] (Q) -- (R) -- (x) -- cycle;
    \draw[line width=0, fill=\negcolor] (S) -- (P) -- (x) -- cycle;
\fi

\ifareas
    \draw (3,1.5) node{$A$};
    \draw (4.5,3) node{$B$};
    \draw (3,4.5) node{$C$};
    \draw (1.3,3) node{$D$};
\fi

\ifnewareas
    \draw (3,1.25) node{\scriptsize $D$};
    \draw (4.75,3) node{\scriptsize $B$};
    \draw (3,4.75) node{\scriptsize $C$};
    \draw (1.25,3) node{\scriptsize $A$};
\fi
    \end{tikzpicture}
    
                \end{tabular}
            \end{center}
            This algebraic subdivision can only be illustrated using a geometric subdivision.  The drawing here is obtained from the drawing of $A+B$ by subdividing one of the triangles.
            }
    
        \entry{A+B-C}{}
            {This polynomial is equivalent to a subdivision of $A-B$.}
            {\begin{center}
                \begin{tabular}{ c p{1in} c }
                    Illustrated by & & with drawing \\
                     
    \newareastrue
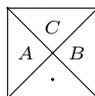
    \begin{tikzpicture}[scale=.2]
    \coordinate (P) at (0,0);
\coordinate (Q) at (6,0);
\coordinate (R) at (6,6);
\coordinate (S) at (0,6);

\draw[fill=black] (P) circle (1pt);
\draw[fill=black] (Q) circle (1pt);
\draw[fill=black] (R) circle (1pt);
\draw[fill=black] (S) circle (1pt);
\draw (P) -- (Q) -- (R) -- (S) -- cycle;

\coordinate (x) at (3,3);
\draw[fill=black] (x) circle (1pt);

\draw (P) -- (x);
\draw (Q) -- (x);
\draw (R) -- (x);
\draw (S) -- (x);

\draw[fill=black] (3,1.25) circle (2pt);

\ifnewareas
    \draw (1.2,3) node{\scriptsize $A$};
    \draw (4.6,3) node{\scriptsize $B$};
    \draw (3,4.75) node{\scriptsize $C$};
\fi
    \end{tikzpicture}
    \newareasfalse
     & &  
    
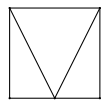
\begin{tikzpicture}[scale=.2]
    \coordinate (P) at (0,0);
\coordinate (Q) at (6,0);
\coordinate (R) at (6,6);
\coordinate (S) at (0,6);

\draw[fill=black] (P) circle (1pt);
\draw[fill=black] (Q) circle (1pt);
\draw[fill=black] (R) circle (1pt);
\draw[fill=black] (S) circle (1pt);
\draw (P) -- (Q) -- (R) -- (S) -- cycle;

\coordinate (x) at (3,0);
\draw[fill=black] (x) circle (1pt);

\draw (S) -- (x) -- (R);
    \end{tikzpicture}
    
                \end{tabular}
            \end{center}
            This algebraic subdivision can be illustrated without geometric subdivisions. 
            }
    
    \degree{Quadratic}
        \entry{A^2+B^2+C^2+2AB+2AC-2BC}{\ace}
            {This polynomial is in the abridged encyclopedia.}
            {\begin{center}
                \begin{tabular}{ c p{1in} c }
                    Illustrated by & & with drawing \\
                     
    \newareastrue
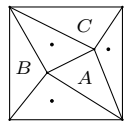
    \begin{tikzpicture}[scale=.25]
    \coordinate (P) at (0,0);
\coordinate (Q) at (6,0);
\coordinate (R) at (6,6);
\coordinate (S) at (0,6);

\draw[fill=black] (P) circle (1pt);
\draw[fill=black] (Q) circle (1pt);
\draw[fill=black] (R) circle (1pt);
\draw[fill=black] (S) circle (1pt);
\draw (P) -- (Q) -- (R) -- (S) -- cycle;

\coordinate (u) at (2,2.5);
\coordinate (v) at (4.5,3.75);
\draw[fill=black] (u) circle (1pt);
\draw[fill=black] (v) circle (1pt);

\draw (P) -- (u) -- (v) -- (R);
\draw (Q) -- (u) -- (S);
\draw (S) -- (v) -- (Q);

\ifcorners
    \draw (P) node[anchor=north east]{$\PP$};
    \draw (Q) node[anchor=north west]{$\QQ$};
    \draw (R) node[anchor=south west]{$\RR$};
    \draw (S) node[anchor=south east]{$\SS$};
\fi

\ifareas
    \draw (2.25,1) node{$A$};
    \draw (4,2.25) node{$B$};
    \draw (5.5,2.75) node{$C$};
    \draw (.75,2.75) node{$D$};
    \draw (2.25,4) node{$E$};
    \draw (4,5) node{$F$};
\fi

\ifinterior
    \draw (u) node[anchor=north]{$u$};
    \draw (v) node[anchor=west]{$v$};
\fi

\ifcolors
    \draw[line width=0, fill=\poscolor] (P) -- (Q) -- (u) -- cycle;
    \draw[line width=0, fill=\poscolor] (v) -- (Q) -- (R) -- cycle;
    \draw[line width=0, fill=\poscolor] (u) -- (v) -- (S) -- cycle;
    \draw[line width=0, fill=\negcolor] (u) -- (Q) -- (v) -- cycle;
    \draw[line width=0, fill=\negcolor] (v) -- (R) -- (S) -- cycle;
    \draw[line width=0, fill=\negcolor] (u) -- (S) -- (P) -- cycle;
\fi

\draw[fill=black] (2.25,1) circle (2pt);
\draw[fill=black] (5.25,3.75) circle (2pt);
\draw[fill=black] (2.25,4) circle (2pt);

\ifnewareas
    \draw (.75,2.75) node{\scriptsize $B$};
    \draw (4,2.25) node{\scriptsize $A$};
    \draw (4,5) node{\scriptsize $C$};
\fi

    \end{tikzpicture}
    \newareasfalse
     & &  
    
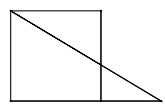
\begin{tikzpicture}[scale=.2]
    \coordinate (P) at (0,0);
\coordinate (Q) at (6,0);
\coordinate (R) at (6,6);
\coordinate (S) at (0,6);

\draw[fill=black] (P) circle (1pt);
\draw[fill=black] (Q) circle (1pt);
\draw[fill=black] (R) circle (1pt);
\draw[fill=black] (S) circle (1pt);
\draw (P) -- (Q) -- (R) -- (S) -- cycle;

\coordinate (u) at (10,0);
\coordinate (v) at (6,2.4);
% \coordinate (v) at (-4,6);
% \coordinate (u) at (0,3.6);
\draw[fill=black] (u) circle (1pt);
\draw[fill=black] (v) circle (1pt);

\draw (P) -- (u) -- (v) -- (R);
\draw (Q) -- (u) -- (S) -- (v) -- cycle;

\ifcorners
    \draw (P) node[anchor=north east]{$\PP$};
    \draw (Q) node[anchor=north]{$\QQ$};
    \draw (R) node[anchor=south west]{$\RR$};
    \draw (S) node[anchor=south east]{$\SS$};
\fi

\ifinterior
    \draw (u) node[anchor=north west]{$u$};
    \draw (v) node[anchor=south west]{$v$};
\fi
    \end{tikzpicture}
     %,\ \fig{.2}{D2ACE-alt}
                \end{tabular}
            \end{center}
            Other than $\sigma$, this is the first polynomial in the encyclopedia that is not illustrated by any drawing whose non-degenerate triangles all have positive real areas.
            }

\volume{4}
    \degree{Linear}
        \entry{A+B+C+D}{}
            {Also denoted $\sigma$.  This polynomial is a subdivision of $A$.}
            {\begin{center}
                \begin{tabular}{ c p{1in} c }
                    Illustrated by & & with drawing \\
                     
    \newareastrue
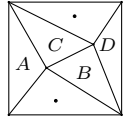
    \begin{tikzpicture}[scale=.25]
    \coordinate (P) at (0,0);
\coordinate (Q) at (6,0);
\coordinate (R) at (6,6);
\coordinate (S) at (0,6);

\draw[fill=black] (P) circle (1pt);
\draw[fill=black] (Q) circle (1pt);
\draw[fill=black] (R) circle (1pt);
\draw[fill=black] (S) circle (1pt);
\draw (P) -- (Q) -- (R) -- (S) -- cycle;

\coordinate (u) at (2,2.5);
\coordinate (v) at (4.5,3.75);
\draw[fill=black] (u) circle (1pt);
\draw[fill=black] (v) circle (1pt);

\draw (P) -- (u) -- (v) -- (R);
\draw (Q) -- (u) -- (S);
\draw (S) -- (v) -- (Q);

\ifcorners
    \draw (P) node[anchor=north east]{$\PP$};
    \draw (Q) node[anchor=north west]{$\QQ$};
    \draw (R) node[anchor=south west]{$\RR$};
    \draw (S) node[anchor=south east]{$\SS$};
\fi

\ifareas
    \draw (2.25,1) node{$A$};
    \draw (4,2.25) node{$B$};
    \draw (5.5,2.75) node{$C$};
    \draw (.75,2.75) node{$D$};
    \draw (2.25,4) node{$E$};
    \draw (4,5) node{$F$};
\fi

\ifinterior
    \draw (u) node[anchor=north]{$u$};
    \draw (v) node[anchor=west]{$v$};
\fi

\ifcolors
    \draw[line width=0, fill=\poscolor] (P) -- (Q) -- (u) -- cycle;
    \draw[line width=0, fill=\poscolor] (v) -- (Q) -- (R) -- cycle;
    \draw[line width=0, fill=\poscolor] (u) -- (v) -- (S) -- cycle;
    \draw[line width=0, fill=\negcolor] (u) -- (Q) -- (v) -- cycle;
    \draw[line width=0, fill=\negcolor] (v) -- (R) -- (S) -- cycle;
    \draw[line width=0, fill=\negcolor] (u) -- (S) -- (P) -- cycle;
\fi

\draw[fill=black] (2.5,.75) circle (2pt);
\draw[fill=black] (3.5,5.25) circle (2pt);

\ifnewareas
    \draw (.75,2.75) node{\scriptsize $A$};
    \draw (4,2.25) node{\scriptsize $B$};
    \draw (2.5,3.75) node{\scriptsize $C$};
    \draw (5.25,3.75) node{\scriptsize $D$};
\fi

    \end{tikzpicture}
    \newareasfalse
     &  &  
    
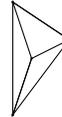
\begin{tikzpicture}[scale=.25]
    \coordinate (P) at (0,0);
\coordinate (Q) at (0,0);
\coordinate (R) at (0,6);
\coordinate (S) at (0,6);

\draw[fill=black] (P) circle (1pt);
\draw[fill=black] (Q) circle (1pt);
\draw[fill=black] (R) circle (1pt);
\draw[fill=black] (S) circle (1pt);
\draw (P) -- (Q) -- (R) -- (S) -- cycle;

\coordinate (u) at (1,3);
\coordinate (v) at (3,4);
\draw[fill=black] (u) circle (1pt);
\draw[fill=black] (v) circle (1pt);

\draw (P) -- (u) -- (v) -- (R);
\draw (Q) -- (u) -- (S) -- (v) -- cycle;
    \end{tikzpicture}
    
                \end{tabular}
            \end{center}
            This algebraic subdivision can only be illustrated using a geometric subdivision.  The drawing here is obtained from the drawing of $A+B$ by subdividing one of the triangles. The constrained triangulation shown is not toroidally irreducible, and it also illustrates the polynomial $A+B-C-D$.
            }
    
        \entry{A+B+C-D}{}
            {This polynomial is equivalent to a subdivision of $A-B$.}
            {\begin{center}
                \begin{tabular}{ c p{1in} c }
                    Illustrated by & & with drawing \\
                     
    \newareastrue
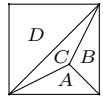
    \begin{tikzpicture}[scale=.2]
    \coordinate (P) at (0,0);
\coordinate (Q) at (6,0);
\coordinate (R) at (6,6);
\coordinate (S) at (0,6);

\draw[fill=black] (P) circle (1pt);
\draw[fill=black] (Q) circle (1pt);
\draw[fill=black] (R) circle (1pt);
\draw[fill=black] (S) circle (1pt);
\draw (P) -- (Q) -- (R) -- (S) -- cycle;

\draw (P) -- (R);

\coordinate (x) at (4,2);
\draw[fill=black] (x) circle (1pt);

\draw (x) -- (P);
\draw (x) -- (Q);
\draw (x) -- (R);

\ifcorners
    \draw (P) node[anchor=north east]{$\PP$};
    \draw (Q) node[anchor=north west]{$\QQ$};
    \draw (R) node[anchor=south west]{$\RR$};
    \draw (S) node[anchor=south east]{$\SS$};
\fi

\ifareas
    \draw (4,2) node{$A$};
    \draw (1.8,4) node{$B$};
\fi

\ifnewareas
    \draw (3.75,1) node{\scriptsize $A$};
    \draw (5.25,2.5) node{\scriptsize $B$};
    \draw (3.5,2.5) node{\scriptsize $C$};
    \draw (1.8,4) node{\scriptsize $D$};
\fi

\ifcolors
    \draw[line width=0, fill=\poscolor] (P) -- (Q) -- (R) -- cycle;
    \draw[line width=0, fill=\negcolor] (R) -- (S) -- (P) -- cycle;
\fi
    \end{tikzpicture}
    \newareasfalse
     &  &  
    
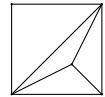
\begin{tikzpicture}[scale=.2]
    \coordinate (P) at (0,0);
\coordinate (Q) at (6,0);
\coordinate (R) at (6,6);
\coordinate (S) at (0,6);

\draw[fill=black] (P) circle (1pt);
\draw[fill=black] (Q) circle (1pt);
\draw[fill=black] (R) circle (1pt);
\draw[fill=black] (S) circle (1pt);
\draw (P) -- (Q) -- (R) -- (S) -- cycle;

\draw (P) -- (R);

\coordinate (x) at (4,2);
\draw[fill=black] (x) circle (1pt);

\draw (x) -- (P);
\draw (x) -- (Q);
\draw (x) -- (R);

\ifcorners
    \draw (P) node[anchor=north east]{$\PP$};
    \draw (Q) node[anchor=north west]{$\QQ$};
    \draw (R) node[anchor=south west]{$\RR$};
    \draw (S) node[anchor=south east]{$\SS$};
\fi

\ifareas
    \draw (4,2) node{$A$};
    \draw (1.8,4) node{$B$};
\fi

\ifnewareas
    \draw (3.75,1) node{\scriptsize $A$};
    \draw (5.25,2.5) node{\scriptsize $B$};
    \draw (3.5,2.5) node{\scriptsize $C$};
    \draw (1.8,4) node{\scriptsize $D$};
\fi

\ifcolors
    \draw[line width=0, fill=\poscolor] (P) -- (Q) -- (R) -- cycle;
    \draw[line width=0, fill=\negcolor] (R) -- (S) -- (P) -- cycle;
\fi
    \end{tikzpicture}
    
                \end{tabular}
            \end{center}
            This algebraic subdivision can only be illustrated using a geometric subdivision.  Here it is shown illustrated by an unconstrained triangulation with a subdivision.
            }
    
        \entry{A+B-C-D}{}
            {This polynomial is equivalent to a subdivision of $A-B$.}
            {\begin{center}
                \begin{tabular}{ c p{1in} c }
                    Illustrated by & & with drawing \\
                     
    \newareastrue
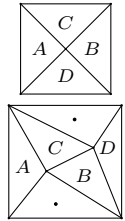
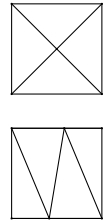
    \begin{tikzpicture}[scale=.2]
    \coordinate (P) at (0,0);
\coordinate (Q) at (6,0);
\coordinate (R) at (6,6);
\coordinate (S) at (0,6);

\draw[fill=black] (P) circle (1pt);
\draw[fill=black] (Q) circle (1pt);
\draw[fill=black] (R) circle (1pt);
\draw[fill=black] (S) circle (1pt);
\draw (P) -- (Q) -- (R) -- (S) -- cycle;

\coordinate (x) at (3,3);
\draw[fill=black] (x) circle (1pt);

\draw (P) -- (x);
\draw (Q) -- (x);
\draw (R) -- (x);
\draw (S) -- (x);

\ifcorners
    \draw (P) node[anchor=north east]{$\PP$};
    \draw (Q) node[anchor=north west]{$\QQ$};
    \draw (R) node[anchor=south west]{$\RR$};
    \draw (S) node[anchor=south east]{$\SS$};
\fi

\ifcolors
    \draw[line width=0, fill=\poscolor] (P) -- (Q) -- (x) -- cycle;
    \draw[line width=0, fill=\poscolor] (R) -- (S) -- (x) -- cycle;
    \draw[line width=0, fill=\negcolor] (Q) -- (R) -- (x) -- cycle;
    \draw[line width=0, fill=\negcolor] (S) -- (P) -- (x) -- cycle;
\fi

\ifareas
    \draw (3,1.25) node{$A$};
    \draw (4.75,3) node{$B$};
    \draw (3,4.75) node{$C$};
    \draw (1.25,3) node{$D$};
\fi

\ifnewareas
    \draw (3,1.25) node{\scriptsize $D$};
    \draw (4.75,3) node{\scriptsize $B$};
    \draw (3,4.75) node{\scriptsize $C$};
    \draw (1.25,3) node{\scriptsize $A$};
\fi
    \end{tikzpicture}
    \newareasfalse
     & &  
    \begin{tikzpicture}[scale=.2]
    \coordinate (P) at (0,0);
\coordinate (Q) at (6,0);
\coordinate (R) at (6,6);
\coordinate (S) at (0,6);

\draw[fill=black] (P) circle (1pt);
\draw[fill=black] (Q) circle (1pt);
\draw[fill=black] (R) circle (1pt);
\draw[fill=black] (S) circle (1pt);
\draw (P) -- (Q) -- (R) -- (S) -- cycle;

\coordinate (x) at (3,3);
\draw[fill=black] (x) circle (1pt);

\draw (P) -- (x);
\draw (Q) -- (x);
\draw (R) -- (x);
\draw (S) -- (x);

\ifcorners
    \draw (P) node[anchor=north east]{$\PP$};
    \draw (Q) node[anchor=north west]{$\QQ$};
    \draw (R) node[anchor=south west]{$\RR$};
    \draw (S) node[anchor=south east]{$\SS$};
\fi

\ifcolors
    \draw[line width=0, fill=\poscolor] (P) -- (Q) -- (x) -- cycle;
    \draw[line width=0, fill=\poscolor] (R) -- (S) -- (x) -- cycle;
    \draw[line width=0, fill=\negcolor] (Q) -- (R) -- (x) -- cycle;
    \draw[line width=0, fill=\negcolor] (S) -- (P) -- (x) -- cycle;
\fi

\ifareas
    \draw (3,1.25) node{$A$};
    \draw (4.75,3) node{$B$};
    \draw (3,4.75) node{$C$};
    \draw (1.25,3) node{$D$};
\fi

\ifnewareas
    \draw (3,1.25) node{\scriptsize $D$};
    \draw (4.75,3) node{\scriptsize $B$};
    \draw (3,4.75) node{\scriptsize $C$};
    \draw (1.25,3) node{\scriptsize $A$};
\fi
    \end{tikzpicture}
     \\
                     
    \newareastrue
    \begin{tikzpicture}[scale=.25]
    
    \end{tikzpicture}
    \newareasfalse
     & &  
    \begin{tikzpicture}[scale=.2]
    \coordinate (P) at (0,0);
\coordinate (Q) at (6,0);
\coordinate (R) at (6,6);
\coordinate (S) at (0,6);

\draw[fill=black] (P) circle (1pt);
\draw[fill=black] (Q) circle (1pt);
\draw[fill=black] (R) circle (1pt);
\draw[fill=black] (S) circle (1pt);
\draw (P) -- (Q) -- (R) -- (S) -- cycle;

\coordinate (u) at (2.5,0);
\coordinate (v) at (3.5,6);
\draw[fill=black] (u) circle (1pt);
\draw[fill=black] (v) circle (1pt);

\draw (P) -- (u) -- (v) -- (R);
\draw (Q) -- (u) -- (S) -- (v) -- cycle;
    \end{tikzpicture}
     
                \end{tabular}
            \end{center}
            This algebraic subdivision can be illustrated without geometric subdivisions.  The second constrained triangulation shown here is not toroidally irreducible, and it also illustrates the polynomial $\sigma$.
            }

    \degree{Quadratic}
        \entry{A^2+B^2+C^2-D^2+2AB+2AC-2BC}{\fourtwo}
            {This polynomial is in the abridged encyclopedia.}
            {\begin{center}
                \begin{tabular}{ c p{1in} c }
                    Illustrated by & & with drawing \\
                     
    \newareastrue
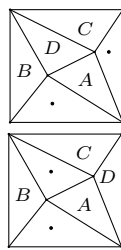
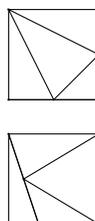
    \begin{tikzpicture}[scale=.25]
    \coordinate (P) at (0,0);
\coordinate (Q) at (6,0);
\coordinate (R) at (6,6);
\coordinate (S) at (0,6);

\draw[fill=black] (P) circle (1pt);
\draw[fill=black] (Q) circle (1pt);
\draw[fill=black] (R) circle (1pt);
\draw[fill=black] (S) circle (1pt);
\draw (P) -- (Q) -- (R) -- (S) -- cycle;

\coordinate (u) at (2,2.5);
\coordinate (v) at (4.5,3.75);
\draw[fill=black] (u) circle (1pt);
\draw[fill=black] (v) circle (1pt);

\draw (P) -- (u) -- (v) -- (R);
\draw (Q) -- (u) -- (S);
\draw (S) -- (v) -- (Q);

\ifcorners
    \draw (P) node[anchor=north east]{$\PP$};
    \draw (Q) node[anchor=north west]{$\QQ$};
    \draw (R) node[anchor=south west]{$\RR$};
    \draw (S) node[anchor=south east]{$\SS$};
\fi

\ifareas
    \draw (2.25,1) node{$A$};
    \draw (4,2.25) node{$B$};
    \draw (5.5,2.75) node{$C$};
    \draw (.75,2.75) node{$D$};
    \draw (2.25,4) node{$E$};
    \draw (4,5) node{$F$};
\fi

\ifinterior
    \draw (u) node[anchor=north]{$u$};
    \draw (v) node[anchor=west]{$v$};
\fi

\ifcolors
    \draw[line width=0, fill=\poscolor] (P) -- (Q) -- (u) -- cycle;
    \draw[line width=0, fill=\poscolor] (v) -- (Q) -- (R) -- cycle;
    \draw[line width=0, fill=\poscolor] (u) -- (v) -- (S) -- cycle;
    \draw[line width=0, fill=\negcolor] (u) -- (Q) -- (v) -- cycle;
    \draw[line width=0, fill=\negcolor] (v) -- (R) -- (S) -- cycle;
    \draw[line width=0, fill=\negcolor] (u) -- (S) -- (P) -- cycle;
\fi

\draw[fill=black] (2.25,1) circle (2pt);
\draw[fill=black] (5.25,3.75) circle (2pt);

\ifnewareas
    \draw (.75,2.75) node{\scriptsize $B$};
    \draw (4,2.25) node{\scriptsize $A$};
    \draw (4,5) node{\scriptsize $C$};
    \draw (2.25,4) node{\scriptsize $D$};
\fi

    \end{tikzpicture}
    \newareasfalse
     & &  
    \begin{tikzpicture}[scale=.2]
    \coordinate (P) at (0,0);
\coordinate (Q) at (6,0);
\coordinate (R) at (6,6);
\coordinate (S) at (0,6);

\draw[fill=black] (P) circle (1pt);
\draw[fill=black] (Q) circle (1pt);
\draw[fill=black] (R) circle (1pt);
\draw[fill=black] (S) circle (1pt);
\draw (P) -- (Q) -- (R) -- (S) -- cycle;

\coordinate (u) at (3,0);
\coordinate (v) at (6,3);
\draw[fill=black] (u) circle (1pt);
\draw[fill=black] (v) circle (1pt);

\draw (P) -- (u) -- (v) -- (R);
\draw (Q) -- (u) -- (S) -- (v) -- cycle;
    \end{tikzpicture}
     \\
                     
    \newareastrue
    \begin{tikzpicture}[scale=.25]
    \coordinate (P) at (0,0);
\coordinate (Q) at (6,0);
\coordinate (R) at (6,6);
\coordinate (S) at (0,6);

\draw[fill=black] (P) circle (1pt);
\draw[fill=black] (Q) circle (1pt);
\draw[fill=black] (R) circle (1pt);
\draw[fill=black] (S) circle (1pt);
\draw (P) -- (Q) -- (R) -- (S) -- cycle;

\coordinate (u) at (2,2.5);
\coordinate (v) at (4.5,3.75);
\draw[fill=black] (u) circle (1pt);
\draw[fill=black] (v) circle (1pt);

\draw (P) -- (u) -- (v) -- (R);
\draw (Q) -- (u) -- (S);
\draw (S) -- (v) -- (Q);

\ifcorners
    \draw (P) node[anchor=north east]{$\PP$};
    \draw (Q) node[anchor=north west]{$\QQ$};
    \draw (R) node[anchor=south west]{$\RR$};
    \draw (S) node[anchor=south east]{$\SS$};
\fi

\ifareas
    \draw (2.25,1) node{$A$};
    \draw (4,2.25) node{$B$};
    \draw (5.5,2.75) node{$C$};
    \draw (.75,2.75) node{$D$};
    \draw (2.25,4) node{$E$};
    \draw (4,5) node{$F$};
\fi

\ifinterior
    \draw (u) node[anchor=north]{$u$};
    \draw (v) node[anchor=west]{$v$};
\fi

\ifcolors
    \draw[line width=0, fill=\poscolor] (P) -- (Q) -- (u) -- cycle;
    \draw[line width=0, fill=\poscolor] (v) -- (Q) -- (R) -- cycle;
    \draw[line width=0, fill=\poscolor] (u) -- (v) -- (S) -- cycle;
    \draw[line width=0, fill=\negcolor] (u) -- (Q) -- (v) -- cycle;
    \draw[line width=0, fill=\negcolor] (v) -- (R) -- (S) -- cycle;
    \draw[line width=0, fill=\negcolor] (u) -- (S) -- (P) -- cycle;
\fi

\draw[fill=black] (2.25,1) circle (2pt);
\draw[fill=black] (2.25,4) circle (2pt);

\ifnewareas
    \draw (.75,2.75) node{\scriptsize $B$};
    \draw (4,2.25) node{\scriptsize $A$};
    \draw (4,5) node{\scriptsize $C$};
    \draw (5.25,3.75) node{\scriptsize $D$};
\fi

    \end{tikzpicture}
    \newareasfalse
     & &  
    \begin{tikzpicture}[scale=.2]
    \coordinate (P) at (0,0);
\coordinate (Q) at (6,0);
\coordinate (R) at (6,6);
\coordinate (S) at (0,6);

\draw[fill=black] (P) circle (1pt);
\draw[fill=black] (Q) circle (1pt);
\draw[fill=black] (R) circle (1pt);
\draw[fill=black] (S) circle (1pt);
\draw (P) -- (Q) -- (R) -- (S) -- cycle;

\coordinate (u) at (2,0);
\coordinate (v) at (1,3);
\draw[fill=black] (u) circle (1pt);
\draw[fill=black] (v) circle (1pt);

\draw (P) -- (u) -- (v) -- (R);
\draw (Q) -- (u) -- (S);
\draw (S) -- (v) -- (Q);
    \end{tikzpicture}
     \\
                \end{tabular}
            \end{center}
            }

        \entry{A^2+B^2+C^2+D^2-2AB+2AC+2AD-2BC+2BD-2CD}{\fourone}
            {This polynomial is in the abridged encyclopedia.}
            {\begin{center}
                \begin{tabular}{ c p{1in} c }
                    Illustrated by & & with drawing \\
                     
    \newareastrue
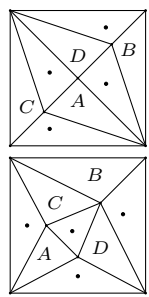
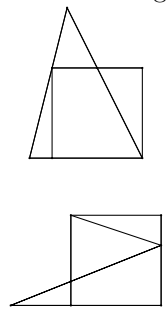
    \begin{tikzpicture}[scale=.3]
    \coordinate (P) at (0,0);
\coordinate (Q) at (6,0);
\coordinate (R) at (6,6);
\coordinate (S) at (0,6);

\draw[fill=black] (P) circle (1pt);
\draw[fill=black] (Q) circle (1pt);
\draw[fill=black] (R) circle (1pt);
\draw[fill=black] (S) circle (1pt);
\draw (P) -- (Q) -- (R) -- (S) -- cycle;

\coordinate (w) at (1.5,1.5);
\coordinate (x) at (3,3);
\coordinate (y) at (4.5,4.5);
\draw[fill=black] (w) circle (1pt);
\draw[fill=black] (x) circle (1pt);
\draw[fill=black] (y) circle (1pt);

\draw (P) -- (w) -- (x) -- (y) -- (R);
\draw (Q) -- (w) ;
\draw (Q) -- (x) ;
\draw (Q) -- (y) ;
\draw (S) -- (w) ;
\draw (S) -- (x) ;
\draw (S) -- (y) ;

\ifcorners
    \draw (P) node[anchor=north east]{$\PP$};
    \draw (Q) node[anchor=north west]{$\QQ$};
    \draw (R) node[anchor=south west]{$\RR$};
    \draw (S) node[anchor=south east]{$\SS$};
\fi

\ifcolors
    \draw[line width=0, fill=\poscolor] (P) -- (Q) -- (w) -- cycle;
    \draw[line width=0, fill=\poscolor] (x) -- (Q) -- (y) -- cycle;
    \draw[line width=0, fill=\poscolor] (R) -- (S) -- (y) -- cycle;
    \draw[line width=0, fill=\poscolor] (x) -- (S) -- (w) -- cycle;
    \draw[line width=0, fill=\negcolor] (w) -- (Q) -- (x) -- cycle;
    \draw[line width=0, fill=\negcolor] (y) -- (Q) -- (R) -- cycle;
    \draw[line width=0, fill=\negcolor] (y) -- (S) -- (x) -- cycle;
    \draw[line width=0, fill=\negcolor] (w) -- (S) -- (P) -- cycle;
\fi

\draw[fill=black] (1.75,.75) circle (2pt);
\draw[fill=black] (4.24,5.25) circle (2pt);
\draw[fill=black] (4.25,2.75) circle (2pt);
\draw[fill=black] (1.75,3.25) circle (2pt);

\ifnewareas
    \draw (3,2) node{\scriptsize $A$};
    \draw (5.25,4.25) node{\scriptsize $B$};
    \draw (.75,1.75) node{\scriptsize $C$};
    \draw (3,4) node{\scriptsize $D$};
\fi

    \end{tikzpicture}
    \newareasfalse
     & &  
    \begin{tikzpicture}[scale=.2]
    \coordinate (P) at (0,0);
\coordinate (Q) at (6,0);
\coordinate (R) at (6,6);
\coordinate (S) at (0,6);

\draw[fill=black] (P) circle (1pt);
\draw[fill=black] (Q) circle (1pt);
\draw[fill=black] (R) circle (1pt);
\draw[fill=black] (S) circle (1pt);
\draw (P) -- (Q) -- (R) -- (S) -- cycle;

\coordinate (w) at (-1.5,0);
\coordinate (x) at (1,10);
\coordinate (y) at (3,6);
\draw[fill=black] (w) circle (1pt);
\draw[fill=black] (x) circle (1pt);
\draw[fill=black] (y) circle (1pt);

\draw (P) -- (w) -- (x) -- (y) -- (R);
\draw (Q) -- (w) ;
\draw (Q) -- (x) ;
\draw (Q) -- (y) ;
\draw (S) -- (w) ;
\draw (S) -- (x) ;
\draw (S) -- (y) ;
    \end{tikzpicture}
     \\
                     
    \newareastrue
    \begin{tikzpicture}[scale=.3]
    %
% The quadratic graph for this is 
%
%   o   o
%   |   |
%   |   |
%   o---o
%
%
\coordinate (P) at (0,0);
\coordinate (Q) at (6,0);
\coordinate (R) at (6,6);
\coordinate (S) at (0,6);

\draw[fill=black] (P) circle (1pt);
\draw[fill=black] (Q) circle (1pt);
\draw[fill=black] (R) circle (1pt);
\draw[fill=black] (S) circle (1pt);
\draw (P) -- (Q) -- (R) -- (S) -- cycle;

\coordinate (w) at (3,1.6);
\coordinate (x) at (1.6,3);
\coordinate (y) at (4,4);
\draw[fill=black] (w) circle (1pt);
\draw[fill=black] (x) circle (1pt);
\draw[fill=black] (y) circle (1pt);

\draw (w) -- (y) -- (x) -- cycle;
\draw (Q) -- (w) -- (P);
\draw (Q) -- (y);
\draw (P) -- (x) -- (S);
\draw (S) -- (y) -- (R) ;

\ifcorners
    \draw (P) node[anchor=north east]{$\PP$};
    \draw (Q) node[anchor=north west]{$\QQ$};
    \draw (R) node[anchor=south west]{$\RR$};
    \draw (S) node[anchor=south east]{$\SS$};
\fi

\ifcolors
    \draw[line width=0, fill=\poscolor] (P) -- (Q) -- (w) -- cycle;
    \draw[line width=0, fill=\poscolor] (w) -- (x) -- (y) -- cycle;
    \draw[line width=0, fill=\poscolor] (x) -- (S) -- (P) -- cycle;
    \draw[line width=0, fill=\negcolor] (P) -- (w) -- (x) -- cycle;
    \draw[line width=0, fill=\negcolor] (w) -- (Q) -- (y) -- cycle;
    \draw[line width=0, fill=\negcolor] (y) -- (Q) -- (R) -- cycle;
    \draw[line width=0, fill=\negcolor] (R) -- (S) -- (y) -- cycle;
    \draw[line width=0, fill=\negcolor] (x) -- (S) -- (y) -- cycle;
\fi

\draw[fill=black] (3,.75) circle (2pt);
\draw[fill=black] (.75,3) circle (2pt);
\draw[fill=black] (2.75,2.75) circle (2pt);
\draw[fill=black] (5,3.5) circle (2pt);

\ifnewareas
    \draw (1.5,1.75) node{\scriptsize $A$};
    \draw (3.75,5.25) node{\scriptsize $B$};
    \draw (2,4) node{\scriptsize $C$};
    \draw (4,2) node{\scriptsize $D$};
\fi
    \end{tikzpicture}
    \newareasfalse
     & &  
    \begin{tikzpicture}[scale=.2]
    \coordinate (P) at (0,0);
\coordinate (Q) at (6,0);
\coordinate (R) at (6,6);
\coordinate (S) at (0,6);

\draw[fill=black] (P) circle (1pt);
\draw[fill=black] (Q) circle (1pt);
\draw[fill=black] (R) circle (1pt);
\draw[fill=black] (S) circle (1pt);
\draw (P) -- (Q) -- (R) -- (S) -- cycle;

\coordinate (w) at (-4,0);
\coordinate (x) at (0,1.6);
\coordinate (y) at (6,4);
\draw[fill=black] (w) circle (1pt);
\draw[fill=black] (x) circle (1pt);
\draw[fill=black] (y) circle (1pt);

\draw (w) -- (y) -- (x) -- cycle;
\draw (y) -- (Q) -- (w) -- (P) -- (x) -- (S) -- (y) -- (R) ;
    \end{tikzpicture}
     \\
                \end{tabular}
            \end{center}
            }

        \entry{A^2+B^2+C^2+D^2+2AB+2AC+2AD+2BC-2BD-2CD}{\aceothersub}
            {This polynomial is equivalent to the subdivision $\left[\ace\right](A,D,B+C)$.}
            {\begin{center}
                \begin{tabular}{ c p{1in} c }
                    Illustrated by & & with drawing \\
                     
    \newareastrue
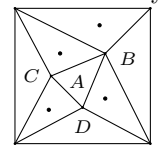
    \begin{tikzpicture}[scale=.3]
    %
% The quadratic graph for this is 
%
%   o   o
%   |
%   |
%   o---o
%
%
\coordinate (P) at (0,0);
\coordinate (Q) at (6,0);
\coordinate (R) at (6,6);
\coordinate (S) at (0,6);

\draw[fill=black] (P) circle (1pt);
\draw[fill=black] (Q) circle (1pt);
\draw[fill=black] (R) circle (1pt);
\draw[fill=black] (S) circle (1pt);
\draw (P) -- (Q) -- (R) -- (S) -- cycle;

\coordinate (w) at (3,1.6);
\coordinate (x) at (1.6,3);
\coordinate (y) at (4,4);
\draw[fill=black] (w) circle (1pt);
\draw[fill=black] (x) circle (1pt);
\draw[fill=black] (y) circle (1pt);

\draw (w) -- (y) -- (x) -- cycle;
\draw (Q) -- (w) -- (P);
\draw (Q) -- (y);
\draw (P) -- (x) -- (S);
\draw (S) -- (y) -- (R) ;

\ifcorners
    \draw (P) node[anchor=north east]{$\PP$};
    \draw (Q) node[anchor=north west]{$\QQ$};
    \draw (R) node[anchor=south west]{$\RR$};
    \draw (S) node[anchor=south east]{$\SS$};
\fi

\ifcolors
    \draw[line width=0, fill=\poscolor] (P) -- (Q) -- (w) -- cycle;
    \draw[line width=0, fill=\poscolor] (w) -- (x) -- (y) -- cycle;
    \draw[line width=0, fill=\poscolor] (x) -- (S) -- (P) -- cycle;
    \draw[line width=0, fill=\negcolor] (P) -- (w) -- (x) -- cycle;
    \draw[line width=0, fill=\negcolor] (w) -- (Q) -- (y) -- cycle;
    \draw[line width=0, fill=\negcolor] (y) -- (Q) -- (R) -- cycle;
    \draw[line width=0, fill=\negcolor] (R) -- (S) -- (y) -- cycle;
    \draw[line width=0, fill=\negcolor] (x) -- (S) -- (y) -- cycle;
\fi

\draw[fill=black] (2,4) circle (2pt);
\draw[fill=black] (4,2) circle (2pt);
\draw[fill=black] (1.5,1.5) circle (2pt);
\draw[fill=black] (3.75,5.25) circle (2pt);

\ifnewareas
    \draw (2.75,2.75) node{\scriptsize $A$};
    \draw (5,3.75) node{\scriptsize $B$};
    \draw (.75,3) node{\scriptsize $C$};
    \draw (3,.75) node{\scriptsize $D$};
\fi
    \end{tikzpicture}
    \newareasfalse
     & &  
    
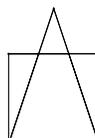
\begin{tikzpicture}[scale=.2]
    \coordinate (P) at (0,0);
\coordinate (Q) at (6,0);
\coordinate (R) at (6,6);
\coordinate (S) at (0,6);

\draw[fill=black] (P) circle (1pt);
\draw[fill=black] (Q) circle (1pt);
\draw[fill=black] (R) circle (1pt);
\draw[fill=black] (S) circle (1pt);
\draw (P) -- (Q) -- (R) -- (S) -- cycle;

\coordinate (w) at (3,9);
\coordinate (x) at (2,6);
\coordinate (y) at (4,6);
\draw[fill=black] (w) circle (1pt);
\draw[fill=black] (x) circle (1pt);
\draw[fill=black] (y) circle (1pt);

\draw (w) -- (y) -- (x) -- cycle;
\draw (y) -- (Q) -- (w) -- (P) -- (x) -- (S) -- (y) -- (R) ;
    \end{tikzpicture}
     \\
                \end{tabular}
            \end{center}
            This algebraic subdivision can be illustrated without geometric subdivisions. 
            The canonical 2-coloring of this $\T$ is not inherited from the unconstrained $T$.}

        \entry{A^2+B^2+C^2+D^2+2AB+2AC+2AD+2BC+2BD-2CD}{\acesub}
            {This polynomial is equivalent to the subdivision $\left[\ace\right](A+B,C,D)$.}
            {\begin{center}
                \begin{tabular}{ c p{1in} c }
                    Illustrated by & & with drawing \\
                     
    \newareastrue
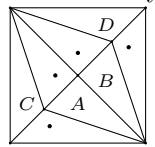
    \begin{tikzpicture}[scale=.3]
    \coordinate (P) at (0,0);
\coordinate (Q) at (6,0);
\coordinate (R) at (6,6);
\coordinate (S) at (0,6);

\draw[fill=black] (P) circle (1pt);
\draw[fill=black] (Q) circle (1pt);
\draw[fill=black] (R) circle (1pt);
\draw[fill=black] (S) circle (1pt);
\draw (P) -- (Q) -- (R) -- (S) -- cycle;

\coordinate (w) at (1.5,1.5);
\coordinate (x) at (3,3);
\coordinate (y) at (4.5,4.5);
\draw[fill=black] (w) circle (1pt);
\draw[fill=black] (x) circle (1pt);
\draw[fill=black] (y) circle (1pt);

\draw (P) -- (w) -- (x) -- (y) -- (R);
\draw (Q) -- (w) ;
\draw (Q) -- (x) ;
\draw (Q) -- (y) ;
\draw (S) -- (w) ;
\draw (S) -- (x) ;
\draw (S) -- (y) ;

\ifcorners
    \draw (P) node[anchor=north east]{$\PP$};
    \draw (Q) node[anchor=north west]{$\QQ$};
    \draw (R) node[anchor=south west]{$\RR$};
    \draw (S) node[anchor=south east]{$\SS$};
\fi

\ifcolors
    \draw[line width=0, fill=\poscolor] (P) -- (Q) -- (w) -- cycle;
    \draw[line width=0, fill=\poscolor] (x) -- (Q) -- (y) -- cycle;
    \draw[line width=0, fill=\poscolor] (R) -- (S) -- (y) -- cycle;
    \draw[line width=0, fill=\poscolor] (x) -- (S) -- (w) -- cycle;
    \draw[line width=0, fill=\negcolor] (w) -- (Q) -- (x) -- cycle;
    \draw[line width=0, fill=\negcolor] (y) -- (Q) -- (R) -- cycle;
    \draw[line width=0, fill=\negcolor] (y) -- (S) -- (x) -- cycle;
    \draw[line width=0, fill=\negcolor] (w) -- (S) -- (P) -- cycle;
\fi

\draw[fill=black] (1.75,.75) circle (2pt);
\draw[fill=black] (3,4) circle (2pt);
\draw[fill=black] (5.25,4.25) circle (2pt);
\draw[fill=black] (2,3) circle (2pt);

\ifnewareas
    \draw (3,1.75) node{\scriptsize $A$};
    \draw (4.25,2.75) node{\scriptsize $B$};
    \draw (.75,1.75) node{\scriptsize $C$};
    \draw (4.25,5.25) node{\scriptsize $D$};
\fi

    \end{tikzpicture}
    \newareasfalse
     & &  
    
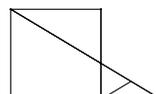
\begin{tikzpicture}[scale=.2]
    \coordinate (P) at (0,0);
\coordinate (Q) at (6,0);
\coordinate (R) at (6,6);
\coordinate (S) at (0,6);

\draw[fill=black] (P) circle (1pt);
\draw[fill=black] (Q) circle (1pt);
\draw[fill=black] (R) circle (1pt);
\draw[fill=black] (S) circle (1pt);
\draw (P) -- (Q) -- (R) -- (S) -- cycle;

\coordinate (u) at (10,0);
\coordinate (v) at (6,2.4);
% \coordinate (v) at (-4,6);
% \coordinate (u) at (0,3.6);
\draw[fill=black] (u) circle (1pt);
\draw[fill=black] (v) circle (1pt);

\draw (P) -- (u) -- (v) -- (R);
\draw (Q) -- (u) -- (S) -- (v) -- cycle;

\ifcorners
    \draw (P) node[anchor=north east]{$\PP$};
    \draw (Q) node[anchor=north]{$\QQ$};
    \draw (R) node[anchor=south west]{$\RR$};
    \draw (S) node[anchor=south east]{$\SS$};
\fi

\ifinterior
    \draw (u) node[anchor=north west]{$u$};
    \draw (v) node[anchor=south west]{$v$};
\fi

\coordinate (w) at (8,1.2);
\draw (w) -- (Q);
    \end{tikzpicture}
     \\
                \end{tabular}
            \end{center}
            This algebraic subdivision can only be illustrated using a geometric subdivision.
            The canonical 2-coloring of this $\T$ is not inherited from the unconstrained $T$, unless the subdivision is resolved first.
            }

    \degree{Quartic}
        %  This is a manual "entry"
        \noindent
        $\mathbf{A^4+4A^3B-4A^3C+4A^3D+6A^2B^2-4A^2BC+4A^2BD+6A^2C^2-4A^2CD }$ \\ 
        $\mathbf{+6A^2D^2+4AB^3+4AB^2C-4AB^2D-4ABC^2-40ABCD-4ABD^2-4AC^3 }$ \\
        $\mathbf{-4AC^2D+4ACD^2+4AD^3+B^4+4B^3C-4B^3D+6B^2C^2-4B^2CD+6B^2D^2 }$ \\
        $\mathbf{+4BC^3+4BC^2D-4BCD^2-4BD^3+C^4+4C^3D+6C^2D^2+4CD^3+D^4 }$
        
        \begin{addmargin}[.25in]{.25in}
            {This polynomial is in the abridged encyclopedia.}
            {\begin{center}
                \begin{tabular}{ c p{1in} c }
                    Illustrated by & & with drawing \\
                     
    \newareastrue
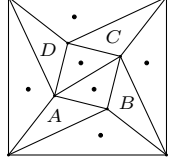
    \begin{tikzpicture}[scale=.35]
    \coordinate (P) at (0,0);
\coordinate (Q) at (6,0);
\coordinate (R) at (6,6);
\coordinate (S) at (0,6);

\draw[fill=black] (P) circle (1pt);
\draw[fill=black] (Q) circle (1pt);
\draw[fill=black] (R) circle (1pt);
\draw[fill=black] (S) circle (1pt);
\draw (P) -- (Q) -- (R) -- (S) -- cycle;

\coordinate (w) at (1.75,2.25);
\coordinate (x) at (3.75,1.75);
\coordinate (y) at (4.25,3.75);
\coordinate (z) at (2.25,4.25);
\draw[fill=black] (w) circle (1pt);
\draw[fill=black] (x) circle (1pt);
\draw[fill=black] (y) circle (1pt);
\draw[fill=black] (z) circle (1pt);

\draw (w) -- (x) -- (y) -- (z) -- cycle;
\draw (w) -- (y);

\draw (P) -- (w);
\draw (P) -- (x);
\draw (Q) -- (x);
\draw (Q) -- (y);
\draw (R) -- (y);
\draw (R) -- (z);
\draw (S) -- (z);
\draw (S) -- (w);

\ifcorners
    \draw (P) node[anchor=north east]{$\PP$};
    \draw (Q) node[anchor=north west]{$\QQ$};
    \draw (R) node[anchor=south west]{$\RR$};
    \draw (S) node[anchor=south east]{$\SS$};
\fi

\ifcolors
    \draw[line width=0,fill=\poscolor] (P) -- (Q) -- (x) -- cycle;
    \draw[line width=0,fill=\poscolor] (x) -- (Q) -- (y) -- cycle;
    \draw[line width=0,fill=\poscolor] (y) -- (R) -- (z) -- cycle;
    \draw[line width=0,fill=\poscolor] (P) -- (x) -- (w) -- cycle;
    \draw[line width=0,fill=\poscolor] (w) -- (x) -- (y) -- cycle;
    \draw[line width=0,fill=\poscolor] (w) -- (y) -- (z) -- cycle;
    \draw[line width=0,fill=\poscolor] (w) -- (z) -- (S) -- cycle;
    \draw[line width=0,fill=\poscolor] (z) -- (R) -- (S) -- cycle;
    \draw[line width=0,fill=\negcolor] (Q) -- (R) -- (y) -- cycle;
    \draw[line width=0,fill=\negcolor] (S) -- (P) -- (w) -- cycle;
\fi

\draw[fill=black] (2.5,5.25) circle (2pt);
\draw[fill=black] (3.5,.75) circle (2pt);
\draw[fill=black] (5.25,3.5) circle (2pt);
\draw[fill=black] (.75,2.5) circle (2pt);
\draw[fill=black] (3.25,2.5) circle (2pt);
\draw[fill=black] (2.75,3.5) circle (2pt);

\ifnewareas
    \draw (1.75,1.5) node{\scriptsize $A$};
    \draw (4.5,2) node{\scriptsize $B$};
    \draw (4,4.5) node{\scriptsize $C$};
    \draw (1.5,4) node{\scriptsize $D$};
\fi
    \end{tikzpicture}
    \newareasfalse
     & &  
    
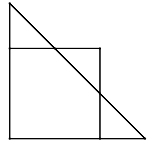
\begin{tikzpicture}[scale=.2]
    \coordinate (P) at (0,0);
\coordinate (Q) at (6,0);
\coordinate (R) at (6,6);
\coordinate (S) at (0,6);

\draw[fill=black] (P) circle (1pt);
\draw[fill=black] (Q) circle (1pt);
\draw[fill=black] (R) circle (1pt);
\draw[fill=black] (S) circle (1pt);
\draw (P) -- (Q) -- (R) -- (S) -- cycle;

\coordinate (w) at (0,9);
\coordinate (x) at (9,0);
\coordinate (y) at (6,3);
\coordinate (z) at (3,6);
\draw[fill=black] (w) circle (1pt);
\draw[fill=black] (x) circle (1pt);
\draw[fill=black] (y) circle (1pt);
\draw[fill=black] (z) circle (1pt);

\draw (w) -- (x) -- (y) -- (z) -- cycle;
\draw (w) -- (y);

\draw (P) -- (w);
\draw (P) -- (x);
\draw (Q) -- (x);
\draw (Q) -- (y);
\draw (R) -- (y);
\draw (R) -- (z);
\draw (S) -- (z);
\draw (S) -- (w);
    \end{tikzpicture}
     \\
                \end{tabular}
            \end{center}
            This is the first polynomial in the encyclopedia to fail the ``positivity'' conjecture made in \cite{chapter3}.  That conjecture applies only to area polynomials of unconstrained triangulations.
            }
        \end{addmargin} 
        \par
        \bigskip
        % end entry
        
\noindent
\hrule
\hrule`

%           \newpage
            \bigskip
            
    \section{The abridged encyclopedia}\label{sec:redencyc}
    
        Suppose that $T$ is a triangulation and we create a new triangulation $T'$ by subdividing some triangle $\Delta$ of $T$ into triangles $\{\Delta_j\}$.  Then it is easy to see that the polynomial $p_{T'}$ is obtained from the polynomial $p_T$ by replacing the variable $A_{\Delta}$ with the sum of the variables $A_{\Delta_j}$. Indeed, one can see this algebraic relation among some of the polynomials of $\encyc$ listed above.  For example, several of the quadratic polynomials listed in $\encyc_4$ can be obtained from polynomials in $\encyc_3$ by replacing one of the three variables $X$ by a sum $X'+X''$ of two variables. 
        Eliminating polynomials obtained in this way produces an {\it abridged encyclopedia} $\redencyc$.  It will turn out that we can recover $\encyc$ from $\redencyc$.
        
        \begin{defn}[Algebraic subdivision]
            Let $p$ be a polynomial in $n$ variables.  An {\em algebraic subdivision} of $p$ is a polynomial in $n+1$ variables obtained by substituting $X' + X''$ for $X$, where $X$ is any one of the variables of $p$.  
        \end{defn}
        
        \begin{defn}[Abridged encyclopedia]
            Let $l$ be a positive integer.   The {\em abridged encyclopedia} $\redencyc_l$ is the set of (equivalence classes of) polynomials $q\in \encyc_l$ such that $q$ is not an algebraic subdivision of any polynomial in $\encyc_{l-1}$.   We define $\redencyc=\cup_i \redencyc_i$.
        \end{defn}

        \begin{table}[h]%\label{table:vol1-4}
            \begin{tabular}{p{.7in}p{.7in}p{3.5in}}
                \toprule
                Volume & Quadratic code & Definition \\
                \midrule
                \addlinespace[5pt]
                     $\redencyc_1$ & - & $A$\\
                \addlinespace[5pt]
                \midrule
                \addlinespace[5pt]
                     $\redencyc_2$ & - & $A-B$\\
                \addlinespace[5pt]
                \midrule
                \addlinespace[5pt]
                     $\redencyc_3$ & \ace & $A^2+2AB+2AC+B^2-2BC+C^2$\\
                \addlinespace[5pt]
                \midrule
                \addlinespace[5pt]
                     $\redencyc_4$ & \fourone & $A^2-2AB+2AC+2AD+B^2-2BC+2BD+C^2-2CD+D^2$\\
                \addlinespace[5pt]
                     & \fourtwo & $A^2+2AB+2AC+B^2-2BC+C^2-D^2$\\
                \addlinespace[5pt]
                     & - & $A^4 + 4A^3B - 4A^3C + 4A^3D + 6A^2B^2 - 4A^2BC + 4A^2BD + 6A^2C^2 - 4A^2CD + 6A^2D^2 + 4AB^3 + 4AB^2C - 4AB^2D - 4ABC^2 - 40ABCD - 4ABD^2 - 4AC^3 - 4AC^2D + 4ACD^2 + 4AD^3 + B^4 + 4B^3C - 4B^3D + 6B^2C^2 - 4B^2CD + 6B^2D^2 + 4BC^3 + 4BC^2D - 4BCD^2 - 4BD^3 + C^4 + 4C^3D + 6C^2D^2 + 4CD^3 + D^4$ \\
                \addlinespace[5pt]
                \bottomrule
                \addlinespace[5pt]
            \end{tabular}
            \caption{The first four volumes of the abridged encyclopedia $\redencyc=\cup\redencyc_k$ of area relations.}
            \label{table:vol1-4}
        \end{table}
        
        Recall that in defining the area encyclopedia $\encyc$, we worked modulo the equivalence relation which allowed renaming of variables and multiplication by nonzero constants. This allowed us the slight abuse of calling the elements of $\encyc$ ``polynomials'' rather than equivalence classes.  This equivalence relation is easily seen to be compatible with algebraic subdivision, so by a similar abuse, we will speak of polynomials in the abridged encyclopedia.
        
        The abridged encyclopedia $\redencyc$ is defined by eliminating algebraic subdivisions from $\encyc$.  We now show that $\encyc$ is closed under both algebraic subdivision and its inverse, which implies that $\encyc$ can be recovered from $\redencyc$ exactly by taking the polynomials of $\redencyc$ and repeatedly subdividing.

        The key observation here is
        
        \begin{lem}
            Let $\qstar$ be an algebraic subdivision of $q$.  Then $\qstar\in \encyc$ if and only if $q\in\encyc$.
        \end{lem}
        
        \begin{proof}
            We may assume that $\qstar(A_0, A_1, \dots, A_{l+1})= q(A_0 + A_1,  A_2,\dots,  A_l)$. First assume that $q\in\encyc$.  Note that $q$ is irreducible, and hence $\qstar$ is irreducible.
            
            Since $q$ is in $\encyc_l$, there is a triangulation $T$ with $n$ triangles and an ordering of these triangles such that
            $q$ is an irreducible factor of 
            $r(A_1, A_2, \dots, A_l) = p_T(A_1, A_2,\dots, A_l, 0, \dots, 0)$.
            Now subdivide the triangle $\Delta_1$ of $T$ into three triangles $\Delta_1', \Delta_1'', \Delta_1'''$, and let $T'$ be the resulting triangulation.
            We have that the polynomial for $T'$ is 
            $$p_{T'}(A_1', A_1'', A_1''', A_2,\dots,A_l, \dots, A_n) = p_T(A_1' + A_1'' + A_1''', A_2,\dots,A_l,\dots, A_n).$$ 
            Substituting  $A_1''' = 0$ as well as $A_{l+1},\dots,  A_n = 0$, we see that $p_{T'}$ has an $(l+1)$-variable specialization $r'$ defined by
            \begin{align*}
                r' = p_{T'}(A_1' A_1'', 0, A_2,\dots,A_l, 0, \dots,0)    &= p_T(A_1' + A_1'', A_2,\dots,A_l, 0,\dots, 0)\\
                &= r( A_1' + A_1'', A_2, \dots,A_l).
            \end{align*} 
            Now $q$ divides $r$ implies $q*$ divides $r'$.  
            Hence $q*$ is in $\encyc$ by definition. 
            
            Now we prove the converse.  Assume $\qstar\in\encyc$.
            This implies $\qstar$, hence also $q$, is irreducible.  Since $\qstar\in\encyc$, there is a triangulation $T$ and a specialization $r(A_0, A_1, ..., A_l) = p_T(A_0, A_1, \dots, A_l, 0, \dots, 0)$ such that $\qstar$ is a factor of $r$.  Hence $q(A_0 + A_1,  A_2,\dots, A_{l-1}, A_l)$ is a factor of $r$, which is to say that there is a polynomial $u$ such that
            $$q(A_0 + A_1,  A_2,\dots, A_{l-1},  A_l) = r(A_0, A_1, \dots, A_l) u(A_0, A_1, \dots, A_l).$$
            Setting $A_0= 0$, we see that $q$ divides $r':=r(0, A_1, A_2,\dots, A_l)$.  Since $r$ is a specialization of $p_T$, so is $r'$.  Thus $q$ is an irreducible factor of a specialization of $p_T$.  Hence $q\in\encyc$.
        \end{proof}
        
        \begin{prop}
                 The encyclopedia $\encyc$ is the closure of the abridged encyclopedia $\redencyc$ under the operation of algebraic subdivision.  That is, a polynomial $q$ is in $\encyc$ if and only if $q$ can be obtained from some polynomial in $\redencyc$ by repeated algebraic subdivision. 
        \end{prop}
        
        \begin{proof}
            Suppose that $q$ can be obtained from some polynomial $q'\in \redencyc$ by repeated algebraic subdivision.  Then since $q'\in \encyc$, the lemma implies that $q \in \encyc$. 
            
            Conversely, take any polynomial $q$ in $\encyc_l$. Create a maximal sequence of polynomials $q_0=q, q_1, \dots, q_k$, where each $q_i$ is an algebraic subdivision of $q_{i+1}$. The lemma implies that $q_k\in \encyc$. By the maximality of the chain, $q_k$ is not an algebraic subdivision of any polynomial.  Hence $q_k\in \redencyc$, and our work is done.
        \end{proof}
    
        We point out that algebraic subdivision need not correspond to geometric subdivision. For example the dissection shown in Figure \ref{fig:introillustrate} has polynomial $A-B+C$ which is an algebraic subdivision of $A-B$.  However the dissection itself does not contain a subdivision.  Nevertheless, the proposition ensures that if we observe a subdivided polynomial in $\encyc$, then we may merge the variables and stay in $\encyc$.

    \section{Linears and quadratics}

        It is easy to see that $\encyc$ is closed under specialization; thus Theorem \ref{thm:encyc2} or Theorem \ref{thm:encyc} implies that the only linear polynomials in $\redencyc$ are $A$ and $A-B$.
        
        Similarly, our calculations impose considerable restrictions on the set of quadratic polynomials that could be in $\encyc$.  Suppose $q\in\encyc_l$.  We know from Theorem \ref{thm:encyc2} that all leading terms occur with the same coefficient up to sign.  If $q$ is quadratic then Theorem \ref{thm:encyc2} allows us to analyze the cross terms as well:  as $e+f=2$, the term $A_iA_j$ occurs either with coefficient $0$ (if $A_i^2$ and $A_j^2$ have opposite signs) or else with coefficient plus or minus twice the leading coefficient (if $A_i^2$ and $A_j^2$ have the same sign).  

        We can therefore encode any quadratic polynomial $q\in\encyc$ using a graph as follows:  there are vertices of two colors corresponding to the canonically 2-colored variables, and an edge connects two vertices of the same color if and only if the sign of the cross term differs from the sign of the leading terms.  The 6-variable quadratic polynomial for the triangulation $T_2$, for example, is encoded by the graph \six.  
        
        Any graph with 2-colored vertices and no edges joining vertices of different colors represents a quadratic polynomial that is compatible with Theorem \ref{thm:encyc2}.
        However, not every polynomial encoded by such a graph occurs in $\encyc$.  Some are reducible, but others simply don't appear.  For example, $\encyc_3$ does not include the irreducible quadratic corresponding to the graph \threeforbidden\ .
        
        Using this graphical encoding, specializing corresponds to deleting a subset of vertices, i.e., taking an induced subgraph.  (Algebraic subdivision is also easily described in terms of the graph.)
        So the fact that the polynomial \threeforbidden does not occur in $\encyc_3$ implies that no quadratic in $\encyc$ can have a triangle in its encoding graph.  Said differently, the irreducible polynomial $A^2+B^2+C^2-2AB-2AC-2BC$ is not a factor of any specialization of any area polynomial.

        \begin{question}
            Does $\redencyc$ contain finitely quadratic polynomials? Generally, does $\redencyc$ contain finitely many polynomials of each degree? 
        \end{question} 

        The first cubics appear in $\encyc_5$, for instance as specializations of $p_T$ for $T=T_3$ and $T=T_{2,1}$.

\part*{Appendix:  Nuggets on wheels}
    
    % \section{Nuggets on wheels}
        Imagine a wheel with $n$ spokes.  That is, take a graph consisting of a cycle with vertices $S_1, \dots S_n$ and a hub vertex $R$ connected to each $S_i$.  Map this graph to the plane, and suppose that the vertices are moved in a way such that the $S_i$ are brought together, all approaching the same point in the plane, and such that $R$ approaches a different point. We show that the area of the $n$-gon $\pentagon$ formed by the $S_i$ becomes small in magnitude relative to the sum of the absolute values of the areas of the $n$ triangles that comprise $\pentagon$.
        
            \begin{nug}  
                Let $n\geq 3$, and let  $S_1(s), \dots,  S_n(s), R(s)$ be continuous paths in $\C^2$ parameterized by $s\in (-\epsilon, \epsilon)$.  Let $\pentagon(s)$ denote the $n$-gon determined by $S_1(s), \dots,  S_n(s)$. Suppose that the $S_i(s)$ all approach the same point in $\C^2$ as $s \to 0$, and that $R(s)$ approaches a different point of $\C^2$.   Then 
                \begin{equation*}
                    \lim_{s \to 0} \frac{|\area(\pentagon(s))|}{\sum\limits_{i=1}^n |\area(S_i(s) S_{i+1}(s) R(s))|} = 0,
                \end{equation*}
                provided that the denominator does not vanish in a neighborhood of 0. 
            \end{nug}
            
            Note that here $\area(\pentagon(s))$ is defined to be the sum of the areas $S_i(s) S_{i+1}(s) P$ for an arbitrary point $P$.  This is a quadratic function of the coordinates of the $S_i(s)$ and is independent of $P$.  Also as usual subscripts are interpreted mod $n$.
            
            \begin{proof}
                We first prove the $n=3$ case and then proceed by induction on $n$. Without loss of generality, we may assume that the $S_i$ approach $(0,0)$ and $R$ approaches $(1,0)$. For each $s$, let $L_s$ be the line through $S_3(s)$ and $R(s)$.  Note that for $s$ sufficiently close to $0$ the line $L_s$ is not vertical.  For $i=1,2$ let $V_i(s)$ be the result of projecting $S_i(m)$ vertically onto $L_s$; precisely, $V_i(s)$ lies on $L_s$ and has the same $x$-coordinate as $S_i(s)$. 
                The area of triangle $S_1S_2S_3$ may be expressed formally as the sum of the areas of the five triangles
                $$ S_1 S_2 V_2, S_2 S_3 V_2, S_1 V_1 S_3, S_1 V_2 V_1, S_3 V_1 V_2 $$
                (Here and in what follows, we have omitted the parameter $s$ from the notation; we rely on the reader to keep in mind that the $S_i, V_i,$ and $R$ are all implicitly functions of $s$.)
                The last triangle has area 0, and the other four each have a segment which is vertical and a width which is vanishingly small. This makes it easy to evaluate their areas and show that each has vanishingly small area compared with one of the $SSR$ triangles, hence compared with the sum.  For example, for the first triangle above we have
                for any fixed $t>0$ that
                $$|\area(S_1 S_2 V_2)| = \frac12 |x(S_1) - x(S_2) |\cdot | y(S_2) - y(V_2) |.$$
                
                Comparing this to  
                $$|\area(S_2 S_3 R)| = \frac12 |x(R) - x(S_3) |\cdot | y(S_2) - y(V_2) |,$$
                we see that the ratio
                $$\frac{|\area(S_1 S_2 V_2)|}{|\area(S_2 S_3 R)|} = \frac{|x(S_1) - x(S_2) |}{|x(R) - x(S_3) |} \to 0 \quad{\mbox{ as }} t\to 0.$$
                Similar calculations hold for the other triangles, and the $n=3$ case is readily established. 
                
                For $n\geq 4$, we split the $n$-gon $S_1, \dots S_n$ into the $(n-1)$-gon $S_1, \dots S_{n-1}$ and the triangle $S_1 S_{n-1} S_n$.  It suffices to show that each of these areas when divided by the sum of the $SSR$ areas tends to $0$. For the triangle, we use the $n=3$ case to get 
                $$\frac{|\area(S_1S_{n-1}S_n)|}{|\area(S_1 R S_{n-1})| +|\area(S_{n-1} R S_n)|+|\area(S_n R S_1)|} \to 0.$$
                Since the four triangles involved in this fraction formally sum to zero and the area in the numerator is small compared with the others, it follows that
                $$\frac{|\area(S_1S_{n-1}S_n)|}{|\area(S_{n-1} R S_n)|+|\area(S_n R S_1)|} \to 0, $$
                and hence the area of triangle $S_1S_{n-1}S_n$ divided by the sum of the $SSR$ areas approaches $0$.  
                The analogous statement for the $(n-1)$-gon is similar.  By induction, we have at each $t>0$
                $$\frac{|\area(S_1 S_2 \cdots S_{n-1})|}{\sum_{i=1}^{n-2}|\area(S_i R  S_{i+1})| + |\area(S_{n-1}RS_1)|}\to 0.$$
                Again, we can use the fact that the areas of all these regions formally sum to $0$ together with the statement that the numerator is small compared with the others to assert that
                $$\frac{|\area(S_1, \dots S_{n-1})|}{\sum_{i=1}^{n-2}|\area(S_i R S_{i+1})|}\to 0.$$
                This completes the proof.
            \end{proof}

            We will make use of Nugget 1 in the following slightly stronger form.
            
            \begin{cor}\label{cor:nugget}
                With notation of the previous nugget, fix any subset $I\subset \{1, \dots, n\}$ with at least $3$ elements, and let $\Delta(s)$ denote the polygon determined by the points $S_i(s)$ where $i\in I$.   Then 
                \begin{equation*}
                    \lim_{s \to 0} \frac{|\area(\Delta(s))|}{\sum\limits_{i=1}^n |\area(S_i(s) S_{i+1}(s) R(s))|} = 0,
                \end{equation*}
                provided that the denominator does not vanish in a neighborhood of 0.  
            \end{cor}
             
            \begin{proof}
                It is enough to handle the case $I=\{1,\dots,k\}$ with $k<n$. By the nugget applied to $\Delta$, we have
                $$\frac{|\area(\Delta)|}{\sum_{i=1}^{k-1}|\area(S_i R  S_{i+1})| + |\area(S_{k}RS_1)|}\to 0.$$
                The $k$ triangles involved in this equation and $\Delta$ formally sum to $0$. It follows that 
                $$\frac{|\area(\Delta)|}{\sum_{i=1}^{k-1}|\area(S_i R  S_{i+1})| }\to 0,$$
                which yields the desired conclusion.
            \end{proof}

\addcontentsline{toc}{part}{References}
\bibliographystyle{plain}
\bibliography{References}

\end{document}